\documentclass[12pt,reqno]{amsart}
\usepackage{amsmath,amsfonts,amsthm,amsopn,amssymb}
\usepackage{cite,marginnote}
\usepackage{color,enumitem,graphicx}
\usepackage[colorlinks=true,urlcolor=blue,
citecolor=red,linkcolor=blue,linktocpage,pdfpagelabels,
bookmarksnumbered,bookmarksopen]{hyperref}
\usepackage[english]{babel}

\usepackage[left=2.9cm,right=2.9cm,top=2.8cm,bottom=2.8cm]{geometry}

\numberwithin{equation}{section}

\makeindex

\newtheorem{Remark}{Remark}[section]

\newtheorem{theo}{Theorem}[section]

\newtheorem{lemma}{Lemma}[section]
\newtheorem{iteration lemma}{iteration Lemma}[section]

\newcommand{\s}{\section}

\newcommand{\R}{\mathbb R}

\newcommand{\bt}{\begin{theorem}}
\newcommand{\et}{\end{theorem}}
\newcommand{\bl}{\begin{lemma}}
\newcommand{\el}{\end{lemma}}
\newcommand{\bd}{\begin{definition}}
\newcommand{\ed}{\end{definition}}
\newcommand{\bc}{\begin{corollary}}
\newcommand{\ec}{\end{corollary}}
\newcommand{\bp}{\begin{proof}}
\newcommand{\ep}{\end{proof}}
\newcommand{\bx}{\begin{example}}
\newcommand{\ex}{\end{example}}
\newcommand{\bi}{\begin{exercise}}
\newcommand{\ei}{\end{exercise}}
\newcommand{\bo}{\begin{proposition}}
\newcommand{\eo}{\end{proposition}}
\newcommand{\br}{\begin{remark}}
\newcommand{\er}{\end{remark}}
\newcommand{\beq}{\begin{equation}}
\newcommand{\eeq}{\end{equation}}
\newcommand{\ba}{\begin{align}}
\newcommand{\ea}{\end{align}}
\newcommand{\bn}{\begin{enumerate}}
\newcommand{\en}{\end{enumerate}}
\newcommand{\bg}{\begin{align*}}
\newcommand{\bcs}{\begin{cases}}
\newcommand{\ecs}{\end{cases}}

\newcommand{\bean}{\begin{eqnarray*}}
\newcommand{\eean}{\end{eqnarray*}}

\def\R{\mathbb{R}}

\def\bd{\mathrm{bd}\,}

\title[Single-peak and multi-peak solutions]{Single-peak and multi-peak solutions for Hamiltonian elliptic systems in dimension two}

\author[H. Zhang]{Hui Zhang}
\author[M. B.\ Yang]{Minbo Yang}
\author[J. J. Zhang]{Jianjun Zhang}
\author[X.~X.~Zhong]{Xuexiu Zhong}

\address[H.\ Zhang]{\newline\indent Department of  Mathematics,
Jinling Institute of Technology,
\newline\indent
Nanjing 211169, PR China
\newline\indent and
\newline\indent Department of Mathematics,
Nanjing University,
\newline\indent
Nanjing 210093, PR China}
\email{\href{mailto:huihz0517@126.com}{huihz0517@126.com}}

\address[M. B.\ Yang]{\newline\indent Department of Mathematics
\newline\indent
Zhejiang Normal University
\newline\indent
Jinhua 321004, PR China}
\email{\href{mailto:mbyang@zjnu.edu.cn}{mbyang@zjnu.edu.cn}}

\address[J. J. \ Zhang]{\newline\indent College of Mathematics and Statistics, Chongqing Jiaotong University,
\newline\indent
Chongqing 400074, PR China}
\email{\href{mailto:zhangjianjun09@tsinghua.org.cn}{zhangjianjun09@tsinghua.org.cn}}

\address[X.~X.~Zhong]{\newline\indent South China Research Center for Applied Mathematics and Interdisciplinary Studies
\newline\indent
South China Normal University
\newline\indent
Guangzhou 510631, PR China}
\email{\href{mailto:zhongxuexiu1989@163.com}{zhongxuexiu1989@163.com}}

\thanks{(1) Corresponding author: Minbo Yang ({\tt mbyang@zjnu.edu.cn})}

\thanks{(2) Hui Zhang was supported by China Postdoctoral Science Foundation (No.2021M691527). Minbo Yang was supported by NSFC(No.11971436, No.12011530199). Xuexiu Zhong was supported by the NSFC (No.11801581), Guangdong Basic and Applied Basic Research Foundation (2021A1515010034),Guangzhou Basic and Applied Basic Research Foundation(No.202102020225). Jianjun Zhang was supported by NSFC (No.12371109,11871123).}


\begin{document}

\begin{abstract}
 This paper is concerned with the Hamiltonian elliptic system in dimension two\begin{equation*}\aligned
\left\{ \begin{array}{lll}
-\epsilon^2\Delta u+V(x)u=g(v)\ & \text{in}\quad \mathbb{R}^2,\\
-\epsilon^2\Delta v+V(x)v=f(u)\ & \text{in}\quad \mathbb{R}^2,
\end{array}\right.\endaligned
\end{equation*}
where $V\in C(\mathbb{R}^2)$ has local minimum points, and $f,g\in C^1(\mathbb{R})$ are assumed to be of exponential growth in the sense of Trudinger-Moser inequality. When $V$ admits one or several local strict minimum points, we show the existence and concentration of single-peak and multi-peak semiclassical states respectively, as well as strong convergence and exponential decay. In addition, positivity of solutions and uniqueness of local maximum points of solutions are also studied. Our theorems extend the results of Ramos and Tavares [Calc. Var. 31 (2008) 1-25], where $f$ and $g$ have polynomial growth. It seems that it is the first attempt to obtain multi-peak semiclassical states for Hamiltonian elliptic system with exponential growth.
\end{abstract}
\maketitle

\s{Introduction and main results}
\renewcommand{\theequation}{1.\arabic{equation}}
Consider the time-dependent system of coupled Schr\"{o}dinger equations
\begin{equation}\label{1.4}\aligned
\left\{ \begin{array}{lll}
i\hbar\frac{\partial \varphi}{\partial t}=-\frac{\hbar^2}{2m}\Delta\varphi+W(x)\varphi-\frac{\partial H(\varphi,\psi)}{\partial \psi}, \ & x\in\mathbb{R}^N,\ t>0,\\
i\hbar\frac{\partial \psi}{\partial t}=-\frac{\hbar^2}{2m}\Delta\psi+W(x)\psi-\frac{\partial H(\varphi,\psi)}{\partial \varphi}, \ & x\in\mathbb{R}^N,\ t>0,
\end{array}\right.\endaligned
\end{equation}
where $N\geq2$, $(\varphi,\psi)$ represents the wave function of the state of an electron, $i$ is the
imaginary unit, $m$ is the mass of a particle, $\hbar$ is the Planck constant, $W(x)$ is a continuous potential, and $H$ is a coupled nonlinear function modeling various types of interaction
effect among many particles. The problem (\ref{1.4}) arises in many fields of physics, especially in nonlinear optics and
Bose-Einstein condensates theory. One of the most interesting problem for (\ref{1.4}) is to look for standing waves, i.e.
$$(\varphi(x,t),\psi(x,t))=\bigl(u(x)e^{-iE/\hbar t}, v(x)e^{-iE/\hbar t}\bigr).$$
Assume $\frac{\partial H(\varphi,\psi)}{\partial \psi}=e^{-iE/\hbar t}\frac{\partial H(u,v)}{\partial v}$, $\frac{\partial H(\varphi,\psi)}{\partial \varphi}=e^{-iE/\hbar t}\frac{\partial H(u,v)}{\partial u},$
then $(\varphi,\psi)$ satisfies (\ref{1.4}) if and only if $(u,v)$ is a solution of the system
\begin{equation}\label{1.3}\aligned
\left\{ \begin{array}{lll}
-\epsilon^2\Delta u+V(x)u=\frac{\partial H(u,v)}{\partial v}, \ & x\in\mathbb{R}^N,\\
-\epsilon^2\Delta v+V(x)v=\frac{\partial H(u,v)}{\partial u}, \ & x\in \mathbb{R}^N,
\end{array}\right.\endaligned
\end{equation}
where $\epsilon^2=\frac{\hbar^2}{2m}$ and $V(x)=W(x)-E$. The existence, multiplicity and concentration of solutions of (\ref{1.3}) are widely investigated and the scenario changes remarkably from the higher
dimensional case $N\geq3$ to the planar case $N=2$. In particular, $N=2$ affects the notion of critical growth which is the maximal admissible growth for the nonlinearities to preserve the variational structure of the problem, for details see \cite{Cassani1,ZhangCPDE} for coupled Schr\"{o}dinger systems and see \cite{TaoLiuYang,AlvesYang,ZhangDO} for a single Schr\"{o}dinger equation for instance. In the case that $\epsilon=1$ in (\ref{1.3}), the results on the existence and properties of solutions can be seen in \cite{BD,DY,LY,BDR,BS,DING} for $N\geq3$ and \cite{DG1,Cassani1,Cassani2,BDRT} for $N=2$ for example.

In the present paper, we are interested in the case $\epsilon\rightarrow0$, which is called semiclassical problem and  the associated solutions are called semiclassical states, which possesses an important physical interest in
describing the translation from quantum to classical mechanics.
The semiclassical states of (\ref{1.3}) has been intensively studied and most of the results are focused on the case $N\geq3$.
 Alves et al. \cite{alves} dealt with the system in the whole space
\begin{equation}\label{1.10}\aligned
\left\{ \begin{array}{lll}
-\epsilon^2\Delta u+u=Q(x)|v|^{p-1}v\ & \text{in}\quad \mathbb{R}^N,\\
-\epsilon^2\Delta v+v=K(x)|u|^{q-1}u\ & \text{in}\quad \mathbb{R}^N,
\end{array}\right.\endaligned
\end{equation}
where $2<p,q<2N/(N-2)$, $Q$ and $K$ are positive bounded functions, they proved that (\ref{1.10}) admits a family of solutions concentrating at a point $x_0\in\mathbb{R}^N$ where
a related functional realizes its minimum energy.
Subsequently, Ramos and  Soares \cite{Ramos1} treated the system
\begin{equation}\label{1.5}\aligned
\left\{ \begin{array}{lll}
-\epsilon^2\Delta u+V(x)u=g(v)\ & \text{in}\quad \mathbb{R}^N,\\
-\epsilon^2\Delta v+V(x)v=f(u)\ & \text{in}\quad \mathbb{R}^N,
\end{array}\right.\endaligned
\end{equation}
where $N\geq3$, $V\in C(\mathbb{R}^N)$ satisfies
\begin{equation} \label{1.6} 0<\min_{x\in\mathbb{R}^N}V(x)<\liminf_{|x|\rightarrow\infty}V(x)\leq+\infty,\end{equation}
and $f,g\in C^1(\mathbb{R})$ are polynomial-type nonlinearities having superlinear and subcritical growth at infinity. They showed the existence of ground states of (\ref{1.5}) which concentrate at global minimum points of $V$ as $\epsilon\rightarrow0$.
Later,  Ramos and Tavares \cite{Ramos2} considered (\ref{1.5}) in an open domain $\Omega$ of $\mathbb{R}^N$, and replaced the global condition (\ref{1.6}) by assuming $V(x)\geq a>0$ is locally H\"{o}lder continuous, and
\begin{equation*}\aligned \inf_{\Lambda_i}V<\inf_{\partial\Lambda_i}V,\ \text{where}\ \Lambda_i\subset\Omega  \text{ are bounded and mutually disjoint}, i=1,\cdots,k.\endaligned\end{equation*}
By means of a reduction method
and a penalization technique, the authors in \cite{Ramos2} showed that (\ref{1.5}) possesses positive solutions $u_\epsilon, v_\epsilon$ and both of them have $k$ local maximum points $x_{i,\epsilon}\in \Lambda_i$, $i=1,\cdots,k$, as $\epsilon\rightarrow0$.

On the contrary, the related results of semiclassical states for (\ref{1.3}) in the plane are few.
In $\mathbb{R}^2$, the natural growth restriction on the nonlinearity is
defined by

\noindent{\bf Definition 1.1.} A function $h:\mathbb{R}\rightarrow\mathbb{R}$ is of subcritical exponential growth  if $$\lim_{|t|\rightarrow+\infty}\frac{|h(t)|}{e^{\alpha t^2}}=0, \quad \forall\alpha>0,$$ and critical exponential growth if there exists $\alpha_0>0$
such that
 $$
\lim_{|t|\rightarrow+\infty}\frac{|h(t)|}{e^{\alpha t^2}}
=\left\{\begin{array}{lll}
0\ & \text{if}\  \alpha>\alpha_0,\\
+\infty\ & \text{if}\ \alpha<\alpha_0.
\end{array}\right.
$$
Moreover, Definition 1.1 is derived from Trudinger-Moser inequality as follows.
\begin{lemma}\label{l1.2} (\cite{Cao}) If $\alpha>0$ and $u\in H^1(\mathbb{R}^2)$, then
$\int_{\mathbb{R}^2}\bigl(e^{\alpha u^2}-1\bigl)dx<+\infty.$
Moreover, if $|\nabla u|^2_2\leq1$, $|u|_2\leq M<+\infty$, and $\alpha<\alpha_0=4\pi$, then there exists a constant $C$, which depends only on $M$ and $\alpha$, such that
\begin{equation*}\int_{\mathbb{R}^2}\bigl(e^{\alpha u^2}-1\bigl)dx<C(M,\alpha).\end{equation*}
\end{lemma}

In \cite{ZhangCPDE}, Cassani and Zhang dealt with the system in dimension two
\begin{equation}\label{1.1}\aligned
\left\{ \begin{array}{lll}
-\epsilon^2\Delta u+V(x)u=g(v)\ & \text{in}\quad \mathbb{R}^2,\\
-\epsilon^2\Delta v+V(x)v=f(u)\ & \text{in}\quad \mathbb{R}^2,
\end{array}\right.\endaligned
\end{equation}
where $f,g\in C(\mathbb{R})$ are superlinear at infinity and of critical exponential growth. Under the global condition (\ref{1.6}) they showed the
existence of ground states concentrating around
global minimum points of the potential $V$ by the method of generalized Nehari manifold. In addition, Cassani and Zhang \cite{ZhangCPDE} studied the positivity and decay of ground states by virtue of a priori estimates of solutions for the limit system.

Hence, compared with related results of system (\ref{1.5}) with polynomial growth, it is natural to ask that:

\vskip 0.05true cm
{\it (Q$_1$) Whether system (\ref{1.5}) with exponential growth  has single-peak solutions when $V$ has one local strict minimum point?}

{\it (Q$_2$) Whether system (\ref{1.5}) with exponential growth  has multi-peak solutions when $V$ has several local strict minimum points?}

Observe that, there is few results of (\ref{1.5}) when $f,g$ are asymptotically linear. So another interesting problem is:

{\it (Q$_3$) What happens to system (\ref{1.5}) when $f, g$ are asymptotically linear?}

To answer question (Q$_2$), we shall introduce a technical condition, see (H$_7$) below, which is given by Ramos and Tavares\cite{Ramos2}. In addition,  in \cite{Ramos2}, (H$_7$) is removed by suitable truncated functions and the uniformly boundedness of solutions in $L^\infty$. So one  open problem is:

 {\it (Q$_4$) Whether the condition (H$_7$) below can be removed in looking for multi-peak solutions of system (\ref{1.5}) with exponential growth?}

In this paper, we are interested in these questions and give some answers
to them. The assumptions on $V$, $f$ and $g$ are given as follows.
\vskip 0.1 true cm
\noindent(V$_1$) $V\in C(\mathbb{R}^2)\cap L^\infty(\mathbb{R}^2)$ and  $\inf_{\mathbb{R}^2}V>0$.
\vskip 0.1 true cm
\noindent(V$_2$) There exists an open bounded set $\Lambda_0\subset\mathbb{R}^2$ with smooth boundary  such that
$V_0:=\inf_{\Lambda_0}V<\inf_{\partial\Lambda_0}V.$
\vskip 0.05 true cm
\noindent(H$_1$) $f,g\in C^1(\mathbb{R},\mathbb{R})$, $f(t)=o(t)$ and $g(t)=o(t)$ as $t\rightarrow0$.\\
\noindent(H$_2$) $t^2f'(t)>f(t)t\ \text{and}\  t^2g'(t)>g(t)t,$ $\forall t\neq0$.
\vskip 0.1 true cm
\noindent(H$_3$) (i) $g'$ and $f'$ are of subcritical exponential growth in the sense of Definition 1.1;\\
\indent\ \ (ii) there is $\theta>2$ such that $\theta F(t)\leq f(t)t$ and $\theta G(t)\leq g(t)t$ for any $t\neq0$.

We are also concerned with the asymptotically linear case and replace (H$_3$) by the following condition.
\vskip 0.1 true cm
\noindent(H$'_3$) (i)\ $|V|_{L^\infty(\mathbb{R}^2)}<\lim_{|t|\rightarrow+\infty}
\frac{f(t)}{t}=\lim_{|t|\rightarrow+\infty}
\frac{g(t)}{t}=l_0<+\infty$;\\
\indent\ \ (ii) $f(t)t-2F(t)\rightarrow+\infty$ and $g(t)t-2G(t)\rightarrow+\infty$ as $|t|\rightarrow+\infty$.

Firstly, we give the  results about single-peak solutions.
\begin{theo}\label{t1.1} {\it Suppose that (V$_1$), (V$_2$), (H$_1$), (H$_2$) and either (H$_3$) or (H$'_3$) are satisfied. Then for sufficiently small $\epsilon>0$,\\
\noindent(1) system (\ref{1.1})
has a nontrivial solution $(\varphi_\epsilon,\psi_\epsilon)$ in $H^1(\mathbb{R}^2)\times H^1(\mathbb{R}^2)$;\\
\noindent(2) If additionally $V$ is locally H\"{o}lder continuous, then \\
\noindent(i) $|\varphi_\epsilon|$, $|\psi_\epsilon|$ and $|\varphi_\epsilon|+|\psi_\epsilon|$ possess global maximum points $x^1_\epsilon, x^2_\epsilon, x_\epsilon\in\Lambda_0$ respectively, such that
 $$V(x_\epsilon)\rightarrow V_0, \ \text{and}\ V(x^i_\epsilon)\rightarrow V_0, \ i=1,2, \ \text{as}\ \epsilon\rightarrow0.$$}
{\it\noindent(ii)  $\bigl(\varphi_\epsilon(\epsilon x+x_\epsilon),\psi_\epsilon(\epsilon x+x_\epsilon)\bigr)$ and $\bigl(\varphi_\epsilon(\epsilon x+x^i_\epsilon),\psi_\epsilon(\epsilon x+x^i_\epsilon)\bigr),\ i=1,2$, converge  in $H^1(\mathbb{R}^2)\times H^1(\mathbb{R}^2)$ to ground states of
\begin{equation}\label{1.1.1}\aligned
\left\{ \begin{array}{lll}
-\Delta u+V_0u=g(v)\ & \text{in}\quad \mathbb{R}^2,\\
-\Delta v+V_0v=f(u)\ & \text{in}\quad \mathbb{R}^2.
\end{array}\right.\endaligned
\end{equation}
(iii)
$|\varphi_\epsilon(x)|+|\psi_\epsilon(x)|\leq C e^{-\frac c\epsilon|x-x_{\epsilon}|}$ for some $C,c>0$.}
\end{theo}
\begin{Remark}\label{r1.2} We shall use some ideas in \cite{DINGXU1,DINGXU2} to show Theorem \ref{t1.1}, and some remarks are as follows.

(1) We would like to point out that, even if $f,g$ have polynomial growth as in \cite{DINGXU1,DINGXU2}, our system (\ref{1.1}) is different from the Reaction-diffusion system or Dirac equation considered in \cite{DINGXU1,DINGXU2}, whose functional is
$$\Phi(z)=\frac12\|z^+\|^2-\frac12\|z^-\|^2+\frac12\int_{\mathbb{R}^N}V(x)z^2
-\int_{\mathbb{R}^N}F(z),$$
different from (\ref{2.2}) below. Note that in general the term $\int_{\mathbb{R}^N}uv$ in (\ref{2.2}) is not one-sign for any $u,v\in H^1(\mathbb{R}^2)$. This results in some new difficulties, for example, the comparison of energy and the monotonicity of least energy cannot be obtained in a standard way. Hence, we cannot argue trivially as \cite{DINGXU1,DINGXU2} even if $f,g$ have polynomial growth as in \cite{DINGXU1,DINGXU2}.

(2)  The proof of the results in \cite{DINGXU1} highly relies on the following inequality (also for $f$)\begin{equation*}\label{bu8}|g(u)|^{\frac p{p-1}}\leq Cg(u)u\ \text{for some}\ p\in(2,2N/(N-2)),\ \text{if} \ |u|\geq1,\end{equation*}
which is invalid since here $f,g$ are of exponential growth. This difficulty will be overcome by some new tricks and an important inequality (\ref{bu1}) below.
\end{Remark}

In order to investigate the positivity of solutions and the uniqueness of maximum points of
 solutions to (\ref{1.1}), we introduce the following conditions.
\vskip 0.05 true cm
\noindent(H$_4$) There exist $p,q>1$ such that $f(t)\geq t^q$ and $g(t)\geq t^p$ for small $t>0$.
\vskip 0.05 true cm
\noindent(H$_5$) $f(t)=g(t)=0$ if $t\leq0$.
\vskip 0.05 true cm
\begin{theo}\label{t1.2} {\it Assume (V$_1$), (V$_2$) and (H$_1$)-(H$_3$) are satisfied and $V$ is locally H\"{o}lder continuous. \\
\noindent (1) If in addition (H$_4$) holds and $f,g$ are odd, then $\varphi_\epsilon,\ \psi_\epsilon$, $x^1_\epsilon$, $x^2_\epsilon$ obtained in Theorem \ref{t1.1} satisfy\\
\noindent(i) $\varphi_\epsilon>0,\ \psi_\epsilon>0$ in $\mathbb{R}^2$;\\
\noindent(ii) $x^1_\epsilon$, $x^2_\epsilon$ are unique, and
$\lim_{\epsilon\rightarrow0}|x^1_\epsilon-x^2_\epsilon|/\epsilon=0.$

\noindent (2) If in addition (H$_5$) holds, then the results of the conclusion (1) remain true.}\end{theo}

\begin{theo}\label{t1.3} {\it Assume (V$_1$), (V$_2$), (H$_1$), (H$_2$) and (H$'_3$) are satisfied and $V$ is locally H\"{o}lder continuous. If in addition (H$_5$) holds, then the results of Theorem \ref{t1.2}-(1) remain true.}\end{theo}

\begin{Remark}\label{r1.3} On one hand, Theorems \ref{t1.1}  and \ref{t1.3} about asymptotically linear nonlinearity are also new for (\ref{1.1}) with polynomial growth in the higher dimensional case $N\geq3$. On the other hand, we believe that the results of \cite{ZhangCPDE} about superlinear nonlinearities can be extended to the asymptotically linear case using similar arguments.
\end{Remark}

Below we state the results of multi-peak solutions. We assume
\vskip 0.05 true cm
\noindent(V$_3$) There exist bounded domains $\Lambda_i$, mutually disjoint, compactly contained in $\mathbb{R}^N$, $i=1,...,k$, such that $V_i:=\inf_{\Lambda_i}V<\inf_{\partial\Lambda_i}V$.
\vskip 0.05 true cm
\noindent(H$_6$) There exists $\delta'>0$ such that $(1+\delta')f(t)t\leq f'(t)t^2$ and $(1+\delta')g(t)t\leq g'(t)t^2$ for all $t\neq0$.
\vskip 0.05 true cm
\noindent(H$_7$) For every $\mu>0$, there exists $C_\mu>0$ such that
$$|f(s)t|+|g(t)s|\leq \mu(s^2+t^2)+C_\mu(f(s)s+g(t)t), \quad\ s,t\in\mathbb{R}.$$
\begin{Remark} The assumption (H$_7$) is technical in dealing with multi-peak solutions of problem (\ref{1.5}). Actually, in \cite{Ramos2}, one truncated argument was adapted to guarantee the modified nonlinearity satisfies condition (H$_7$), together with uniformly $L^\infty$-estimates of solutions, condition (H$_7$) was removed. In contrast to polynomial growth, it seems that exponential growth case becomes much more tough. Indeed, for general nonlinearity with exponential growth, the modified nonlinearity does not satisfy (H$_7$) by using similar truncation argument.
\end{Remark}

\begin{theo}\label{t1.4}{\it Assume (V$_1$), (V$_3$), (H$_1$), (H$_3$)-(i), (H$_6$) and (H$_7$) are satisfied.  Then for sufficiently small $\epsilon>0$,\\
\noindent(1) system (\ref{1.1})
has a nontrivial solution $(\varphi_\epsilon,\psi_\epsilon)$ in $H^1(\mathbb{R}^2)\times H^1(\mathbb{R}^2)$;\\
\noindent(2) If additionally $V$ is locally H\"{o}lder continuous, then \\
\noindent(i) $|\varphi_\epsilon|$, $|\psi_\epsilon|$ and $|\varphi_\epsilon|+|\psi_\epsilon|$ possess $k$ local maximum points $x_{i,\epsilon}$, $y_{i,\epsilon}$ and $z_{i,\epsilon}$ respectively in $\Lambda_i$ for $i=1,...,k$, and $|\varphi_\epsilon|+|\psi_\epsilon|$ does not have  local maximum points in $\mathbb{R}^2\backslash{\cup^k_{i=1}\Lambda_i}$. Moreover, $$V(x_{i,\epsilon})\rightarrow\inf_{\Lambda_i}V, \ V(y_{i,\epsilon})\rightarrow\inf_{\Lambda_i}V, \ V(z_{i,\epsilon})\rightarrow\inf_{\Lambda_i}V, \ \text{as}\ \epsilon\rightarrow0.$$
\\
\noindent(ii)
$$\|\varphi_\epsilon(\epsilon x)-\Sigma^k_{i=1}\bar{u}_i(\cdot-x_{i,\epsilon}/\epsilon)\|_{H^1(\mathbb{R}^2)}\rightarrow0, \ \|\psi_\epsilon(\epsilon x)-\Sigma^k_{i=1}\bar{v}_i(\cdot-x_{i,\epsilon}/\epsilon)\|_{H^1(\mathbb{R}^2)}\rightarrow0,$$
with $(\bar{u}_i,\bar{v}_i)$ ground states of
\begin{equation}\label{bu2}\aligned
\left\{ \begin{array}{lll}
-\Delta u+V_iu=g(v)\ & \text{in}\quad \mathbb{R}^2,\\
-\Delta v+V_iv=f(u)\ & \text{in}\quad \mathbb{R}^2,
\end{array}\right.\endaligned
\end{equation}
and the above remains true with $x_{i,\epsilon}$ replacing by $y_{i,\epsilon}$ or $z_{i,\epsilon}$.
\\
\noindent(iii) $|\varphi_\epsilon(x)|+|\psi_\epsilon(x)|\leq Ce^{-\frac{c}{\epsilon}|x-z_{i,\epsilon}|}$,\ $\forall x\in \mathbb{R}^2\backslash{\cup_{j\neq i}\Lambda_j}$ for some $C, c>0$.
}\end{theo}

\begin{theo}\label{t1.5} {\it Let the assumptions of Theorem \ref{t1.4} and (H$_5$) be satisfied. Then for sufficiently small $\epsilon>0$, $\varphi_\epsilon$, $\psi_\epsilon$, $x_{i,\epsilon}$, $y_{i,\epsilon}$ and $z_{i,\epsilon}$ obtained in Theorem \ref{t1.4} satisfy\\
\noindent(1) $\varphi_\epsilon>0$, $\psi_\epsilon>0$ in $\mathbb{R}^2$.\\
\noindent(2) $x_{i,\epsilon}$, $y_{i,\epsilon}$ and $z_{i,\epsilon}$ are unique for any $i=1,...,k$. In addition, $\lim_{\epsilon\rightarrow0}|x_{i,\epsilon}-y_{i,\epsilon}|/\epsilon=0$ for any $i=1,...,k$.}\end{theo}
\begin{Remark}\label{r1.4}We shall adapt the ideas in \cite[Chapter 5]{HUGO-PHD} (see also \cite{Ramos2}) to show Theorems \ref{t1.4} and \ref{t1.5}, while we are required to deal with the difficulties brought by exponential growth terms. Besides the difficulty mentioned in Remark \ref{r1.2}-(2), we also need to overcome the following difficulties:
\\
\noindent(1) in \cite[Chapter 5]{HUGO-PHD}, the author use the relation between $\|\nabla u\|_{L^2(\mathbb{R}^2)}$ and $\|u\|_{L^{2^*}(\mathbb{R}^2)}$ to show the least energy is bounded from below, while inspired by \cite{Cassani2}, we make use of the relation between $\|\nabla u\|_{L^2(\mathbb{R}^2)}$ and $\int_{\mathbb{R}^2}|u|^2(e^{u^2}-1)$, together with some new estimates to get the desired results.\\
\noindent(2) Compared with \cite[Chapter 5]{HUGO-PHD}, another distinct feature is: the limit profile of solutions are established. For this, we expand the regions where the cut-off functions are equal to one, and take advantage of limit decomposition of solutions to show Theorem \ref{t1.4}-(2)(ii).
\end{Remark}

Finally we consider a special case of $f,g$ satisfying (H$_3$)-(i) and then (H$_7$) shall be removed.
\vskip 0.1 true cm
(H$_8$)  $\lim_{|t|\rightarrow\infty}{f'(t)}/{e^{t^{p}}}=0$ and $\lim_{|t|\rightarrow\infty}{g'(t)}/{e^{t^{q}}}=0$ for some $0<p,q<2$.
\vskip 0.1 true cm
Without loss of generality, assume $\max\{p,q\}=1$.

\begin{theo}\label{t1.6} {\it Assume (V$_1$), (V$_3$), (H$_1$), (H$_6$) and (H$_8$) are satisfied and $V$ is locally H\"{o}lder continuous, then the results of Theorem \ref{t1.4} remain true. If additionally (H$_5$) is satisfied, then the results of Theorem \ref{t1.5} remain true. }\end{theo}

The paper is organized as
follows. In Section 2 we introduce the variational framework. In Section 3 we study autonomous systems. In Section 4, we prove Theorems \ref{t1.1}, \ref{t1.2} and \ref{t1.3}. In Section 5, we prove Theorems \ref{t1.4} and \ref{t1.5}. In Section 6, we show Theorem \ref{t1.6}.

\section{Variational setting}
\renewcommand{\theequation}{2.\arabic{equation}}
In this paper we use the following
notations. Denote the norm in $L^r(\mathbb{R}^2)$
($1\leq r\leq\infty$) by $|\cdot|_r$, the
norm in $L^p(\Omega)$ with $\Omega\subseteq\mathbb{R}^2$ is denoted by $|\cdot|_{p,\Omega}$. For the subset $\Omega\subset\mathbb{R}^2$, denote $\mathbb{R}^2\backslash{\Omega}$ by $\Omega^c$. For simplicity, denote $\int_{\mathbb{R}^2}f(x)dx$ by $\int_{\mathbb{R}^2}f(x)$, denote $\frac{\partial h(x,y)}{\partial x}$ and $\frac{\partial h(x,y)}{\partial y}$ by $h'_1(x,y)$ and $h'_2(x,y)$ respectively. Without loss of generality, we assume $0\in \Lambda_0$ in the condition (V$_2$). $C, C_1, C_2, C',C''$,... may denote positive constants whose precise value
can change from line to line.

By making the change of variable $x\rightarrow\epsilon x$, the problem (\ref{1.1})
turns out to be
\begin{equation}\label{2.1}\aligned
\left\{ \begin{array}{lll}
-\Delta u+V(\epsilon x)u=g(v)\ & \text{in}\quad \mathbb{R}^2,\\
-\Delta v+V(\epsilon x)v=f(u)\ & \text{in}\quad \mathbb{R}^2.
\end{array}\right.\endaligned
\end{equation}
Let $H^1(\mathbb{R}^2)$ be the Sobolev space endowed with the inner product and norm
$$(u,v)_\epsilon=\int_{\mathbb{R}^2}(\nabla u\nabla v+V(\epsilon x)uv),\ \|u\|^2_{\epsilon}=(u,u)_{\epsilon}, \ u,v\in H^1(\mathbb{R}^2).$$
The standard norm in $H^1(\mathbb{R}^2)$ is denoted by $\|u\|^2=\int_{\mathbb{R}^2}(|\nabla u|^2+u^2)$,
and $E:=H^1(\mathbb{R}^2)\times H^1(\mathbb{R}^2)$  be the Sobolev space endowed with the inner product
$$(z_1,z_2)_{1,\epsilon}=(u_1,u_2)_{\epsilon}+(v_1,v_2)_{\epsilon},\quad z_i=(u_i,v_i)\in E,\ i=1,2.$$
It is easy to see that there is a space decomposition of $E$ that $E=E^+\oplus E^-$, where
$$E^+=\{(u,u): u\in H^1(\mathbb{R}^2)\}, \ E^-=\{(u,-u): u\in H^1(\mathbb{R}^2)\}.$$
For each $z=(u,v)\in E$, one has
$$z=z^++z^-=\bigl((u+v)/2,(u+v)/2\bigr)+\bigl((u-v)/2,(v-u)/2\bigr).$$
The functional of (\ref{2.1}) is
\begin{equation}\label{2.2}\aligned \mathcal{L}_\epsilon(z)&=\int_{\mathbb{R}^2}[\nabla u\nabla v+V(\epsilon x)uv]-\int_{\mathbb{R}^2}[F(u)+G(v)]\\
&=\frac12\|z^+\|^2_{1,\epsilon}
-\frac12\|z^-\|^2_{1,\epsilon}-\int_{\mathbb{R}^2}[F(u)+G(v)].\endaligned
\end{equation}
\begin{lemma}\label{l2.0} If (H$_1$) and (H$_2$) are satisfied, then
\begin{equation*}f(t)t>2F(t)>0\ \text{and}\  g(t)t>2G(t)>0,\ \ \quad\forall t\neq0.\end{equation*}
If (H$_1$) and (H$_6$) are satisfied, then
\begin{equation*} f(t)t\geq (2+\delta')F(t)>0\ \text{and}\  g(t)t\geq (2+\delta')G(t)>0, \quad\forall t\neq0.
\end{equation*}
\end{lemma}
\subsection{The modified problem}
Note that $f(t),\ g(t)$ may not be equal to zero if $t<0$, so it is necessary to give a new penalized nonlinearity different from \cite{DINGXU1,HUGO-PHD}. For $i=1,2$, fix small numbers $a_i>0, b_i<0$, in such a way that $f'(a_1), f'(b_1)\leq \frac{\inf_{\mathbb{R}^2}V}{2}$,  $f'(t)\geq f'(a_1)$ for any $t\geq a_1$, $f'(t)\geq f'(b_1)$ for any $t\leq b_1$, and $g'(a_2), g'(b_2)\leq \frac{\inf_{\mathbb{R}^2}V}{2}$,  $g'(t)\geq g'(a_2)$ for any $t\geq a_2$, $g'(t)\geq g'(b_2)$ for any $t\leq b_2$.
Set
$$\aligned\tilde{f}(t)=\left\{\begin{array}{lll}f(t),\ \ &\text{if} \ 0\leq t\leq a_1\ \text{or}\ b_1\leq t\leq0;\\ f'(a_1)t+f(a_1)-f'(a_1)a_1,\ \ &\text{if} \ t>a_1;\\ f'(b_1)t+f(b_1)-f'(b_1)b_1,\ \ &\text{if} \ t< b_1,
\end{array}\right.\endaligned$$
$$\aligned\tilde{g}(t)=\left\{\begin{array}{lll}g(t),\ \ &\text{if} \ 0\leq t\leq a_2\ \text{or}\ b_2\leq t\leq0;\\ g'(a_2)t+g(a_2)-g'(a_2)a_2,\ \ &\text{if} \ t>a_2;\\
f'(b_2)t+f(b_2)-f'(b_2)b_2,\ \ &\text{if} \ t< b_2,\end{array}\right.\endaligned$$
Then we introduce
$$\bar{f}(x,t)=\chi_{\Lambda_0}(x)f(t)
+(1-\chi_{\Lambda_0}(x))\tilde{f}(t),\quad \bar{g}(x,t)=\chi_{\Lambda_0}(x)g(t)
+(1-\chi_{\Lambda_0}(x))\tilde{g}(t),$$
and $\bar{F}(x,t)=\int^t_0\bar{f}(x,s)ds$, $\bar{G}(x,t)=\int^t_0\bar{g}(x,s)ds$. The relevant properties of $\bar{f}$ and $\bar{g}$ are
displayed in the next lemma, whose proof is elementary.
\begin{lemma}\label{l2.1} The function $\bar{f}(x,t)$ ( and also $\bar{g}(x,t)$) satisfies:\\
\noindent(H$''_1$) $\bar{f}(x,t)=o(t)$\ uniformly in\
$x\in\mathbb{R}^2$, and  $|\bar{f}(x,t)|\leq |f(t)|$ for all $x\in\mathbb{R}^2$ and $t\neq0$;\\
\noindent(H$''_2$) (i) if (H$_3$) is satisfied, then $0< \theta \bar{F}(x,t)\leq \bar{f}(x,t)t,$ for all $x\in \Lambda_0$ and $t\neq0$;\\
\  \ \indent\  (ii) if (H$'_3$) is satisfied, then
 $0< 2 \bar{F}(x,t)< \bar{f}(x,t)t,$ for all $x\in \Lambda_0$ and $t\neq0$;\\
\noindent(H$''_3$) $0< 2\bar{F}(x,t)<\bar{f}(x,t)t\leq \frac{\inf_{\mathbb{R}^2}V}{2}t^2,$ for all $x\not\in\Lambda_0$ and $t\neq0$;\\
\noindent(H$''_4$) $\bar{f}(x,t)$ is nondecreasing in $t\in \mathbb{R}$ for all $x\in\mathbb{R}^2$; \\
\noindent(H$''_5$) set $\hat{F}(x,t)=\frac12\bar{f}(x,t)t-\bar{F}(x,t)$, for all $(x,t)\in\mathbb{R}^2\times(0,+\infty).$ Then $\hat{F}(x,t)$ is nondecreasing in $t\in (0,+\infty)$ and nonincreasing in $t\in(-\infty,0)$  for all $x\in\mathbb{R}^2$ and  $\hat{F}(x,t)\rightarrow+\infty$ uniformly in $x$ as $|t|\rightarrow+\infty$;\\
\noindent(H$''_6$) for some (arbitrarily small) $\delta=\delta(a_1,a_2,b_1,b_2)>0$
\begin{equation}\label{2.2.1}
|\bar{f}(x,t)|\leq\delta|t|, \ |\bar{f}'_2(x,t)|\leq\delta, \quad\text{for all}\ x\not\in{\Lambda_0}\ \text{and}\ t\neq0;\end{equation}
\end{lemma}

Now we establish the modified problem
\begin{equation}\label{2.3}\aligned
\left\{ \begin{array}{lll}
-\Delta u+V(\epsilon x)u=\bar{g}(\epsilon x, v)\ & \text{in}\quad \mathbb{R}^2,\\
-\Delta v+V(\epsilon x)v=\bar{f}(\epsilon x, u)\ & \text{in}\quad \mathbb{R}^2.
\end{array}\right.\endaligned
\end{equation}
Denote $\Lambda^\epsilon_0=\{x\in\mathbb{R}^2:\epsilon x\in \Lambda_0\}$, then the solution $(u,v)$ of (\ref{2.3}) with $|u(x)|\leq \min\{a_1,-b_1\}$ and $|v(x)|\leq \min\{a_2,-b_2\}$ for each $x\in {\mathbb{R}^2\backslash{\Lambda^\epsilon_0}}$ is also a solution of (\ref{2.1}). The functional of (\ref{2.3}) is
$${\Phi}_\epsilon(z)=\int_{\mathbb{R}^2}(\nabla u\nabla v+V(\epsilon x)uv)-\int_{\mathbb{R}^2}[\bar{G}(\epsilon x,v)+\bar{F}(\epsilon x,u)],\quad \forall z=(u,v)\in E,$$
 and ${\Phi}'_\epsilon\in C^1(E, \mathbb{R})$.

 \subsection{Functional reduction}
Since the functional $\Phi_\epsilon$ is strongly indefinite, we shall apply the reduction approach, see for example \cite{Ramos1,DINGXU2}, to look for critical  points of $\Phi_\epsilon$. More precisely, we shall  reduce the strongly indefinite functional $\Phi_\epsilon$ to a
functional on $E^+$.

 For any fixed $(u,u)\in E^+$, let $\phi_{(u,u)}:E^-\rightarrow\mathbb{R}$ defined by $$\phi_{(u,u)}(v,-v)=\Phi_\epsilon(u+v,u-v).$$ Then
\begin{equation}\label{3.1}\aligned
\phi_{(u,u)}(v,-v)
&=\|u\|^2_\epsilon-\|v\|^2_\epsilon-\int_{\mathbb{R}^2}[\bar{F}(\epsilon x,u+v)+\bar{G}(\epsilon x,u-v)]\leq\|u\|^2_\epsilon-\|v\|^2_\epsilon.
\endaligned
\end{equation}
Moreover, for any $(v,-v)$, $(w,-w)\in E^-$, by (H$''_4$) we have
\begin{equation}\label{3.2}\aligned
\langle\phi''_{(u,u)}(v,-v),(w,-w)\rangle&=-2\|w\|^2_\epsilon-
\int_{\mathbb{R}^2}[\bar{g}'_2(\epsilon x,u-v)+\bar{f}'_2(\epsilon x,u+v)]w^2\leq-2\|w\|^2_\epsilon.
\endaligned
\end{equation}
Due to (\ref{3.1}) and (\ref{3.2}), there exists a unique ${h}_\epsilon(u)\in H^1(\mathbb{R}^2)$ such that
$$\Phi_\epsilon(u+{h}_\epsilon(u),u-{h}_\epsilon(u))=\max_{v\in H^1(\mathbb{R}^2)}\Phi_\epsilon(u+v,u-v).$$
Consequently, the operator ${h}_\epsilon: H^1(\mathbb{R}^2)\rightarrow H^1(\mathbb{R}^2)$ is well defined, and
\begin{equation}\label{3.4}\langle\Phi'_\epsilon(u+{h}_\epsilon(u),
u-{h}_\epsilon(u)),(\varphi,-\varphi)\rangle=0, \quad \forall \varphi\in H^1(\mathbb{R}^2).\end{equation}

\begin{lemma}\label{l2.4} Let (V$_1$), (H$_1$), (H$_2$), (H$_3$) or (H$'_3$) hold. Then ${h}_\epsilon\in C^1(H^1(\mathbb{R}^2), H^1(\mathbb{R}^2))$.
\end{lemma}
{\bf Proof}: Define the map
$\mathcal{H}:E\times E^-\rightarrow E^-$, $\mathcal{H}((u,v),(\psi,-\psi))=P\circ\Phi'_\epsilon(u+\psi,u-\psi),$
where $P$ is the projection from $E$ to $E^-$.
Note that the partial derivative of $\mathcal{H}$ with respect to the second variable is
$$\mathcal{H}'_2((u,v),(\psi,-\psi))(\varphi,-\varphi)=P\circ \Phi''_\epsilon(u+\psi,v-\psi)(\varphi,-\varphi):=T(\varphi,-\varphi).$$
Then
$$\langle T(\phi,-\phi),(\varphi,-\varphi)\rangle=-2(\phi,\varphi )_\epsilon-\int_{\mathbb{R}^2}[\bar{f}'_2(\epsilon x,u+\psi)+\bar{g}'_2(\epsilon x,u-\psi)]\phi\varphi,\quad\forall (\phi,\varphi)\in E.$$
We claim that $T$ is one-to-one, since if $T(\phi,-\phi)=(0,0)$, then
$\langle T(\phi,-\phi),(\phi,-\phi) \rangle=0$. Using (H$''_4$) we have $\phi=0$. On the other hand, we have
$$(id+T)((\phi,-\phi),(\varphi,-\varphi))=
\int_{\mathbb{R}^2}\bar{f}'_2(\epsilon x,u+\psi)\phi\varphi+\int_{\mathbb{R}^2}\bar{g}'_2(\epsilon x,v-\psi)\phi\varphi,$$
where $id: E^-\rightarrow(E^-)^*$, $id((\phi,-\phi),(\varphi,-\varphi))=2(\phi,\varphi)_\epsilon$ for all
$\phi,\varphi\in H^1(\mathbb{R}^2)$.
In view of (H$''_1$) and (H$_3$)-(i), for any fixed $\alpha>0$, then for any $\delta>0$, there exists $C_\delta>0$ such that
\begin{equation}\label{2.2.0}|\bar{f}'_2(x,t)|\leq \delta+C_\delta(e^{\alpha t^2}-1), \ \forall (x,t)\in\mathbb{R}^2\times\mathbb{R}.\end{equation}
If $(\phi_n,-\phi_n)\rightharpoonup (\phi_0,-\phi_0)$ in $E^-$, then it is easy to see
$$\aligned\int_{\mathbb{R}^2}|\bar{f}'_2(\epsilon x,u+\psi)||\phi_n-\phi_0||\varphi|
=o_n(1)\|\varphi\|_\epsilon.\endaligned$$
Similar results hold for $\bar{g}'_2$. So $id+T$ is compact. From the Fredholm alternative we conclude that $T$ is isomorphism.
Thus the implicit function theorem implies $h_\epsilon$ is of class $C^1$.\ \ \ \ $\Box$

By Lemma \ref{l2.4}, we consider the reduced functional $J_\epsilon: H^1(\mathbb{R}^2)\rightarrow\mathbb{R}$ defined by
\begin{equation}\label{2.4} J_\epsilon(w):=
\Phi_\epsilon(w+{h}_\epsilon(w),w-{h}_\epsilon(w))\end{equation}
which is of class $C^1$. $(J_\epsilon, h_\epsilon)$ is called the reduced couple of $\Phi_\epsilon$ for simplicity. Moreover, by (\ref{3.4}) we have
\begin{equation}\label{2.2.2}\aligned\langle J'_\epsilon(w),\phi\rangle&= \langle\Phi'_\epsilon(w+{h}_\epsilon(w),w-{h}_\epsilon(w)),
(\phi+{h}'_\epsilon(w)\phi,\phi-{h}'_\epsilon(w)\phi)\rangle\\
&= \langle\Phi'_\epsilon(w+{h}_\epsilon(w),w-{h}_\epsilon(w)),
(\phi,\phi)\rangle, \quad \forall\phi\in H^1(\mathbb{R}^2).\endaligned\end{equation}

\begin{lemma}\label{l2.5}(1) The map $$\eta_\epsilon: H^1(\mathbb{R}^2)\rightarrow E: u\rightarrow (u+{h}_\epsilon(u), u-{h}_\epsilon(u))$$
is a homeomorphism between critical points of $J_\epsilon$ and $\Phi_\epsilon$, and $\eta^{-1}_\epsilon: E\rightarrow H^1(\mathbb{R}^2)$ is given by $\eta^{-1}_\epsilon(u,v)=\frac{u+v}{2}$.\\
\noindent(2) If $\{w_n\}$ is a (PS)$_c$ sequence with any $c\in\mathbb{R}$ for $J_\epsilon$, then $\{(w_n+{h}_\epsilon(w_n),w_n-{h}_\epsilon(w_n))\}$ is
a (PS)$_c$ sequence for $\Phi_\epsilon$. \\
\noindent(3) If $f$ and $g$ satisfy (H$'_3$), and  $\{w_n\}$ is a (Ce)$_c$ sequence with any $c\in\mathbb{R}$ for $J_\epsilon$, then $\{(w_n+{h}_\epsilon(w_n),w_n-{h}_\epsilon(w_n))\}$ is
a (Ce)$_c$ sequence for $\Phi_\epsilon$.
\end{lemma}
{\bf Proof}: The proof of (1) is easy and we omit. For any $(\phi,\varphi)\in E$, by (\ref{3.4}) and (\ref{2.2.2}) we get
$$\aligned&\langle \Phi'_\epsilon(w_n+{h}_\epsilon(w_n),w_n-{h}_\epsilon(w_n)),
(\phi,\varphi)\rangle=\langle J'_\epsilon(w_n),({\phi+\varphi})/{2}\rangle.
\endaligned$$
Then the conclusion (2) holds true. Regarding (3), firstly  note that
\begin{equation}\label{3.1.3}\aligned 0\leq \Phi_\epsilon(u+{h}_\epsilon(u),u-{h}_\epsilon(u))-\Phi_\epsilon(u,u)
\leq
\int_{\mathbb{R}^2}[\bar{F}(\epsilon x,u)+\bar{G}(\epsilon x,u)]-\|{h}_\epsilon(u)\|^2_\epsilon.
\endaligned\end{equation}
 If (H$'_3$) is satisfied, by (\ref{3.1.3}) we know
$\|{h}_\epsilon(u)\|^2_\epsilon\leq C\|u\|^2_\epsilon.$
Together with the conclusion (2),  the conclusion (3) yields. \ \ \  \ $\Box$

According to Lemma \ref{l2.5}, it suffices to look for critical points of $J_\epsilon$ and we shall show that $J_\epsilon$ possesses the mountain-pass structure.
\begin{lemma}\label{l2.6} There are $r>0$ and $\tau>0$ both independent of $\epsilon$, such that $J_\epsilon|_{S_r}\geq\tau$, where
$S_r=\{u\in H^1(\mathbb{R}^2):\|u\|_\epsilon=r\}$.\end{lemma}
{\bf Proof}:  By (H$_3$)-(i), it is easy to see that $g$ and $f$ are also of subcritical exponential growth. Then by (H$''_1$), for any fixed $\alpha>0$, then for any $\delta>0$ and any $q\geq1$, there exist $C(\delta,q)$ such that
\begin{equation}\label{2.1.0} |\bar{f}(\epsilon x,t)|\leq |f(t)|\leq \delta|t|+C(\delta,q)|t|^{q-1}(e^{\alpha t^2}-1),\quad \forall(x,t)\in \mathbb{R}^2\times\mathbb{R},\end{equation}
and the same inequality holds for $g$ and $\bar{g}$.
For any $w\in H^1(\mathbb{R}^2)$, by (H$''_2$), (\ref{2.1.0}), the H\"{o}lder inequality and Trudinger-Moser inequality we obtain
\begin{equation*}\aligned J_\epsilon(w)&\geq\Phi_\epsilon(w,w)\geq\|w\|^2_\epsilon-2\delta|w|^2_2-C(\delta,q)
|w|^q_{2q}\int_{\mathbb{R}^2}\bigl(e^{2\alpha w^2}-1\bigr)\\
&\geq\|w\|^2_\epsilon-2\delta|w|^2_2-C|w|^q_{2q}.
\endaligned\end{equation*}
Choosing $q>2$ in (\ref{2.1.0}), for some $r,\tau>0$ we have $J_\epsilon(w)\geq\tau>0$ when $\|w\|_\epsilon=r$. \ \ \ \ \ $\Box$

\begin{lemma}\label{l3.6} For any $u\in H^1(\mathbb{R}^2)\backslash\{0\}$, the following results hold.\\
\noindent(1) If (H$_3$) is satisfied, then\\
\ \indent(i) $J_\epsilon(tu)\rightarrow-\infty$ as $t\rightarrow+\infty$ if $suppu \cap \Lambda^\epsilon_0\neq\emptyset$;\\
\ \indent(ii) $J_\epsilon(tu)\rightarrow-\infty$ or $J_\epsilon(tu)\rightarrow+\infty$ as $t\rightarrow+\infty$ if $supp u\subset \mathbb{R}^2\backslash{\Lambda^\epsilon_0}$.\\
\noindent(2) If (H$'_3$) is satisfied, then $J_\epsilon(tu)\rightarrow-\infty$ or $J_\epsilon(tu)\rightarrow+\infty$ as $t\rightarrow+\infty$.
\end{lemma}

{\bf Proof}: (1) (i) Assume $supp u\cap \Lambda^\epsilon_0\neq\emptyset$.  Note that
$$\aligned J_\epsilon(tu)
&\leq t^2\|u\|^2_\epsilon-\|{h}_\epsilon(tu)\|^2_\epsilon-
\int_{\Lambda^\epsilon_0}{F}(tu+{h}_\epsilon(tu))-
\int_{\Lambda^\epsilon_0}{G}(tu-{h}_\epsilon(tu)).\endaligned$$
By (H$_3$)-(ii) and (H$_1$), for any $\delta>0$ there exist $C_{1,\delta}$ and $C_{2,\delta}>0$ such that
$$F(s)\geq C_{1,\delta} |s|^\theta-{\delta} s^2, \ G(s)\geq C_{2,\delta} |s|^\theta-{\delta} s^2,\ \forall s\neq0.$$
Then
$$\aligned J_\epsilon(tu)
\leq&({1+2\delta})t^2\|u\|^2_\epsilon-\min\{C_{1,\delta},C_{2,\delta}\}
\int_{\Lambda^\epsilon_0}2^\theta t^\theta|u|^\theta.
\endaligned$$
Hence
$J_\epsilon(tu)\rightarrow -\infty$ as $t\rightarrow+\infty$.

(ii) Suppose $supp u\subset \mathbb{R}^2\backslash{\Lambda^\epsilon_0}$.
If $J_\epsilon(tu)\rightarrow+\infty$, we are done. Otherwise,
we may assume $\sup_{t\geq0}J_\epsilon(tu)=M<+\infty$.
For $r>0$, by (H$''_5$) we infer
\begin{equation}\label{2.12.0}\aligned
&\int_{\mathbb{R}^2}\hat{F}(\epsilon x, tu+{h}_\epsilon(tu))\geq
 \min\{\hat{F}(\epsilon x,rt),\hat{F}(\epsilon x,-rt)\}meas(\Omega),
\endaligned\end{equation}
where $\Omega=\{x\in\mathbb{R}^2:|u+{{h}_\epsilon(tu)}/{t}|\geq r\}$.
Similar results hold for $\hat{G}(\epsilon x, tu-{h}_\epsilon(tu))$.
Taking into account of (\ref{3.1.3}), $supp u\subset\mathbb{R}^2\backslash{\Lambda^\epsilon_0}$ and (H$''_3$) we know
$\|{h}_\epsilon(tu)\|^2_\epsilon\leq \frac{\inf_{\mathbb{R}^2}V}{2}t^2|u|^2_2.$
Then $\{\frac{{h}_\epsilon(tu)}{t}\}_{t>0}\subset H^1(\mathbb{R}^2)$ is bounded, which combines with the fact that $u\neq0$ imply that either $meas\{x\in\mathbb{R}^2:
|u-{{h}_\epsilon(tu)}/{t}|\geq r\}\geq \delta$  or $meas(\Omega)\geq \delta$ with some $\delta>0$ for all $t>0$ provided $r>0$ small. If the former holds, we treat $\hat{G}(\epsilon x, tu-{h}_\epsilon(tu))$ to get a contradiction similarly, and below we assume $meas(\Omega)\geq \delta$. By (H$''_5$) and (\ref{2.12.0}) we get
$$\frac12\langle J'_\epsilon(tu),tu\rangle\leq J_\epsilon(tu)-\int_{\mathbb{R}^2}\hat{F}(\epsilon x, tu+{h}_\epsilon(tu))\leq -M.$$
So
$J_\epsilon(tu)=\int^t_0\langle J'_\epsilon(tu),u\rangle dt\rightarrow-\infty$, as\ $t\rightarrow+\infty,$ contradicts the hypothesis $J_\epsilon(tu)\leq M$.

(2) If (H$'_3$) is satisfied, then taking similar arguments as in the proof of the conclusion (1)-(ii), we obtain desired results.
\ \ \ \ \ $\Box$

In the same way as the proof of Claim 1 in \cite[Lemma 4.6]{DINGXU2}, we have the following result.
\begin{lemma}\label{l3.7} Under the assumptions of Lemma \ref{l2.4}, if $u\in H^1(\mathbb{R}^2)\backslash\{0\}$ satisfies $\langle J'_\epsilon(u),u\rangle=0$, then $\langle J''_\epsilon(u)u,u\rangle<0$.
\end{lemma}

\section{The autonomous systems}
\renewcommand{\theequation}{3.\arabic{equation}}
\subsection{The limit system}
For any $\mu>0$, consider the limit system
\begin{equation}\label{3.0}\aligned
\left\{ \begin{array}{lll}
-\Delta u+\mu u=g(v)\ & \text{in}\quad \mathbb{R}^2,\\
-\Delta v+\mu v=f(u)\ & \text{in}\quad \mathbb{R}^2,
\end{array}\right.\endaligned
\end{equation}
 and define another inner product in $H^1(\mathbb{R}^2)$ and in $E$ by
 $(u,v)_\mu=\int_{\mathbb{R}^2}(\nabla u\nabla v+\mu uv),$ and $$\ (z_1,z_2)_{1,\mu}=(u_1,u_2)_{\mu}+(v_1,v_2)_{\mu},\ z_i=(u_i,v_i)\in E, i=1,2,$$
 respectively.
The functional  of (\ref{3.0}) is
$$\Phi_\mu(z)=\frac12(\|z^+\|^2_{1,\mu}-\|z^-\|^2_{1,\mu})
-\int_{\mathbb{R}^2}[G(v)+F(u)], \ \forall z=z^++z^-\in E.$$
Define
\begin{equation}\label{3.1.0}c_\mu:=\inf_{\mathcal{K}_\mu}\Phi_\mu,\quad \text{where}\ \mathcal{K}_\mu=\bigl\{ z\in E\backslash\{(0,0)\}:\ \Phi'_\mu(z)=0\bigr\}.\end{equation}
We would like to point out that similar to the arguments in \cite{ZhangCPDE} for the superlinear case and in \cite{LY} for the asymptotically linear case, we can show that  (\ref{3.0})  has a ground state and weaken $f$ and $g$ to be merely continuous.
In this paper, to give another characterization of the least energy $c_\mu$ to restore the compactness of (\ref{2.3}), we assume $f, g\in C^1$ and introduce the mapping ${h}_\mu$.
As (\ref{3.4}) and (\ref{2.4}), we define $h_\mu, J_\mu$ respectively and call $(J_\mu, h_\mu)$ the reduced couple of $\Phi_\mu$.

  \vskip 0.1 true cm
{\bf 3.1.1. The super-linear case.}
\begin{lemma}\label{l3.2.0}
Let $\mu>0$ and (H$_1$)-($H_3$) hold. Then (\ref{3.0}) admits a ground state in $E$. Moreover,
$c_\mu=\bar{c}_\mu=\bar{\bar{{c}}}_\mu>0,$
where \begin{equation*}\label{3.18.0}\bar{c}_\mu:=\inf_{\nu\in \Gamma_\mu}\max_{t\in[0,1]}J_\mu(\nu(t))\ \text{with}\ \Gamma_\mu=\{\nu\in C([0,1], H^1(\mathbb{R}^2)):\nu(0)=0, J_\mu(\nu(1))<0\},\end{equation*} and
$$\bar{\bar{{c}}}_\mu:=\inf_{u\in H^1(\mathbb{R}^2)\backslash\{0\}}\max_{t>0}J_\mu(tu)=\inf_{u\in N_\mu}J_\mu(u),$$
with $N_\mu:=\{u\in H^1(\mathbb{R}^2)\backslash\{0\}: \langle J'_\mu(u),u\rangle=0\}$.
\end{lemma}
{\bf Proof}: One easily has that Lemma \ref{l2.6} and Lemma \ref{l3.6} (1)-(i)  hold for $J_\epsilon$ replacing $J_\mu$. Then taking easier proof than \cite[Theorem 1.3]{DG1}, where the authors consider system (\ref{3.0}) with critical exponential growth, we can show (\ref{3.0}) admits a ground state in $E$. Moreover, by a standard argument we have $c_\mu=\bar{c}_\mu=\bar{\bar{{c}}}_\mu$.\ \ \ $\Box$

{\bf 3.1.2. The asymptotically linear case.} In this subsection, let $\mu\in(0,|V|_\infty]$.
Choose $u\in C^\infty_0(\mathbb{R}^2)$ such that $|u|_2=1$. Set $u_t(x):=tu(tx)$, $t>0$. Then $|u_t|_2=1$ and $|\nabla u_t|_2\rightarrow0$ as $t\rightarrow0$. Hence, by (H$'_3$) there exists $t_0>0$ such that $|\nabla u_{t_0}|^2_2+\mu|u_{t_0}|^2_2-l_0|u_{t_0}|^2_2<0.$
Let $e=\frac{1}{\sqrt{2}}\bigl(\frac{u_{t_0}}{\|u_{t_0}\|_\mu}, \frac{u_{t_0}}{\|u_{t_0}\|_\mu}\bigr)\in E^+,$ and $E_e=E^-\oplus\mathbb{R}e$, arguing as \cite[Lemma 6.4]{DING} and Lemma \ref{l3.6} (1)-(ii) respectively, we have the following result.

\begin{lemma}\label{l3.4.0}(1) If $z\in E_e$ and $\|z\|_{1,\mu}\rightarrow+\infty$, then $\Phi_\mu(z)\rightarrow-\infty$  and $\sup_{z\in E_e} \Phi_\mu(z)<+\infty$.\\
\noindent(2) For any $u\in H^1(\mathbb{R}^2)\backslash\{0\}$, then either $J_\mu(tu)\rightarrow+\infty$ or $J_\mu(tu)\rightarrow-\infty$ as $t\rightarrow+\infty$.\end{lemma}

\begin{lemma}\label{l3.3} Let $\mu\in(0,|V|_\infty]$, (H$_1$), (H$_2$) and (H$'_3$) hold. Then
$c_\mu>0$ is achieved and $c_\mu=\bar{c}_\mu=\bar{\bar{{c}}}_\mu$, where $\bar{c}_\mu$ and $\bar{\bar{{c}}}_\mu$ are defined as in Lemma \ref{l3.2.0}.
\end{lemma}
{\bf Proof}: By Lemma \ref{l3.4.0}, we know $J_\mu(te)\rightarrow-\infty$ as $t\rightarrow+\infty$. Then the minimax value $\bar{c}_\mu$ is well defined.
 Let $u\in H^1(\mathbb{R}^2)\backslash\{0\}$ we find that the function $t\mapsto J_\mu(tu)$ has at most one nontrivial critical point $t=t(u)>0$. Hence, if we denote
$$\mathcal{N}_\mu:=\bigl\{t(u)u: u\in H^1(\mathbb{R}^2)\backslash\{0\}, t(u)<+\infty\bigr\},$$ then
$\bar{\bar{{c}}}_\mu=\inf_{z\in \mathcal{N}_\mu}J_\mu(z).$ In view of Lemma \ref{l6.1} below with $\Phi_\epsilon$ replacing by $\Phi_\mu$ and Lemma \ref{l2.5}-(3),
and taking standard arguments as the proof of  Lemma \ref{l3.2.0}, we can show that $c_\mu>0$ is achieved and $c_\mu=\bar{c}_\mu=\bar{\bar{{c}}}_\mu$.
\ \ \ \ \ \ \ $\Box$

To investigate the properties of solutions of (\ref{2.3}), we require in addition the following results of ground states for (\ref{3.0}), whose proof is similar as in  \cite{ZhangCPDE} and we omit.
\begin{lemma}\label{lbu} \noindent (1) Let $\mu>0$ and (H$_1$)-(H$_3$) hold. \\
\noindent(i) If in addition (H$_4$) holds and $f$ and $g$ are odd, then for any ground state $(u,v)$ of system (\ref{3.0}), one has $u,v\in C^2(\mathbb{R}^2)$, $u,v>0$ are radially symmetric with respect to the origin and are radially decreasing. Moreover,
$\Delta u(0)<0$, $\Delta v(0)<0$. In addition, there exist $C, c>0$ such that $|u(x)|+|v(x)|\leq Ce^{-c|x|}$.\\
\noindent(ii) If in addition (H$_5$) holds, then the results of (i) remain true.\\
\noindent(2) Let $\mu\in(0,|V|_\infty]$ and (H$_1$), (H$_2$), (H$'_3$) hold. If in addition (H$_5$) holds, then the results of the conclusion (1)-(i) remain true.
\end{lemma}
\subsection{Modified autonomous systems}
For both the superlinear and asymptotically linear cases, we also introduce another autonomous system. We consider the system
\begin{equation}\label{3.3}\aligned
\left\{ \begin{array}{lll}
-\Delta u+u={\mu}^{-1}g({v})\ & \text{in}\quad \mathbb{R}^2,\\
-\Delta {v}+{v}={\mu}^{-1}f({u})\ & \text{in}\quad \mathbb{R}^2.
\end{array}\right.\endaligned
\end{equation} If $(u_\mu,v_\mu)$ is the solution of (\ref{3.0}), then
$\bigl({u}_\mu(\cdot/\sqrt{\mu}),
{v}_\mu(\cdot/\sqrt{\mu})\bigr)$ is a solution of (\ref{3.3}).
Moreover, $\Phi_{\mu}(u_\mu,v_\mu)=\hat{\Phi}_{\mu}\bigl({u}_\mu(\cdot/\sqrt{\mu}),
{v}_\mu(\cdot/\sqrt{\mu})\bigr)$, where $\hat{\Phi}_{\mu}$ is the associated functional of (\ref{3.3})
\begin{equation*}\hat{\Phi}_{\mu}(u,v)=
\int_{\mathbb{R}^2}\bigl(\nabla u\nabla v+ uv\bigr)-\frac{1}{\mu}\int_{\mathbb{R}^2}\bigl[G(v)+F(u)\bigr],\ \forall (u,v)\in E.\end{equation*}
Define $\hat{c}_{\mu}:=\inf_{\hat{\mathcal{K}}_\mu}\hat{\Phi}_{\mu}$, where
$\hat{\mathcal{K}}_\mu=\{z\in E\backslash\{(0,0)\}: \hat{\Phi}'_{\mu}(z)=0\}$.
Then
$\hat{c}_{\mu}={c}_{\mu}$,
where ${c}_{\mu}$ is defined in (\ref{3.1.0}).
\begin{lemma}\label{l5.0}
The map $\mu\mapsto c_\mu$ is increasing.
\end{lemma}
{\bf Proof}:  As (\ref{3.4}) and (\ref{2.4}),
$(\hat{J}_{\mu},\hat{h}_{\mu})$ denote the reduction couple of $\hat{\Phi}_{\mu}$. Similar to Lemmas \ref{l3.2.0} and \ref{l3.3},  we have
$\hat{c}_{\mu}=\inf_{u\in H^1(\mathbb{R}^2)\backslash\{0\}}\max_{t>0}\hat{J}_{\mu}(tu).$
 By virtue of the above characterization,
it is easy to see that  $\hat{c}_\mu$ is increasing in $\mu$. Note that $c_\mu=\hat{c}_\mu$. So $c_{\mu}$ is increasing in $\mu$.\ \ \ \ $\Box$

\section{Proof of Theorems \ref{t1.1}, \ref{t1.2} and \ref{t1.3}}
\renewcommand{\theequation}{4.\arabic{equation}} In this section, we shall give the proof of Theorems \ref{t1.1}, \ref{t1.2} and \ref{t1.3}.
\begin{lemma}\label{l2.2} Let (V$_1$) and (H$_1$)-(H$_3$) hold. Then for  each $\epsilon>0$, the (PS)$_c$ condition of $\Phi_\epsilon$ with $c\in\mathbb{R}$ holds. \end{lemma}
{\bf Proof}: Let $\{z_n=(u_n,v_n)\}$ be the (PS)$_c$ sequence of $\Phi_\epsilon$. Firstly we show that $\{z_n\}$ is bounded in $E$. Note that \begin{equation}\label{3.14}\aligned c+o_n(1)\|z_n\|_{1,\epsilon}&=\Phi_\epsilon(z_n)-\frac12\langle \Phi'_\epsilon(z_n),(v_n,u_n)\rangle\geq\frac{\theta-2}{2\theta}
\int_{\Lambda^\epsilon_0}\bigl(f(u_n)u_n+g(v_n)v_n\bigr).
\endaligned\end{equation}
By $\langle \Phi'_\epsilon(z_n),(0,v_n)\rangle=o_n(1)\|v_n\|_\epsilon$ we have
\begin{equation}\label{3.1.1}
\|v_n\|_\epsilon=\int_{\mathbb{R}^2}\bar{f}(\epsilon x,u_n)\frac{v_n}{\|v_n\|_\epsilon}+o_n(1).\end{equation}
In view of (H$''_1$), there exist $\beta>0$ and $C_\beta>0$ such that
$|\bar{f}(\epsilon x,t)|\leq C_\beta e^{\beta t^2}$ for all $(x,t)\in \mathbb{R}^2\times\mathbb{R}.$
Set
$$\Lambda^1_{n,\epsilon}=\{x\in \Lambda^\epsilon_0: |\bar{f}(\epsilon x, u_n)|/ C_\beta\geq e^{\frac14}\},\quad \Lambda^2_{n,\epsilon}=\{x\in \Lambda^\epsilon_0: |\bar{f}(\epsilon x, u_n)|/ C_\beta\leq e^{\frac14}\}.$$
By (H$''_1$) there exists $C_1>0$ such that for all $n$,
$|\bar{f}(\epsilon x, u_n)|\leq C_1|u_n|$, for $x\in \Lambda^2_{n,\epsilon}$. As in \cite[Lemma 3.2]{DDR}, there holds
\begin{equation}\label{bu1}\aligned
st\leq
\left\{ \begin{array}{lll}
(e^{t^2}-1)+|s|(\log|s|)^{\frac12},\ & \ t\in\mathbb{R}\text{ and }\ |s|\geq e^{\frac14};\\
(e^{t^2}-1)+\frac12s^2,\ & \ t\in\mathbb{R}\text{ and }\ |s|\leq e^{\frac14}.
\end{array}\right.\endaligned
\end{equation}
Applying the above inequality with $t=\frac{v_n}{\|v_n\|_\epsilon}$ and $s=\frac{\bar{f}(\epsilon x, u_n)}{C_\beta}$, from Trudinger-Moser inequality we get
$$\aligned\Bigl|\int_{\Lambda^\epsilon_0}\bar{f}(\epsilon x,u_n)\frac{v_n}{\|v_n\|_\epsilon}\Bigr|
\leq& C_\beta\int_{\Lambda^1_{n,\epsilon}}\frac{|\bar{f}(\epsilon x,u_n)|}{C_\beta}\Bigl[\log\bigl(\frac{|\bar{f}(\epsilon x,u_n)|}{C_\beta}\bigr)\Bigr]^{\frac12}\\
&+\frac12\int_{\Lambda^2_{n,\epsilon}}\frac{1}{C_\beta}\bar{f}^2(\epsilon x,u_n)+C_\beta
\int_{\Lambda^\epsilon_0}\bigl[e^{\frac{v^2_n}{\|v_n\|^2}}-1
\bigr]
\\ \leq& C_2+\bigl(\beta^{\frac12}+\frac{C_1}{2C_\beta}\bigr)
\int_{\Lambda^\epsilon_0}f(u_n)u_n.\endaligned$$
Moreover, by (H$''_3$) we have
$\Bigl|\int_{(\Lambda^\epsilon_0)^c}\bar{f}(\epsilon x,u_n)\frac{v_n}{\|v_n\|_\epsilon}\Bigr|
\leq\frac{\|u_n\|_\epsilon}{2}.$
From (\ref{3.1.1}) it follows that
$$\|v_n\|_\epsilon\leq C_2+\bigl(\beta^{\frac12}+\frac{C_1}{2C_\beta}\bigr)\int_{\Lambda^\epsilon_0}
f(u_n)u_n
+\frac{1}{2}\|u_n\|_\epsilon.$$
Similarly, the above inequality holds with $(v_n,u_n)$ replacing by $(u_n,v_n)$.
Using (\ref{3.14}) we get
$$\frac{1}{2}(\|u_n\|_\epsilon+\|v_n\|_\epsilon)\leq C_2+C\bigl(\beta^{\frac12}+\frac{C_1}{2C_\beta}\bigr)
\bigl(c+o_n(1)\|z_n\|_{1,\epsilon}\bigr).$$
Therefore, $\{z_n\}$ is bounded in $E$ and we assume $z_n\rightharpoonup z_0=(u_0,v_0)$ in $E$.

Next we show that $z_n\rightarrow z_0$ in $E$.
 In view of \cite[Lemma 2.1]{DDR} and (\ref{2.1.0}), one easily has that  $\Phi'_\epsilon$ is  weakly
sequentially continuous. Then $\Phi'_\epsilon(z_0)=0$. Since $\Lambda^\epsilon_0$ is bounded for any fixed $\epsilon>0$, using (\ref{2.1.0}) one easily has that
$$\int_{\Lambda^\epsilon_0}\bigl({f}(u_n )-{f}(u_0 )\bigr)(v_n-v_0)=o_n(1),\quad\int_{\Lambda^\epsilon_0}\bigl({g}(v_n )-{g}(v_0 )\bigr)(u_n-u_0)=o_n(1).$$
In addition, by (\ref{2.2.1}) we infer
$$\Bigl|\int_{(\Lambda^\epsilon_0)^c}\bigl(\bar{f}(\epsilon x,u_n )-\bar{f}(\epsilon x,u_0 )\bigr)(v_n-v_0)\Bigr|\leq \frac{\inf_{\mathbb{R}^2}V}{4}(|u_n-u_0|^2_2+|v_n-v_0|^2_2).$$
Similar inequalities hold for $\bar{g}$. Then
$$o_n(1)=\langle \Phi'_\epsilon(z_n)-\Phi'_\epsilon(z_0),z_n-z_0\rangle
\geq\frac{\inf_{\mathbb{R}^2}V}{2}(|v_n-v_0|^2_2+|u_n-u_0|^2_2)+o_n(1).$$
Thus $z_n\rightarrow z_0$ in $E$.
\ \ \ \ \ $\Box$

\begin{lemma}\label{l6.1} Let (V$_1$), (H$_1$), (H$_2$) and (H$'_3$) hold. Then for each $\epsilon>0$, the (Ce)$_c$ condition of $\Phi_\epsilon$ with $c\in\mathbb{R}$ holds.
\end{lemma}
{\bf Proof}:
  Let $\{z_n=(u_n,v_n)\}$ be the (Ce)$_c$ sequence of $\Phi_\epsilon$ with $c\in\mathbb{R}$. Firstly, we show that $\{z_n\}$ is bounded in $E$. Argue by contradiction we assume $\|z_n\|_{1,\epsilon}\rightarrow+\infty$.

{\bf Case 1}: $\frac{\|u_n\|_\epsilon}{\|v_n\|_\epsilon}\leq1$. Then $\|v_n\|_{\epsilon}\rightarrow+\infty$. Set for $r\geq0$ and $0\leq\rho<r$
$$ d(r)=\inf\{\hat{F}(\epsilon x,s):x\in\mathbb{R}^2 \ \text{and}\ |s|>r\},\ \Omega_n(\rho,r)=\{x\in\mathbb{R}^2:\rho\leq |u_n(x)|\leq r\},$$
 and
 $C^r_\rho=\inf\bigl\{\frac{\hat{F}(\epsilon x,s)}{s^2}:x\in\mathbb{R}^2,\ \rho\leq |s|\leq r\bigr\}.$
By (H$''_5$),  it is easy to see that $|\Omega_n(r,+\infty)|\rightarrow0$ as $r\rightarrow+\infty$ uniformly in $n$, and
$\int_{\Omega_n(\rho,r)}|u_n(x)|^2\leq \frac{C}{C^r_\rho}.$
For any fixed $\delta>0$, by (H$''_1$) we know
$|\bar{f}(\epsilon x,s)|\leq \delta|s|$,  for all $s\in[0,\rho_\delta].$ Then $\Bigl|\int_{\Omega_n(0,\rho_{\delta})}\frac{\bar{f}(\epsilon x,u_n)v_n}{\|v_n\|^2_\epsilon}\Bigr|\leq C\delta,$ $\Bigl|\int_{\Omega_n(\rho_\delta,r_\delta)}\frac{\bar{f}(\epsilon x,u_n)v_n}{\|v_n\|^2_\epsilon}\Bigr|\rightarrow0$ and $\Bigl|\int_{\Omega_n(r_\delta,+\infty)}\frac{\bar{f}(\epsilon x,u_n)v_n}{\|v_n\|^2_\epsilon}\Bigr|\rightarrow0$.
Hence $\int_{\mathbb{R}^2}\frac{\bar{f}(\epsilon x,u_n)v_n}{\|v_n\|^2_\epsilon}\rightarrow0$.
On the other hand, by (\ref{3.1.1})
and $\|v_n\|_\epsilon\rightarrow+\infty$, we infer
$\int_{\mathbb{R}^2}\frac{\bar{f}(\epsilon x,u_n)v_n
}{\|v_n\|^2_\epsilon}\rightarrow1.$
This is a contradiction.
Therefore, $\|z_n\|_{1,\epsilon}$ is bounded.

 {\bf Case 2}: $\frac{\|u_n\|_\epsilon}{\|v_n\|_\epsilon}>1$. If this case occurs, we get $\|u_n\|_\epsilon\rightarrow+\infty$. Proceeding as case 1 with  $\bar{f}(\epsilon x,u_n)$ and $u_n$ replacing by
 $\bar{g}(\epsilon x,v_n)$ and $v_n$ respectively we also get a contradiction.

 Taking same arguments as in Lemma \ref{l2.2} we know $\{z_n\}$ converges strongly in $E$. This ends the proof.\ \ \ \ $\Box$

To show that $J_\epsilon$ has mountain-pass structure, by Lemma \ref{l2.6}, it suffices to show that, there exists $u_0\in H^1(\mathbb{R}^2)$ independent on $\epsilon$, such that $\|u_0\|>r$ and $J_\epsilon(u_0)<0$. For this, we shall investigate the relationship between the system (\ref{2.3}) and its limit system (\ref{1.1.1}). In particular,  we need to study the relation between ${h}_\epsilon(u)$ and ${h}_0(u):={h}_{V_0}(u)$.

\begin{lemma}\label{l4.1} For any $u\in H^1(\mathbb{R}^2)$, ${h}_\epsilon(u)\rightarrow{h}_{0}(u)$  in $H^1(\mathbb{R}^2)$ as $\epsilon\rightarrow0$.\end{lemma}
{\bf Proof}: For any sequence $\epsilon_n\rightarrow 0^+$, put $v_n:={h}_{\epsilon_n}(u)$ and $v_0:={h}_{0}(u)$. We shall prove that
$v_n\rightarrow v_0$ in $H^1(\R^2)$.
Firstly, using (\ref{3.1.3}), it is easy to see that $\{v_n\}\subset H^1(\R^2)$ is bounded. Up to a subsequence, we assume
$v_n\rightharpoonup v^*$ {in} $H^1(\R^2)$.
Recall that
\begin{equation}\label{eq:20220410-e5}\langle\Phi'_{\epsilon_n}(u+v_n,
u-v_n),(\varphi,-\varphi)\rangle=0, \quad \forall \varphi\in H^1(\mathbb{R}^2).\end{equation}
Then \begin{equation}\label{4.5.0}\Phi'_{V_0}(u+v^*,
u-v^*),(\varphi,-\varphi)\rangle=0, \quad \forall \varphi\in H^1(\mathbb{R}^2).\end{equation}
Hence $v^*={h}_{0}(u)=v_0$.
Testing \eqref{eq:20220410-e5} by $\varphi=v_n$ and (\ref{4.5.0}) by $\varphi=v_0$ respectively and then subtracting, we infer
 \begin{equation}\label{eq:20220411-e0}
 \aligned
 2\|v_n-v_0\|_{\epsilon_n}^{2}
 =&o_n(1)+\int_{\R^2}\left[\bar{g}(\epsilon_n x, u-v_n)v_n -g(u-v_0)v_0\right]\\
 &-\int_{\R^2}\left[\bar{f}(\epsilon_n x, u+v_n)v_n -f(u+v_0)v_0\right]:=o_n(1)+\mathbb{A}-\mathbb{B}.
 \endaligned
 \end{equation}
 Note that
 \begin{equation*}\label{eq:20220411-e1}
 \aligned
 \mathbb{B}
 =&\int_{\R^2}\left[\bar{f}(\epsilon_n x, u+v_n)-\bar{f}(\epsilon_n x, u+v_0)\right](v_n-v_0)+\int_{\R^2}\bar{f}(\epsilon_n x, u+v_0)(v_n-v_0)\\
 &+\int_{\R^2}[\bar{f}(\epsilon_n x, u+v_n)v_0 -f(u+v_0)v_0]:=\mathbb{B}_1+\mathbb{B}_2+\mathbb{B}_3.
 \endaligned
 \end{equation*}
Using $(H''_4)$ one easily has $\mathbb{B}_1\geq0$.
Since $v_n\rightharpoonup v_0$ in $H^1(\mathbb{R}^2)$, we infer $\mathbb{B}_3=o_n(1)$. Moreover,
\begin{equation*}\label{eq:20220411-e2}
\aligned
 \mathbb{B}\geq\mathbb{B}_2=o_n(1)+\int_{(\Lambda^{\epsilon_n}_0)^{c}}\bar{f}(\epsilon_n x, u+v_0)(v_n-v_0).
\endaligned
\end{equation*}
 Similarly
 \begin{equation*}\label{eq:20220411-e5}
 \aligned
 \mathbb{A}
\leq \int_{(\Lambda^{\epsilon_n}_0)^{c}}\bar{g}(\epsilon_n x, u-v_0)(v_n-v_0)+o_n(1).
 \endaligned
 \end{equation*}
Combining with $(H''_3)$ and \eqref{eq:20220411-e0}, it follows that
$2\|v_n-v_0\|_{\epsilon_n}^{2}
 \leq o_n(1)+o_n(1)\|v_n-v_0\|_{\epsilon_n}$,
 which implies $v_n\rightarrow v_0$ in $H^1(\R^2)$.
\ \ \  \ $\Box$

\begin{lemma}\label{l4.2}For $\epsilon>0$ small enough, there exists $u_0\in H^1(\mathbb{R}^2)$ (independent of $\epsilon$) such that $\|u_0\|_\epsilon>r$, where $r$ is given in Lemma \ref{l2.6}, and $J_\epsilon(u_0)<0$.\end{lemma}
{\bf Proof}: Taking into account of Lemma \ref{l3.4.0}-(2), Lemma  \ref{l3.6} (1)-(i) with $J_\epsilon$ replacing by $J_\mu$ and Lemma \ref{l4.1}, and adapting similar arguments of \cite[Lemma 3.6-(2)]{DINGXU1}, we get the desired results.\ \ \ \ $\Box$

\begin{lemma}\label{l3.5.0} For $\epsilon$ small enough, the system (\ref{2.3}) has a nontrivial solution.\end{lemma}
{\bf Proof}: For the superlinear case, by Lemma \ref{l2.2} and Lemma \ref{l2.5}-(2),  it is easy to see that $J_\epsilon$ satisfies the (PS) condition. In view of Lemma \ref{l2.6} and Lemma \ref{l4.2}, applying the mountain-pass theorem with (PS) condition, we know
\begin{equation}\label{3.18.0}c_\epsilon=\inf_{v\in \Gamma_\epsilon}\max_{t\in[0,1]}J_\epsilon(\nu(t)),\end{equation}
where $\Gamma_\epsilon=\{\nu\in C([0,1], H^1(\mathbb{R}^2)):\nu(0)=0, J_\epsilon(\nu(1))<0\}$, is a critical value of $J_\epsilon$, also for $\Phi_\epsilon$ and $\tau\leq c_\epsilon<+\infty$ with $\tau$ given in Lemma \ref{l2.6}.

For the asymptotically linear case,
  from Lemma \ref{l6.1} and Lemma \ref{l2.5}-(3) it follows that $J_\epsilon$ satisfies (Ce) condition.
  In view of Lemma \ref{l2.6} and Lemma \ref{l4.2}, applying the mountain-pass theorem with (Ce) condition, we know $c_\epsilon$ given in (\ref{3.18.0}) is a critical value of $\Phi_\epsilon$.
 \ \ \ \ \ $\Box$

Taking similar arguments of \cite[Lemma 3.8]{DINGXU1} and \cite[Lemma 3.9]{DINGXU1}, we have the following lemma.
\begin{lemma}\label{l4.6}
(1) For $\epsilon$ small enough, $c_\epsilon=\inf_{u\in H^1(\mathbb{R}^2)\backslash\{0\}}\max_{t\geq0}J_\epsilon(tu)$, where $c_\epsilon$ is given in (\ref{3.18.0}).\\
\noindent(2) $c_\epsilon\leq c_{V_0}+o_\epsilon(1)$ as $\epsilon\rightarrow0$.
\end{lemma}

\begin{lemma} \label{l4.7} Let $z_{n}=(u_n,v_n)$ be nontrivial solutions of (\ref{2.3}) obtained in Lemma \ref{l3.5.0} with
$\epsilon_n\rightarrow0$. If in addition $V$ is locally H\"{o}lder continuous, then there is $x_{n}\in\Lambda^{\epsilon_n}_0$
such that $V(\epsilon_n x_n)\rightarrow V_0$ and the
sequence $\bar{z}_{n}(x):=z_{n}(x+x_n)$ converges strongly in $E$ to a ground
state $z_0=(u_0,v_0)$ of (\ref{1.1.1}).
 \end{lemma}
{\bf Proof}: By Lemma \ref{l4.6}, we know $c_{\epsilon_n}\leq c_{V_0}+o_n(1)$, then as the proof of Lemmas \ref{l2.2} and \ref{l6.1}, we know $\{z_n\}$ is bounded in $E$ and assume $\|z_n\|^2_{1,\epsilon}\leq M$ for some constant $M>0$. If $\{z_n\}$ is vanishing, i.e.
 $\sup_{y\in\mathbb{R}^2}\int_{B_R(y)}
 (u^2_n+v^2_n)\rightarrow0$, for all $R>0,$
then P.L. Lions compactness lemma implies that $u_n\rightarrow0$ and $v_n\rightarrow0$ in $L^r(\mathbb{R}^2)$ for any $r>2$.
For some $\alpha\in (0,\frac{8\pi}{3M})$, by (\ref{2.1.0}) with $q=1$, for any $\delta>0$, there exists $C_\delta>0$ such that
\begin{equation}\label{5.3}\aligned
\Bigl|\int_{\mathbb{R}^2}\bar{f}(\epsilon x, u_n)v_n\Bigr|&
&\leq\delta|u_n|_2|v_n|_2+C_\delta
|v_n|^{3}_{3}\Bigl(\int_{\mathbb{R}^2}
\bigl[e^{\alpha \frac{3}{2} \|u_n\|^2_\epsilon\frac{u^2_n}{\|u_n\|^2_\epsilon}}-1\bigr]\Bigr)^{\frac{2}{3}}\rightarrow0.
\endaligned\end{equation}
Similarly, we get
$\int_{\mathbb{R}^2}\bar{{g}}(\epsilon x,v_n)u_n\rightarrow0.$
Since $\langle \Phi'_\epsilon(z_n),(v_n,u_n)\rangle=0$, we know
$\|z_n\|^2_{1,\epsilon}\rightarrow0$.
In addition, similar to (\ref{5.3}) we get $\bar{F}(\epsilon x,u_n)\rightarrow0$ and $\bar{G}(\epsilon x,v_n)\rightarrow0$. Then $c_{\epsilon_n}\rightarrow0$, contradicts the fact that $c_{\epsilon_n}\geq \tau>0$.
So $\{z_n\}$ is nonvanishing and there exist $x_n\in\mathbb{R}^2$ and $\delta>0$ such that
\begin{equation}\label{4.9.0}\int_{B_1(x_n)}(u^2_n+v^2_n)\geq\delta.\end{equation}
{\bf Claim 1}: $\epsilon_nx_n\rightarrow x_0\in \Lambda_0$ and $V(x_0)=V_0$.

Letting $\bar{z}_n=z_n(x+x_n)$, then $\bar{z}_n \rightharpoonup z_0=(u_0,v_0)\neq(0,0)$. Moreover, $\bar{z}_n$ satisfies the system
\begin{equation}\label{4.2}\aligned
\left\{ \begin{array}{lll}
-\Delta u+V(\epsilon_n x+\epsilon_n x_n)u=\bar{g}(\epsilon_n x+\epsilon_nx_n,v)\ & \text{in}\quad \mathbb{R}^2,\\
-\Delta v+V(\epsilon_n x+\epsilon_n x_n)v=\bar{f}(\epsilon_n x+\epsilon_n x_n,u)\ & \text{in}\quad \mathbb{R}^2.
\end{array}\right.\endaligned
\end{equation}
According to (\ref{4.9.0}) we discuss into two cases.

{\bf Case 1:} $\int_{B_1(x_n)}\Bigl[(1-\chi_{\Lambda_0}(\epsilon_n x))u^2_n+(1-\chi_{\Lambda_0}(\epsilon_n x))v^2_n\Bigr]\geq\frac{\delta}{2}$.
In this case, $B_1(x_n)\cap \{\mathbb{R}^2\setminus supp\chi_{\Lambda_0}(\epsilon_n x) \} \neq\emptyset$.
Then $|\epsilon_n x_n|\rightarrow+\infty $ or $\epsilon_n x_n\rightarrow x_0\in \mathbb{R}^2\backslash\Lambda_0$, and further we assume that $x_0\not\in \partial\Lambda_0$ since we can argue as the following Case 2 if $x_0\in \partial\Lambda_0$). Observe that $V\in L^\infty(\mathbb{R}^2)$, we may assume $a=\lim_{n\rightarrow\infty}V(\epsilon_n x_n)$. Since $V$ is locally H\"{o}lder continuous, we deduce
$z_0=(u_0,v_0)$ satisfies
\begin{equation*}
-\Delta u_0+au_0=\tilde{g}(v_0),\
-\Delta v_0+av_0=\tilde{f}(u_0)\quad \text{in}\quad \mathbb{R}^2.
\end{equation*}
Then $u_0=v_0=0$, contradicts $(u_0,v_0)\neq(0,0)$.

{\bf Case 2}: $\int_{B_1(x_n)}\bigl[\chi_{\Lambda_0}(\epsilon_n x)u^2_n+\chi_{\Lambda_0}(\epsilon_n x)v^2_n\bigr]\geq\frac{\delta}{2}$. If this case occurs, we have
$B_1(x_n)\cap supp\chi_{\Lambda_0}(\epsilon_n x) \neq\emptyset$. Then we may assume $x_n\in \Lambda^{\epsilon_n}_0$ and $\epsilon_n x_n\rightarrow x_0\in \bar{\Lambda}_0$ as $n\rightarrow\infty$. Since $\chi_{\Lambda_0}(\epsilon_n x+\epsilon_n x_n)$ is bounded, we assume $\chi_{\Lambda_0}(\epsilon_n x+\epsilon_n x_n)\rightharpoonup \chi_\infty$ with $\chi_\infty\in [0,1]$ in $L^p_{loc}(\mathbb{R}^2)$ for some $p>2$. Then
$z_0=(u_0,v_0)$ satisfies
\begin{equation}\label{4.3}
-\Delta u_0+V(x_0)u_0={{g}}_\infty(x,v_0), \
-\Delta v_0+V(x_0)v_0={{f}}_\infty(x,u_0)\quad  \text{in}\quad \mathbb{R}^2,
\end{equation}
where ${g}_\infty(x,s)=\chi_\infty(x)g(s)+(1-\chi_\infty(x))\tilde{g}(s)$, ${f}_\infty(x,s)=\chi_\infty(x)f(s)+(1-\chi_\infty(x))\tilde{f}(s)$.
Denote  ${G}_\infty(x,s):=\int^s_0 {g}_\infty(x,t)dt$, ${F}_\infty(x,s):=\int^s_0 {f}_\infty(x,t)dt$ and the functional
$${\Phi}_{\infty}({u},{v}):=\int_{\mathbb{R}^2}(\nabla u\nabla v+V(x_0)uv)-\int_{\mathbb{R}^2}({F}_\infty(x,u)+{G}_\infty(x,v)).$$
Define as (\ref{2.4}) and (\ref{3.4}), we call $(J_\infty, h_\infty)$ the reduced couple of $\Phi_\infty$. Furthermore, as Lemma \ref{l3.7} we can show that, if $u\in H^1(\mathbb{R}^2)\backslash\{0\}$ satisfies $\langle J''_\infty(u)u,u\rangle<0$. Since $(u_0,v_0)\neq(0,0)$ solves (\ref{4.3}), we infer that $u_0+v_0\neq0$ and $\frac{u_0+v_0}{2}$ is a critical point of $J_\infty$. Moreover, $J_\infty(\frac{u_0+v_0}{2})=\max_{t>0}J_\infty(t\frac{u_0+v_0}{2})$.
Therefore,
 $${\Phi}_{\infty}({u}_0,{v}_0)=
 \sup\Bigl\{{\Phi}_{\infty}
 \bigl((t(u_0+v_0)/2,t(u_0+v_0)/2)+(\phi,-\phi)\bigr): t\geq0, \ \phi\in H^1(\mathbb{R}^2)\Bigr\}.$$
 Note that ${F}_\infty(x,s)\leq F(s)$, ${G}_\infty(x,s)\leq G(s)$, and $V_0\leq V(x_0)$, from Lemma \ref{l5.0}, Lemma \ref{l3.2.0} and Lemma \ref{l3.5.0} we conclude
\begin{equation}\label{4.17}\aligned c_{V_0}&\leq c_{V(x_0)}\leq
 \sup\bigl\{{\Phi}_{V(x_0)}
 \bigl((t(u_0+v_0)/2,t(u_0+v_0)/2)+
 (\phi,-\phi)\bigr): t\geq0,\ \phi\in H^1(\mathbb{R}^2)\bigr\}
 \\&\leq\sup\bigl\{{\Phi}_{\infty}
 \bigl((t(u_0+v_0)/2,t(u_0+v_0)/2)
 +(\phi,-\phi)\bigr): t\geq0,\ \phi\in H^1(\mathbb{R}^2)\bigr\}\\&=
{\Phi}_{\infty}({u}_0,{v}_0).
 \endaligned\end{equation}
On the other hand, from Fatou Lemma it follows that
\begin{equation}\label{4.2.3}\aligned c_{\epsilon_n}&={\Phi}_{\epsilon_n}(z_n)-\frac12\langle {\Phi}'_{\epsilon_n}(z_n),(v_n,u_n)\rangle=\int_{\mathbb{R}^2}(\hat{F}(\epsilon_n x,u_n)+\hat{G}(\epsilon_n x,v_n))\\&=\int_{\mathbb{R}^2}\bigl[\hat{F}(\epsilon_n x+\epsilon_n x_n,u_n(x+x_n))+\hat{G}(\epsilon_n x+\epsilon_n x_n,v_n(x+x_n))\bigr]\\
&\geq \int_{\mathbb{R}^2}[\hat{F}_\infty( x,u_0)+\hat{G}_\infty( x,v_0)]+o_n(1)={\Phi}_{\infty}(u_0,v_0)+o_n(1),
\endaligned\end{equation}
where $\hat{F}_\infty( x,s)=\frac12{f}_\infty(x,s)s-{F}_\infty(x,s)$ and
$\hat{G}_\infty( x,s)=\frac12{g}_\infty(x,s)s-{G}_\infty(x,s)$ for all $(x,s)\in\mathbb{R}^2\times \mathbb{R}$. In view of (\ref{4.2.3}), (\ref{4.17}) and Lemma \ref{l4.6}-(2), we infer $V(x_0)=V_0$, $x_0\in \Lambda_0$ and $\chi_\infty\equiv1$.

{\bf Claim 2}: $\bar{z}_n\rightarrow z_0$ in $E$.

Recall that $\bar{z}_n(x)=(u_n(x+x_n),v_n(x+x_n))$ satisfies system (\ref{4.2}), whose corresponding functional  is denoted as,  for all $(u,v)\in E$,
$$\aligned\tilde{\Phi}_{\epsilon_n}(u,v)
=\int_{\mathbb{R}^2}(\nabla u\nabla v+{\tilde{V}}_{\epsilon_n}(x)uv)
-\int_{\mathbb{R}^2}[\bar{F}(\epsilon_nx+\epsilon_nx_n, u)+\bar{G}(\epsilon_nx+\epsilon_nx_n, v)],\endaligned$$
where $\tilde{V}_{\epsilon_n}(x)=V(\epsilon_n x+\epsilon_nx_n)$.
Moreover, $z_0=(u_0,v_0)$ satisfies system (\ref{1.1.1})
and
\begin{equation}\label{4.11} \Phi_{V_0}(z_0)=c_{V_0}
=c_{\epsilon_n}+o_n(1)={\Phi}_{\epsilon_n}(z_{n})+o_n(1).
\end{equation}
Argue by contradiction we assume $\bar{z}_n\not\rightarrow z_0$ in $E$.  Set $w_{n,1}(x)=z_n(x)-{z_0}(x-x_n)$.
Then $w_n(x)=w_{n,1}(x+x_n)\rightharpoonup(0,0)$ and $w_{n}\not\rightarrow(0,0)$ in $E$. For any $\varphi\in E$, by standard arguments we have
$$\tilde{{\Phi}}_{\epsilon_n}(w_n)=\tilde{{\Phi}}_{\epsilon_n}(\bar{z}_n)
-\tilde{{\Phi}}_{\epsilon_n}(z_0)+o_n(1).$$
$$\langle\tilde{{\Phi}}'_{\epsilon_n}(w_n),\varphi\rangle
=\langle\tilde{{\Phi}}'_{\epsilon_n}(\bar{z}_n),\varphi\rangle
-\langle\tilde{{\Phi}}'_{\epsilon_n}(z_0),\varphi\rangle+o_n(1)\|\varphi\|_{1,\epsilon_n}.$$
Using the fact that  $\tilde{V}_{\epsilon_n}(x)\rightarrow V_0$, $\bar{F}(\epsilon_nx+\epsilon_nx_n,s)\rightarrow \bar{F}(x_0,s)=F(s)$ as $n\rightarrow\infty$ uniformly
on any bounded set of $x$, one easily has
$$\langle\tilde{\Phi}'_{\epsilon_n}(z_0),\varphi\rangle=
\langle\Phi'_{V_0}(z_0),\varphi\rangle+o_n(1)\|\varphi\|_{1,\epsilon_n}.$$
Then
\begin{equation}\label{4.12} \aligned\langle \Phi'_{\epsilon_n}(w_{n,1}),\varphi\rangle
&=\langle \tilde{\Phi}'_{\epsilon_n}(\bar{z}_{n}),\varphi(x+x_n)\rangle-\langle \tilde{\Phi}'_{\epsilon_n}(z_0),\varphi(x+x_n)\rangle+o_n(1)
\|\varphi\|_{1,\epsilon_n}\\
&=0-\langle \Phi'_{V_0}(z_0),\varphi(x+x_n)\rangle+o_n(1)\|\varphi\|_{1,\epsilon_n}
=o_n(1)\|\varphi\|_{1,\epsilon_n}.
\endaligned\end{equation}
Similarly, by (\ref{4.11}) we get
\begin{equation}\label{4.13}{\Phi}_{\epsilon_n}(w_{n,1})
={\Phi}_{\epsilon_n}({z}_{n})
-\Phi_{V_0}(z_0)+o_n(1)=o_n(1).\end{equation}
Then $\{w_{n,1}\}$ is nonvanishing, otherwise, as (\ref{5.3}) we easily infer
$w_{n,1}\rightarrow(0,0)$ in $E$, contradicts with the hypothesis that $w_n\not\rightarrow(0,0)$ in $E$. Hence, there exist $\{x^2_n\}\subset\mathbb{R}^2$ and $\delta>0$ such that
\begin{equation}\label{4.15}\int_{B_1(x^2_n)}[(w^1_{n,1})^2+(w^2_{n,1})^2]
\geq\delta.\end{equation}
Set $w_{n,2}=w_{n,1}(x+x^2_n)$ and assume $w_{n,2}=(w^1_{n,2},w^2_{n,2})\rightharpoonup w_2=(w^1_2,w^2_2)\neq(0,0)$ in $E$. Then $w_{n,2}$ satisfies
\begin{equation}\label{4.14}\tilde{\tilde{\Phi}}'_{\epsilon_n}(w_{n,2})\rightarrow0,\end{equation}
with energy functional
$$\aligned\tilde{\tilde{\Phi}}_{\epsilon_n}(u,v)
=&\int_{\mathbb{R}^2}(\nabla u\nabla v+\tilde{\tilde{V}}_{\epsilon_n}(x)uv)
-\int_{\mathbb{R}^2}[\bar{F}(\epsilon_nx+\epsilon_nx^2_n, u)+\bar{G}(\epsilon_nx+\epsilon_nx^2_n, v)],\endaligned$$
for all  $(u,v)\in E$,
where $\tilde{\tilde{V}}_{\epsilon_n}(x)=V(\epsilon_n x+\epsilon_nx^2_n)$. According to (\ref{4.15}), we next discuss for two cases.

{\bf Case 1}:  $\int_{B_1(x^2_n)}\bigl[(1-\chi_{\Lambda_0}(\epsilon_n x))(w^1_{n,1})^2+(1-\chi_{\Lambda_0}(\epsilon_n x))(w^2_{n,1})^2\bigr]\geq{\delta}/{2}$. As the arguments of case 1 in Claim 1, there will be a contradiction.

{\bf Case 2}:  $\int_{B_1(x^2_n)}\bigl[\chi_{\Lambda_0}(\epsilon_n x)(w^1_{n,1})^2+\chi_{\Lambda_0}(\epsilon_n x)(w^2_{n,1})^2\bigr]\geq{\delta}/{2}$.
If this case occurs, we may assume $\epsilon_nx^2_n\rightarrow x'_0\in \bar{\Lambda}_0$. Since $\chi_{\Lambda_0}(\epsilon x+\epsilon_n x^2_n)$ is bounded, we assume $\chi_{\Lambda_0}(\epsilon_n x+\epsilon_n x^2_n)\rightharpoonup \chi_{1,\infty}$ with $\chi_{1,\infty}\in [0,1]$ in $L^p_{loc}(\mathbb{R}^2)$ for some $p>2$.
Using (\ref{4.14}) we deduce that
$w_2=(w^1_2,w^2_2)$ satisfies
\begin{equation}\label{4.16}\aligned
\left\{ \begin{array}{lll}
-\Delta w^1_2+V(x'_0)w^1_2={g}_{1,\infty}(x,w^2_2)\ & \text{in}\quad \mathbb{R}^2,\\
-\Delta w^2_2+V(x'_0)w^2_2={f}_{1,\infty}(x,w^1_2)\ & \text{in}\quad \mathbb{R}^2,
\end{array}\right.\endaligned
\end{equation}
where ${g}_{1,\infty}(x,s)=\chi_{1,\infty} g(s)+(1-\chi_{1,\infty})\tilde{g}(s)$, ${f}_{1,\infty}(x,s)=\chi_{1,\infty} f(s)+(1-\chi_{1,\infty})\tilde{f}(s)$.
Denote
$\Phi_{1,\infty}$ as the associate energy functional of (\ref{4.16})
$$\Phi_{1,\infty}(u,v)=\int_{\mathbb{R}^2}(\nabla u\nabla v+V(x'_0)uv)-\int_{\mathbb{R}^2}[F_{1,\infty}(x,u)+
G_{1,\infty}(x,v)],\quad\forall (u,v)\in E,$$
where $F_{1,\infty}(x,s)=\int^s_0f_{1,\infty}(x,t)dt$, $G_{1,\infty}(x,s)=\int^s_0g_{1,\infty}(x,t)dt$.
By (\ref{4.12}) and (\ref{4.13}), as (\ref{4.2.3}) we deduce
\begin{equation*} o_n(1)={\Phi}_{\epsilon_n}(w_{n,1})-\frac12\langle {\Phi}'_{\epsilon_n}(w_{n,1}),w_{n,1}\rangle\geq\Phi_{1,\infty}(w_2).
\end{equation*}
However, in the same way as (\ref{4.17}) we get
$c_{V_0}\leq c_{V(x'_0)}\leq \Phi_{1,\infty}(w_2)$. This is a contradiction. So $\bar{z}_n\rightarrow z_0$ in $E$.\ \ \ \ \ \ $\Box$

For solutions  $z_n=(u_n,v_n)$ of (\ref{2.3}) obtained in Lemma \ref{l3.5.0} with
$\epsilon_n\rightarrow0$, and $x_n$  obtained in Lemma \ref{l4.7}, we have the following results.
\begin{lemma}\label{l4.8}
$u_n(x+x_n)\rightarrow0$ and $v_n(x+x_n)\rightarrow0$ uniformly in $n$ as $|x|\rightarrow\infty.$
In addition
$\sup_{n\geq1}(|u_n|_\infty+|v_n|_\infty)<+\infty.$
\end{lemma}
{\bf Proof}: Let $\bar{u}_n(\cdot)=u_n(\cdot+x_n)$ and $\bar{v}_n(\cdot)=v_n(\cdot+x_n)$. Note that $\bar{u}_n$ is a weak solution of the following problem
\begin{equation*}-\Delta U+V(\epsilon_nx+\epsilon_nx_n)U=\bar{g}(\epsilon_nx+\epsilon_nx_n, \bar{v}_n) \ \text{in}\ B_2,\ U-\bar{u}_n\in H^1_0(B_2),\end{equation*}
where $B_2=B_2(0)$. Moreover, letting $B_1=B_1(0)$, for any $p\geq2$, we have
\begin{equation}\label{4.0.1}\aligned
\|\bar{u}_n\|_{W^{2,p}(B_1)}\leq C(|g( \bar{v}_n)|_{p,B_2}+|\bar{u}_n|_{p,B_2}).\endaligned\end{equation}
 By the Sobolev embedding theorem, if $p>2$, we get $\bar{u}_n\in C^{1,\gamma}(\bar{B}_1)$ for some $\gamma\in (0,1)$ and there exists $C$ (independent of $n$) such that
 $\|\bar{u}_n\|_{C^{1,\gamma}(\bar{B}_1)}\leq C\|\bar{u}_n\|_{W^{2,p}(B_1)}.$
Using (\ref{4.0.1})  we get
$\|\bar{u}_n\|_{C^{1,\gamma}(\bar{B}_1)}\leq C(|g( \bar{v}_n)|_{p}+|\bar{u}_n|_{p}).$
Assume $\|s\bar{v}_n+(1-s)v_0\|^2_{\epsilon_n}\leq M$, where $s\in (0,1)$.
For some $\alpha\in (0,\frac{2\pi}{pM})$, the Trudinger-Moser inequality implies that
\begin{equation}\label{4.22.0}\int_{\mathbb{R}^2}\Bigl[e^{2\alpha p(s\bar{v}_n+(1-s)v_0)^2}-1\Bigr]\leq
\int_{\mathbb{R}^2}\bigl[e^{2\alpha pM\frac{(s\bar{v}_n+(1-s)v_0)^2}{\|s\bar{v}_n+(1-s)v_0\|^2_{\epsilon_n}}}-1\bigr]\leq C.\end{equation}
As (\ref{2.2.0}), for the above $\alpha$, for any $\delta>0$ there exists $C_\delta>0$ such that
\begin{equation}\label{4.2.2}\aligned
&\int_{\mathbb{R}^2}|g(\bar{v}_n)-g(v_0)|^p
=\int_{\mathbb{R}^2}|g'(s\bar{v}_n+(1-s)v_0)|^p|\bar{v}_n-v_0|^p\\
\leq&\delta|\bar{v}_n-v_0|^p_p+C_\delta\Bigl(\int_{\mathbb{R}^2}[e^{2\alpha p(s\bar{v}_n+(1-s)v_0)^2}-1]\Bigr)^{\frac12}|\bar{v}_n-v_0|^{p}_{2p} \rightarrow0.\endaligned\end{equation}Hence
\begin{equation}\label{4.1.4}
\sup_{n\geq1}\|\bar{u}_n\|_{C^{1,\gamma}(\bar{B}_1)}<+\infty.
\end{equation}
Below as the proof of \cite[Proposition 2.8]{ZhangCPDE}, we get $\bar{u}_n(x)\rightarrow0$ uniformly as $|x|\rightarrow\infty$.
Combining with (\ref{4.1.4}) we infer that
$\sup_{n\geq1}|u_n|_\infty=\sup_{n\geq1}|\bar{u}_n|_\infty<+\infty$. Similarly, $\sup_{n\geq1}|v_n|_\infty<+\infty$.\ \ \ \ \ $\Box$

As \cite[Proposition 2.7]{ZhangCPDE}, we have the following lemma.
\begin{lemma}\label{l4.9} There exists $\delta>0$ independent of $n$ such that $\min\{|u_n|_\infty,|v_n|_\infty\}\geq \delta$.\end{lemma}

By Lemma \ref{l4.8}, we may assume $\{y_n\}\subset\mathbb{R}^2$ satisfies
\begin{equation}\label{4.29}|u_n(y_n)|+|v_n(y_n)|=
\max_{\mathbb{R}^2}(|u_n(x)|+|v_n(x)|).\end{equation} Then the following result holds.
\begin{lemma}\label{l4.10} $z_n(\cdot+y_n)\rightarrow \hat{z}_0$ in $E$ as $n\rightarrow\infty$, where $\hat{z}_0$ is a ground state of the limit equation (\ref{1.1.1}). Moreover, $z_n(x+y_n)\rightarrow(0,0)$ uniformly in $n$ as $|x|\rightarrow+\infty$.
\end{lemma}
{\bf Proof}: Firstly we claim that there exist $\mu, R_1>0$ such that
\begin{equation}\label{4.2.0} \lim_{n\rightarrow\infty}\int_{B_{R_1}(y_n)}(u^2_n+v^2_n)\geq\mu.\end{equation}
Argue by contradiction we assume that, for any $R>0$,
$\int_{B_{R}(y_n)}(u^2_n+v^2_n)\rightarrow0$ as $n\rightarrow\infty$.
Let $\hat{u}_n(\cdot)=u_n(\cdot+y_n)$ and $\hat{v}_n(\cdot)=v_n(\cdot+y_n)$, then $\hat{u}_n, \hat{v}_n\rightarrow0$ in $L^2_{loc}(\mathbb{R}^2)$ as $n\rightarrow\infty$. Observe that $\hat{u}_n$ is a weak solution of the following problem
\begin{equation*}-\Delta U+V(\epsilon_nx+\epsilon_ny_n)U=\bar{g}(\epsilon_nx+\epsilon_ny_n, \hat{{v}}_n) \ \text{in}\ B_2,\ U-\hat{u}_n\in H^1_0(B_2).\end{equation*}
As in Lemma \ref{l4.8}, we get $\hat{u}_n\in C^{1,\gamma}(\bar{B}_1)$ for some $\gamma\in (0,1)$ and $\|\hat{u}_n\|_{C^{1,\gamma}(\bar{B}_1)}\leq C(|g( \hat{v}_n)|_{p}+|\hat{u}_n|_{p}),$ where $C$ is independent of $n$. Moreover, by (\ref{4.2.2}) we get
$|g(\hat{v}_n)|_{p}\rightarrow|g(v_0)|_p$. Then
$ \sup_{n\geq1}\|\hat{u}_n\|_{C^{1,\gamma}(\bar{B}_1)}<+\infty.$
Since $\hat{u}_n\rightarrow0$ in $L^2(B_1)$, we get ${\hat{u}_n}\rightarrow0$ uniformly in $B_1$. In particular, $\hat{u}_n(0)=u_n(y_n)\rightarrow0$. Similarly, $\hat{v}_n(0)=v_n(y_n)\rightarrow0$. Then by (\ref{4.29}) we obtain
$\max_{x\in\mathbb{R}^2}(|u_n(x)|+|v_n(x)|)
\rightarrow0,$
contradicts with Lemma \ref{l4.9}. Thus (\ref{4.2.0}) holds true.
Then we may assume $z_n(\cdot+y_n)\rightharpoonup \hat{z}_0\neq(0,0)$. Arguing as Claims 1 and 2 in the proof of Lemma \ref{l4.7},  we deduce that $\epsilon_ny_n\rightarrow y_0\in \Lambda_0$, $V(y_0)=V_0$ and $z_n(\cdot+y_n)\rightarrow \hat{z}_0$ in $E$, where $\hat{z}_0$ is a ground state of the limit system (\ref{1.1.1}). Below taking similar arguments as \cite[Proposition 2.8]{ZhangCPDE}, we conclude $z_n(x+y_n)\rightarrow(0,0)$ uniformly in $n$ as $|x|\rightarrow+\infty$.\ \ \ \ \ $\Box$

{\bf Proof of Theorem \ref{t1.1}}
Similar to the arguments of \cite[page 22]{Zhang-du}, we show that the solutions $z_n=(u_n,v_n)$ of the penalized problem (\ref{2.3}) obtained in Lemma \ref{l3.5.0} are actually the solutions of the problem (\ref{2.1}). Combining with Lemmas \ref{l4.7}-\ref{l4.10}, the conclusions of Theorem \ref{t1.1} (1), (2)-(i) and (2)-(ii) hold true.

It suffices to show  Theorem \ref{t1.1} (2)-(iii).
Note that
\begin{equation}\label{bu6}-\Delta \hat{{u}}_n+V(\epsilon_n x+\epsilon_ny_n)\hat{u}_n=g(\hat{v}_n),\ -\Delta \hat{v}_n+V(\epsilon_n x+\epsilon_ny_n)\hat{v}_n=f(\hat{u}_n) \ \text{in}\ \mathbb{R}^2,\end{equation}
where $\hat{u}_n(\cdot)=u_n(\cdot+y_n)$ and $\hat{v}_n(\cdot)=v_n(\cdot+y_n)$.
Testing the first equality of (\ref{bu6}) with $sgn(u_n)$, where $sgn$ is the sign function,
we obtain
$$sgn(\hat{u}_n)[-\Delta \hat{u}_n+V(\epsilon_n x+\epsilon_n y_n) \hat{u}_n]=sgn(\hat{u}_n)g(\hat{v}_n).$$
By \cite[Lemma A]{Kato}, we get $-sgn(\hat{u}_n)\Delta \hat{u}_n\geq-\Delta|\hat{u}_n|$ in the sense of distribution. Thus by Lemma \ref{l4.10} and (H$_1$), for some $R_0>0$ there holds \begin{equation}\label{5.0.1}-\Delta |\hat{u}_n|+V(\epsilon_n x+\epsilon_n y_n)|\hat{u}_n|\leq|g(\hat{v}_n)|\leq\delta |\hat{v}_n|,\quad \text{if}\ |x|\geq R_0.\end{equation}
Similarly
$-\Delta |\hat{v}_n|+V(\epsilon_n x+\epsilon_n y_n)|\hat{v}_n|\leq\delta |\hat{u}_n|$, if $|x|\geq R_0.$
Letting $\delta\leq\frac{\inf_{\mathbb{R}^2}V}{2}$,
there holds
\begin{equation}\label{5.0.2}-\Delta (|\hat{u}_n|+|\hat{v}_n|)+\frac{\inf_{\mathbb{R}^2}V}{2}(|\hat{u}_n|+|\hat{v}_n|)\leq0 \ \text{if} \ |x|\geq R_0.\end{equation} Then as \cite[Page 23]{Zhang-du}, for some $C,c>0$ there holds $|{u}_{n}(x)|+|{v}_{n}(x)|\leq C e^{-c|x-y_n|}$.\ \ \ \ \ $\Box$

{\bf Proof of Theorems \ref{t1.2} and \ref{t1.3}}  By Theorem \ref{t1.1}, we assume $z_n=(u_n,v_n)$ are nontrivial solutions of (\ref{2.1}) with $\epsilon_n\rightarrow0$, $x^1_n$ and $x^2_n$ are maximum points of $u_n$ and $v_n$ respectively. Using Lemma \ref{lbu} and taking similar arguments as in the proof of Propositions 3.13 and 3.14 in \cite{ZhangCPDE}, we can deduce that for $\epsilon_n>0$ small enough, $u_n>0$, $v_n>0$ in $\mathbb{R}^2$,
 $\lim_{n\rightarrow\infty}|x^1_n-x^2_n|=0,$
and $x^1_n$ and $x^2_n$ are unique.
Letting $(\varphi_\epsilon(x),\psi_\epsilon(x))=(u(\frac{x}{\epsilon}),
v(\frac{x}{\epsilon}))$,  the conclusions of Theorems \ref{t1.2} and \ref{t1.3} yield. \ \ \ \ \ \ $\Box$

\section{Proof of Theorems \ref{t1.4} and \ref{t1.5}}
\renewcommand{\theequation}{5.\arabic{equation}}
In this section, we prove Theorems \ref{t1.4} and \ref{t1.5}.
\subsection{Refined Nehari manifold}
In the following we fix bounded domains ${\Lambda}'_i,\ \tilde{\Lambda}_i$, mutually disjoint, such that $\Lambda_i\Subset{\Lambda}'_i\Subset \tilde{\Lambda}_i\Subset\mathbb{R}^2$ for each $i=1,...,k$. We take cut-off functions $\phi_i$ such that $\phi_i=1$ in $\Lambda'_i$ and $\phi_i=0$ in $\mathbb{R}^2\backslash{\tilde{\Lambda}_i}$, and suppose without loss of generality that $|\nabla \phi_i|\leq C\phi^{\frac12}_i$ for some $C>0$ (by replacing if necessary $\phi_i$ by $\phi^2_i$). We also denote $\Lambda:=\cup^k_{i=1}\Lambda_i$, $\Lambda':=\cup^k_{i=1}\Lambda'_i$ and $\tilde{\Lambda}:=\cup^k_{i=1}\tilde{\Lambda}_i$.
For $\tilde{f}$ and $\tilde{g}$ defined in Section 2.1,  we introduce another modified nonlinearity as follows.
For any $(x,t)\in\mathbb{R}^2\times\mathbb{R}$, define
 $$\bar{f}_0(x,t)=\chi_\Lambda(x)f(t)+(1-\chi_\Lambda(x))\tilde{f}(t), \ \bar{g}_0(x,t)=\chi_\Lambda(x)g(t)+(1-\chi_\Lambda(x))\tilde{g}(t),$$
 and denote $\bar{F}_0(x,t)=\int^t_0\bar{f}_0(x,s)ds$, $\bar{G}_0(x,t)=\int^t_0\bar{g}_0(x,s)ds$.
 It is easy to see that, for the above $\bar{f}_0$ and $\bar{g}_0$, the properties as in Lemma \ref{l2.1} with $\Lambda_0$ replacing by $\Lambda$ still hold true. Moreover, using (H$_7$),
 for every $\mu>0$, there exists $C_\mu>0$ such that
 \begin{equation}\label{5.19}
 |\bar{f}_0(x,s)t|+|\bar{g}_0(x,t)s|\leq\mu(s^2+t^2)+C_\mu(
 \bar{f}_0(x,s)s+\bar{g}_0(x,t)t), \forall x, (s,t)\in\mathbb{R}^2.\end{equation}

Now we establish the modified problem
\begin{equation}\label{7.0}\aligned
\left\{ \begin{array}{lll}
-\Delta u+V(\epsilon x)u=\bar{g}_0(\epsilon x, v)\ & \text{in}\quad \mathbb{R}^2,\\
-\Delta v+V(\epsilon x)v=\bar{f}_0(\epsilon x, u)\ & \text{in}\quad \mathbb{R}^2.
\end{array}\right.\endaligned
\end{equation}
The functional of (\ref{2.3}) is
$${\Psi}_\epsilon(z)=\int_{\mathbb{R}^2}(\nabla u\nabla v+V(\epsilon x)uv)-\int_{\mathbb{R}^2}\bar{G}_0(\epsilon x,v)-\int_{\mathbb{R}^2}\bar{F}_0(\epsilon x,u),\quad \forall z=(u,v)\in E.$$
As (\ref{3.4}) and (\ref{2.4}),  we define $\bar{h}_\epsilon$ and $I_\epsilon$ as the reduced couple of $\Psi_\epsilon$. Denote
$\phi^\epsilon_i(x):=\phi_i(\epsilon x)$ and  $\Lambda^\epsilon_i=\Lambda_i/\epsilon$. We also define refined Nehari manifold $M_\epsilon$ by
$$M_\epsilon=\Bigl\{u\in H^1(\mathbb{R}^2): \langle I'_\epsilon(u),u\phi^\epsilon_i\rangle=0,\ \int_{{\Lambda^\epsilon_i}}[u^2+\bar{h}^2_\epsilon(u)]>\epsilon^{{
\frac12}},\ i=1,...,k\Bigr\},$$
and the least energy on $M_\epsilon$ by
$$d_\epsilon=\inf_{M_\epsilon}I_\epsilon.$$
 The proof of non-emptiness of $M_\epsilon$ is complicated, and is postponed in Lemma \ref{p5.15} below. Fix points $x_i\in \Lambda_i$ such that $V(x_i)=\inf_{\Lambda_i}V=V_i$ and consider a fixed nontrivial solution $(u_i,v_i)\in E$ of system (\ref{bu2}),
corresponding to the least energy
$$c_i:=c_{V(x_i)}=\Phi_{V(x_i)}(u_i,v_i).$$
Let
\begin{equation}\label{7.5}u_{i,\epsilon}(x):=\phi^\epsilon_i(x)u_i(x-{x_i}/\epsilon),\quad v_{i,\epsilon}(x):=\phi^\epsilon_i(x)v_i(x-{x_i}/\epsilon),\ w_{i,\epsilon}=(u_{i,\epsilon}+v_{i,\epsilon})/2.\end{equation}
For simplicity, denote \begin{equation}\label{5.6}\phi^\epsilon_{x_i}:=\phi^\epsilon_{i}(x+x_i/\epsilon)=\phi_i(\epsilon x+ x_i),\ A^{\epsilon}_{x_i}:=(A-x_i)/\epsilon\ \text{for any set}\ A\subseteq\mathbb{R}^2.\end{equation}
\begin{lemma}\label{l5.2}(\cite[Lemma 3.3]{Cassani2}) Assume $\Omega\subset\mathbb{R}^2$ is a bounded domain with smooth boundary. Assume $\{u_n\}\subset H^1_0(\Omega)\backslash\{0\}$ is such that $|\nabla u_n|_2\rightarrow0$ as $n\rightarrow\infty$. Then
$$\liminf_{n\rightarrow\infty}\frac{\int_{\Omega}|\nabla u_n|^2}{\Bigl[\int_{\Omega}u^2_n(e^{u^2_n}-1)\Bigr]^{\frac12}}\geq S_4,$$
where $S_4$ is the best Sobolev constant of the embedding $H^1_0(\Omega)\hookrightarrow L^4(\Omega)$.
\end{lemma}
\begin{lemma}\label{p5.15}For every $C_0>0$, there exist $\epsilon_0$, $D_0$, $\eta_0>0$ such that, for every $\epsilon\in(0,\epsilon_0]$, if $w_\epsilon\in M_\epsilon, I_\epsilon(w_\epsilon)\leq C_0$, then
$\|w_\epsilon\|^2+\|\bar{h}_\epsilon(w_\epsilon)\|^2\leq D_0$ and $\int_{\Lambda^\epsilon_i}(w^2_\epsilon+\bar{h}^2_\epsilon(w_\epsilon))
\geq\eta_0,$
for every $i=1,...,k$. Moreover, $d_\epsilon\geq\rho_0$ for some $\rho_0>0$ if $\epsilon\in(0,\epsilon_0]$.
\end{lemma}
{\bf Proof}: For any $w_\epsilon\in M_\epsilon$ satisfying $I_\epsilon(w_\epsilon)\leq C_0$, let $u_\epsilon=w_\epsilon+\bar{h}_\epsilon(w_\epsilon)$ and $v_\epsilon=w_\epsilon-\bar{h}_\epsilon(w_\epsilon)$. Then $(u_\epsilon(\frac x\epsilon),v_\epsilon(\frac x\epsilon))\in \tilde{M}_\epsilon$, where
\begin{equation}\label{5.0.11}\aligned\tilde{M}_\epsilon=\Bigl\{(u,v)\in E:&\ \langle\Psi'_\epsilon(u,v),(\phi,-\phi)\rangle=0,\ \forall \phi\in H^1(\mathbb{R}^2), \\& \langle\Psi'_\epsilon(u,v),(u\phi_i,v\phi_i)\rangle=0, \int_{\Lambda_i}(u^2+v^2)>\epsilon^{N+1}, \forall i=1,...,k\Bigr\}.\endaligned\end{equation}
In the same way as Step 1 and 3 in the proof of \cite[Proposition 5.15]{HUGO-PHD}, we can deduce
$(u_\epsilon,v_\epsilon)$ is bounded in $E$. So $\|w_\epsilon\|^2+\|\bar{h}_\epsilon(w_\epsilon)\|^2\leq D_0$. It suffices to show $\int_{\Lambda^\epsilon_i}[w^2_\epsilon+\bar{h}^2_\epsilon
(w_\epsilon)]\geq\eta_0$ and $d_\epsilon\geq\rho_0$.
Since $f$ and $g$ are of exponential growth, there are some  differences from the argument of \cite[Proposition 5.15]{HUGO-PHD}, where
the authors considered system (\ref{1.5}) which is of polynomial growth.

{\bf Claim 1}. There is $\eta>0$ such that
\begin{equation*}
\int_{\mathbb{R}^2}(|\nabla u_\epsilon|^2+|\nabla v_\epsilon|^2)\phi^\epsilon_i\geq\eta\quad\text{for small} \ \epsilon.
\end{equation*}
Otherwise, assume $\int_{\mathbb{R}^2}(|\nabla u_\epsilon|^2+|\nabla v_\epsilon|^2)\phi^\epsilon_i\rightarrow0$ as $\epsilon\rightarrow0$. Since $(u_\epsilon,v_\epsilon)$ is bounded in $E$, we get
\begin{equation}\label{5.36.1}\aligned
&\int_{\mathbb{R}^2}\bigl(|\nabla (u_\epsilon\phi^\epsilon_i)|^2+|\nabla (v_\epsilon\phi^\epsilon_i)|^2\bigr)\\
=&\int_{\mathbb{R}^2}\Bigl[(|\nabla u_\epsilon|^2+|\nabla v_\epsilon|^2)(\phi^\epsilon_i)^2+2\epsilon (u_\epsilon\nabla u_\epsilon+v_\epsilon\nabla v_\epsilon)\phi^\epsilon_i\nabla \phi^\epsilon_i+\epsilon^2(u^2_\epsilon+v^2_\epsilon)(\nabla \phi^\epsilon_i)^2\Bigr]=o_\epsilon(1).\endaligned\end{equation}
Using $w_\epsilon\in M_\epsilon$, we get
$\langle I'_\epsilon(w_\epsilon),w_\epsilon\phi^\epsilon_i\rangle=0.$
Then
$\langle \Psi'_\epsilon(u_\epsilon,v_\epsilon), (v_\epsilon\phi^\epsilon_i,u_\epsilon\phi^\epsilon_i)\rangle=0$.
Hence
\begin{equation}\label{5.36.0}\aligned &\int_{\mathbb{R}^2}(|\nabla u_\epsilon|^2+|\nabla v_\epsilon|^2)\phi^\epsilon_i+\int_{\mathbb{R}^2}V(\epsilon x)(u^2_\epsilon+v^2_\epsilon)\phi^\epsilon_i\\=&\int_{\mathbb{R}^2}[
\bar{f}_0(\epsilon x,u_\epsilon)v_\epsilon\phi^\epsilon_i+\bar{g}_0(\epsilon x,v_\epsilon)u_\epsilon\phi^\epsilon_i]-
\epsilon\int_{\mathbb{R}^2}(u_\epsilon\nabla {u_\epsilon}+v_\epsilon\nabla {v_\epsilon})\nabla\phi^\epsilon_i.\endaligned\end{equation}
Note that $|\nabla \phi_i|\leq C\phi^{\frac12}_i$ we infer
$$\epsilon\Bigl|\int_{\mathbb{R}^2}u_\epsilon\nabla {u_\epsilon}\nabla\phi^\epsilon_i\Bigr|
\leq\epsilon\int_{\mathbb{R}^2}u^2_\epsilon+\epsilon\int_{\mathbb{R}^2}|\nabla u_\epsilon|^2\phi^\epsilon_i.$$
Moreover, since $\int_{\Lambda^\epsilon_i}(u^2_\epsilon+v^2_\epsilon)>\epsilon^{\frac12}$ we get
$$\epsilon\int_{\mathbb{R}^2}u^2_\epsilon\leq D_0\epsilon^{\frac12}\int_{\mathbb{R}^2}(u^2_\epsilon+v^2_\epsilon)\phi^\epsilon_i
\leq \frac14\int_{\mathbb{R}^2}V(\epsilon x)(u^2_\epsilon+v^2_\epsilon)\phi^\epsilon_i,$$
for small enough $\epsilon>0$. Proceeding in a similar way with the function $v_\epsilon$ we have
$$-\epsilon\int_{\mathbb{R}^2}(u_\epsilon\nabla {u_\epsilon}+v_\epsilon\nabla {v_\epsilon})\nabla\phi^\epsilon_i\leq\frac12\int_{\mathbb{R}^2}
\bigl[|\nabla u_\epsilon|^2+|\nabla v_\epsilon|^2+V(\epsilon x)(u^2_\epsilon+v^2_\epsilon)\bigr]\phi^\epsilon_i.$$
Using (\ref{5.36.0}), (\ref{5.19}) and ($H''_6$) we get
\begin{equation*}
\aligned&\int_{\mathbb{R}^2}(|\nabla u_\epsilon|^2+|\nabla v_\epsilon|^2)\phi^\epsilon_i+\int_{\mathbb{R}^2}
V(\epsilon x)(u^2_\epsilon+v^2_\epsilon)\phi^\epsilon_i\\ \leq& \mu\int_{\mathbb{R}^2}
V(\epsilon x)(u^2_\epsilon+v^2_\epsilon)\phi^\epsilon_i+C_{\mu}
\int_{\mathbb{R}^2}[\bar{f}_0(\epsilon x,u_\epsilon)u_\epsilon\phi^\epsilon_i+\bar{g}_0(\epsilon x,v_\epsilon)v_\epsilon\phi^\epsilon_i]\\
 \leq& \mu\int_{\mathbb{R}^2}
V(\epsilon x)(u^2_\epsilon+v^2_\epsilon)\phi^\epsilon_i+C_{\mu}\delta
\int_{(\Lambda^\epsilon)^c}(u^2_\epsilon+v^2_\epsilon)\phi^\epsilon_i
+C_{\mu}
\int_{\Lambda^\epsilon_i}(f(u_\epsilon)u_\epsilon+g(v_\epsilon)v_\epsilon),
\endaligned\end{equation*}
where $\mu$, $\delta>0$ are arbitrary and $C_{\mu}$ depends on $\mu$.
Letting $\mu=\frac14$ and then letting $\delta=\frac{\inf_{\mathbb{R}^2}V}{4C_{1/4}}$, together with (\ref{2.1.0}) with $q=3$ and $\alpha=1$ we infer\begin{equation}\label{5.36.2}\aligned
&\int_{\mathbb{R}^2}(|\nabla u_\epsilon|^2+|\nabla v_\epsilon|^2)\phi^\epsilon_i+\int_{\mathbb{R}^2}
V(\epsilon x)(u^2_\epsilon+v^2_\epsilon)\phi^\epsilon_i
 \leq C\int_{\Lambda^\epsilon_i}(f(u_\epsilon)u_\epsilon
+g(v_\epsilon)v_\epsilon)\\
\leq&\frac12\int_{\mathbb{R}^2}V(\epsilon x)(|u_\epsilon|^2+|v_\epsilon|^2)\phi^\epsilon_i+
C\int_{\Lambda^\epsilon_i}[|u_\epsilon|^{2}(e^{u^2_\epsilon}-1)+
|v_\epsilon|^{2}(e^{v^2_\epsilon}-1)].\endaligned\end{equation}
In view of (\ref{5.36.1}), from Lemma \ref{l5.2} it follows that
\begin{equation}\label{5.36.3}\aligned
&\int_{\Lambda^\epsilon_i}[|u_\epsilon|^{2}(e^{u^2_\epsilon}-1)+
|v_\epsilon|^{2}(e^{v^2_\epsilon}-1)]\leq \int_{\mathbb{R}^2}\bigl[|u_\epsilon\phi^\epsilon_i|^{2}(e^{ u^2_\epsilon(\phi^\epsilon_i)^2}-1)+
|v_\epsilon\phi^\epsilon_i|^{2}(e^{ v^2_\epsilon(\phi^\epsilon_i)^2}-1)\bigr]\\
\leq& 4S^{-2}_4\bigl(\int_{\mathbb{R}^2}|\nabla (u_\epsilon\phi^\epsilon_i)^2|\bigr)^2+4S^{-2}_4\bigl(\int_{\mathbb{R}^2}|\nabla (v_\epsilon\phi^\epsilon_i)^2|\bigr)^2\\
\leq& C'\Bigl(\int_{\mathbb{R}^2}(|\nabla u_\epsilon|^2+|\nabla v_\epsilon|^2)(\phi^\epsilon_i)^2\Bigr)^2+C'\Bigl(\epsilon^2
\int_{{\mathbb{R}}^2}
(|\nabla\phi^\epsilon_i|^2(u^2_\epsilon+v^2_\epsilon)\Bigr)^2\\
\leq& C'\Bigl(\int_{\mathbb{R}^2}(|\nabla u_\epsilon|^2+|\nabla v_\epsilon|^2)\phi^\epsilon_i\Bigr)^2+C''\epsilon^4\Bigl(
\int_{{\mathbb{R}}^2}
V(\epsilon x)(u^2_\epsilon+v^2_\epsilon)\phi^\epsilon_i\Bigr)^2.
\endaligned\end{equation}
Note that  $\|(u_\epsilon,v_\epsilon)\|\leq D_0$. Then
$\epsilon^4
\int_{{\mathbb{R}}^2}
V(\epsilon x)(u^2_\epsilon+v^2_\epsilon)\phi^\epsilon_i\leq\frac14,$
for sufficiently small $\epsilon>0$.
Combining (\ref{5.36.2}) and (\ref{5.36.3}) we conclude
$$\int_{\mathbb{R}^2}(|\nabla u_\epsilon|^2+|\nabla v_\epsilon|^2)\phi^\epsilon_i+\int_{\mathbb{R}^2}
V(\epsilon x)(u^2_\epsilon+v^2_\epsilon)\phi^\epsilon_i\leq
C\bigl(\int_{\mathbb{R}^2}(|\nabla u_\epsilon|^2+|\nabla v_\epsilon|^2)\phi^\epsilon_i\bigr)^2.$$
Then Claim 1 yields. Moreover, by the first inequality of (\ref{5.36.2}) and Claim 1 we get
\begin{equation}\label{5.36.5}\int_{\Lambda^\epsilon_i}
[f(u_\epsilon)u_\epsilon+g(v_\epsilon)v_\epsilon]\geq\eta'>0,\ \ \forall i.\end{equation}

{\bf Claim 2}. There exists $\eta_0>0$ such that
$$\int_{\Lambda^\epsilon_i}(w^2_\epsilon+\bar{h}^2_\epsilon(w_\epsilon))\geq\eta_0, \ \text{for small}\ \epsilon.$$
Otherwise, we assume there exists $w_\epsilon\in M_\epsilon$ such that
$\int_{\Lambda^\epsilon_i}(w^2_\epsilon+\bar{h}^2_\epsilon(w_\epsilon))
\rightarrow0.$ Then
$\int_{\Lambda^\epsilon_i}(u^2_\epsilon+v^2_\epsilon)\rightarrow0.$
Since $\|(u_\epsilon,v_\epsilon)\|\leq D_0$, letting $0<\alpha<\frac{2\pi}{D^2_0}$, by Trudinger-Moser inequality we get
$\int_{\mathbb{R}^2}(e^{2\alpha u^2_\epsilon}-1)\leq C$, \ $\int_{\mathbb{R}^2}(e^{2\alpha v^2_\epsilon}-1)\leq C.$
By the first inequality of (\ref{5.36.2}), (\ref{2.1.0}) and H\"{o}lder inequality, for any $\delta_0>0$ we have
\begin{equation*}\aligned&\int_{\mathbb{R}^2}(|\nabla u_\epsilon|^2+|\nabla v_\epsilon|^2)\phi^\epsilon_i+
\int_{\mathbb{R}^2}
V(\epsilon x)(u^2_\epsilon+v^2_\epsilon)\phi^\epsilon_i\\
\leq&\delta_0(|u_\epsilon|^2_2+|v_\epsilon|^2_2)+C
\Bigl(\int_{\Lambda^\epsilon_i}|u_\epsilon|^2
\int_{\Lambda^\epsilon_i}(e^{2\alpha u^2_\epsilon}-1)\Bigr)^{\frac12}+C
\Bigl(\int_{\Lambda^\epsilon_i}|v_\epsilon|^2
\int_{\Lambda^\epsilon_i}(e^{2\alpha v^2_\epsilon}-1)\Bigr)^{\frac12}\rightarrow0,\endaligned\end{equation*}
contradicts Claim 1 and thus Claim 2 holds.

{\bf Claim 3}: $d_\epsilon\geq\rho_0$ for some $\rho_0>0$.
Argue by contradiction we assume there exists $\tilde{w}_\epsilon\in M_\epsilon$ such that $I_\epsilon(\tilde{w}_\epsilon)=o_\epsilon(1)$. Letting $\tilde{u}_\epsilon=\tilde{w}_\epsilon+\bar{h}_\epsilon(\tilde{w}_\epsilon)$ and $\tilde{v}_\epsilon=\tilde{w}_\epsilon-\bar{h}_\epsilon(\tilde{w}_\epsilon)$, for small $\delta>0$ we get
\begin{equation}\label{5.36.12}
\aligned
(\tilde{u}_\epsilon,\tilde{v}_\epsilon)_\epsilon
&=\int_{\mathbb{R}^2}[\bar{F}_0(\epsilon x,\tilde{u}_\epsilon)+\bar{G}_0(\epsilon x,\tilde{v}_\epsilon)]+\Psi_\epsilon(\tilde{u}_\epsilon,\tilde{v}_\epsilon)\\
&\leq\frac{1}{2+\delta'}\int_{\Lambda^\epsilon}
\bigl(f(\tilde{u}_\epsilon)\tilde{u}_\epsilon+g(\tilde{v}_\epsilon)
\tilde{v}_\epsilon\bigr)+
\delta\int_{\mathbb{R}^2}(\tilde{u}^2_\epsilon+\tilde{v}^2_\epsilon)
+o_\epsilon(1).\endaligned
\end{equation}
Let $\phi^\epsilon:=\Sigma^k_{i=1}\phi^\epsilon_i$. Since
$\langle\Psi'_\epsilon(\tilde{u}_\epsilon,\tilde{v}_\epsilon),(\phi^\epsilon \tilde{v}_\epsilon,\phi^\epsilon \tilde{u}_\epsilon)\rangle=0$, as (5.27), (5.28) and (5.29) in \cite{HUGO-PHD} we deduce
$$\int_{\Lambda^\epsilon}
[f(\tilde{u}_\epsilon)\tilde{u}_\epsilon+g(\tilde{v}_\epsilon)
\tilde{v}_\epsilon]\leq2( \tilde{u}_\epsilon,\tilde{v}_\epsilon)_\epsilon+
\delta(|\tilde{u}_\epsilon|^2_2+|\tilde{v}_\epsilon|^2_2)+o_\epsilon(1)
\|(\tilde{u}_\epsilon,\tilde{v}_\epsilon)\|^2_{1,\epsilon}.$$
Combining with (\ref{5.36.12}) we get
$$(1-\frac{2}{2+\delta'})\int_{\Lambda^\epsilon}
[f(\tilde{u}_\epsilon)\tilde{u}_\epsilon+g(\tilde{v}_\epsilon)\tilde{v}_\epsilon]\leq
\delta\int_{\mathbb{R}^2}(\tilde{u}^2_\epsilon+\tilde{v}^2_\epsilon)+
o_\epsilon(1)\|(\tilde{u}_\epsilon,\tilde{v}_\epsilon)\|^2_\epsilon
+o_\epsilon(1),$$
contradicts (\ref{5.36.5}) since $\delta$ is small enough and $(\tilde{u}_\epsilon,\tilde{v}_\epsilon)$ is bounded in $E$.\ \ \ \ $\Box$

As \cite[Proposition 5.20]{HUGO-PHD}, we have the following result.
\begin{lemma}\label{l5.3} For $w_\epsilon$ given in Lemma \ref{p5.15}, letting $u_\epsilon=w_\epsilon+\bar{h}_\epsilon(w_\epsilon)$, $v_\epsilon=w_\epsilon-\bar{h}_\epsilon(w_\epsilon)$ and $i\in\{1,...,k\}$, there exist $\epsilon_0,\eta>0$ such that for any $0<\epsilon<\epsilon_0$ and any $\psi\in H^1(\mathbb{R}^2)$
$$\aligned&\langle \Psi'_\epsilon(u_\epsilon,v_\epsilon),(u_\epsilon(\phi^\epsilon_i)^2,v_\epsilon(\phi^\epsilon_i)^2)\rangle+
\langle \Psi''_\epsilon(u_\epsilon,v_\epsilon)(u_\epsilon\phi^\epsilon_i+\psi,v_\epsilon\phi^\epsilon_i-\psi),
(u_\epsilon\phi^\epsilon_i+\psi,v_\epsilon\phi^\epsilon_i-\psi)\rangle\\
&\leq-\eta-\frac{\|\psi\|^2_\epsilon}{2}<0.
\endaligned$$
\end{lemma}

\begin{lemma}\label{p5.14} There exists $\epsilon_0>0$ such that for any $0<\epsilon\leq\epsilon_0$ and $w_{i,\epsilon}$ given in (\ref{7.5}), $i=1,...,k$,
there are $t_{1,\epsilon},...,t_{k,\epsilon}\in[0,2]$ such that the function
$\bar{w}_\epsilon=\Sigma^{k}_{i=1}t_{i,\epsilon}
w_{i,\epsilon}$
satisfies
\begin{equation}\label{5.5.0}\langle I'_\epsilon(\bar{w}_\epsilon),\bar{w}_\epsilon\phi^\epsilon_i\rangle=0,\ \forall i=1,...,k,\
I_\epsilon(\bar{w}_\epsilon)=\Sigma^k_{i=1}c_i+o_\epsilon(1),\ \text{as}\ \epsilon\rightarrow0.\end{equation}
Moreover,
\begin{equation}\label{5.5.1}\int_{\Lambda^\epsilon_i}(\bar{w}^2_\epsilon+
\bar{h}^2_\epsilon(\bar{w}_\epsilon))\geq\eta.\end{equation}
for some $\eta>0$. In particular, $\bar{w}_\epsilon\in M_\epsilon$. \end{lemma}
{\bf Proof}: For every $\bar{t}=(t_1,\cdots,t_k)\in[0,2]\times\ \cdots\times[0,2]$, let
$$u_{\epsilon,\bar{t}}:=\Sigma^{k}_{i=1}t_iu_{i,\epsilon},\quad v_{\epsilon,\bar{t}}:=\Sigma^{k}_{i=1}t_iv_{i,\epsilon},\quad w_{\epsilon,\bar{t}}:=\Sigma^{k}_{i=1}t_iw_{i,\epsilon},$$
where $u_{i,\epsilon},v_{i,\epsilon},w_{i,\epsilon}$ are given in (\ref{7.5}). One easily has that  $u_{\epsilon,\bar{t}}$, $v_{\epsilon,\bar{t}}$ and $w_{\epsilon,\bar{t}}$ are bounded in $H^1(\mathbb{R}^2)$. Let $\bar{h}_{\epsilon,\bar{t}}:=\bar{h}_\epsilon(w_{\epsilon,\bar{t}})$, as (\ref{3.1.3}) we get $\bar{h}_{\epsilon,\bar{t}}$ is bounded in $H^1(\mathbb{R}^2)$.

{\bf Step 1}. In this step, we estimate
\begin{equation}\label{5.43}
\int_{\mathbb{R}^2\backslash{\Lambda^\epsilon}}[|\nabla \psi_{\epsilon,\bar{t}}|^2+V(\epsilon x)\psi^2_{\epsilon,\bar{t}}]\rightarrow0, \ \text{as}\ \epsilon\rightarrow0, \ \text{uniformly in $\bar{t}$},\end{equation}
where
\begin{equation}\label{5.0.10}\psi_{\epsilon,\bar{t}}=\bar{h}_{\epsilon,\bar{t}}-\frac{u_{\epsilon,\bar{t}}-
v_{\epsilon,\bar{t}}}{2}=\bar{h}_{\epsilon,\bar{t}}
-\frac12\Sigma^k_{i=1}t_i(u_{i}-v_{i})\phi_i(\epsilon x).\end{equation}
Since (\ref{3.4}) holds with $\Phi_\epsilon$, $u$ and $h_\epsilon$ replacing by $\Psi_\epsilon$, $w_{\epsilon,\bar{t}}$ and $\bar{h}_{\epsilon,\bar{t}}$, we have
\begin{equation}\label{5.41} -2\Delta \bar{h}_{\epsilon,\bar{t}}+2V(\epsilon x)\bar{h}_{\epsilon,\bar{t}}=-\bar{f}_0(\epsilon x,w_{\epsilon,\bar{t}}+
\bar{h}_{\epsilon,\bar{t}})+\bar{g}_0(\epsilon x,w_{\epsilon,\bar{t}}-
\bar{h}_{\epsilon,\bar{t}}). \end{equation}
Note that $u_{\epsilon,\bar{t}}=v_{\epsilon,\bar{t}}=0$ in $\mathbb{R}^2\backslash{{\tilde{\Lambda}}^\epsilon}$, by (\ref{5.41}) we get
\begin{equation*} -2\Delta \bar{h}_{\epsilon,\bar{t}}+2V(\epsilon x)\bar{h}_{\epsilon,\bar{t}}=-\bar{f}_0(\epsilon x,
\bar{h}_{\epsilon,\bar{t}})+\bar{g}_0(\epsilon x,-
\bar{h}_{\epsilon,\bar{t}}), \ \text{in}\ \mathbb{R}^2\backslash{\tilde{\Lambda}^\epsilon}. \end{equation*}
For any $\xi\in H^1(\mathbb{R}^2)$ satisfying $\xi=1$ in $\mathbb{R}^2\backslash\tilde{\Lambda}^\epsilon$, letting $\xi_\epsilon(x)=\xi(\epsilon x)$ we get
\begin{equation*}\aligned
&2\int_{\mathbb{R}^2\backslash{\tilde{\Lambda}^\epsilon}}(|\nabla \bar{h}_{\epsilon,\bar{t}}|^2+V(\epsilon x)\bar{h}^2_{\epsilon,\bar{t}})\xi_\epsilon\\=&
-2\epsilon\int_{\mathbb{R}^2\backslash\tilde{\Lambda}^\epsilon}
\bar{h}_{\epsilon,\bar{t}}\nabla \bar{h}_{\epsilon,\bar{t}}\nabla \xi_\epsilon -
\int_{\mathbb{R}^2\backslash\tilde{\Lambda}^\epsilon}[\bar{f}_0(\epsilon x,\bar{h}_{\epsilon,\bar{t}})\bar{h}_{\epsilon,\bar{t}}\xi_\epsilon-
\bar{g}_0(\epsilon x,-\bar{h}_{\epsilon,\bar{t}})\bar{h}_{\epsilon,\bar{t}}\xi_\epsilon]
\leq o_\epsilon(1).
\endaligned
\end{equation*}
Using the fact that $u_{\epsilon,\bar{t}}=v_{\epsilon,\bar{t}}=0$ in $\mathbb{R}^2\backslash{\tilde{\Lambda}^\epsilon}$ again, we have
$\int_{\mathbb{R}^2\backslash{\tilde{\Lambda}^\epsilon}}[|\nabla \psi_{\epsilon,\bar{t}}|^2+V(\epsilon x)\psi^2_{\epsilon,\bar{t}}]\rightarrow0,$ as $\epsilon\rightarrow0$.

To show (\ref{5.43}), it suffices to show for any $i=1,\cdots,k$
\begin{equation}\label{5.0.9}
\int_{\tilde{\Lambda}^\epsilon_i\backslash{{\Lambda}^\epsilon_i}}
\bigl[|\nabla \psi_{\epsilon,\bar{t}}|^2+V(\epsilon x)\psi^2_{\epsilon,\bar{t}}\bigr]\rightarrow0,\ \text{as}\ \epsilon\rightarrow0\ \text{uniformly in} \ \bar{t}.
\end{equation}
In fact, denote
$\bar{h}^i_{\epsilon,\bar{t}}=\bar{h}_{\epsilon,\bar{t}}(x+\frac{x_i}{\epsilon})$. Over the set $\mathbb{R}^2\backslash{\cup_{j\neq i}(\tilde{\Lambda}_j)^\epsilon_{x_i}}$, where $(\tilde{\Lambda}_j)^\epsilon_{x_i}=(\tilde{\Lambda}_j-x_i)/\epsilon$ we have
 $$-2\Delta \bar{h}^i_{\epsilon,\bar{t}}+2V(\epsilon x+x_i)\bar{h}^i_{\epsilon,\bar{t}}=\bar{g}_0(\epsilon x+x_i,t_iw_i\phi^\epsilon_{x_i}-\bar{h}^i_{\epsilon,\bar{t}})-\bar{f}_0(\epsilon x+x_i,t_iw_i\phi^\epsilon_{x_i}+\bar{h}^i_{\epsilon,\bar{t}}),$$
 where $\phi^\epsilon_{x_i}$ is given in (\ref{5.6}).
 Since $(u_i,v_i)$ is a solution of system (\ref{bu2}), we have
\begin{equation}\label{5.0.12}-\Delta(u_i-v_i)+V(x_i)(u_i-v_i)=g(v_i)-f(u_i).\end{equation}
 Together with the fact that
 $\|(1-\phi^\epsilon_{x_i})u_i\|\rightarrow0$ and  $\|(1-\phi^\epsilon_{x_i})v_i\|\rightarrow0$, we infer
 $$-\Delta[(u_i-v_i)\phi^\epsilon_{x_i}]+V(\epsilon x+x_i)(u_i-v_i)\phi^\epsilon_{x_i}=g(v_i)-f(u_i)+o_\epsilon(1).$$
 Let $\psi^i_{\epsilon,\bar{t}}=\psi_{\epsilon,\bar{t}}(x+\frac{x_i}{\epsilon})$.
Then by (\ref{5.0.10}) we know
$\psi^i_{\epsilon,\bar{t}}=\bar{h}^i_{\epsilon,\bar{t}}
-\frac12\Sigma^k_{i=1}t_i(u_i-v_i)\phi^\epsilon_{x_i}.
$
Over the set $\mathbb{R}^2\backslash{\cup_{j\neq i}(\tilde{\Lambda}_j)^\epsilon_{x_i}}$, by (\ref{5.41}) there holds
 \begin{equation}\label{5.45}\aligned
 -2\Delta\psi^i_{\epsilon,\bar{t}}+2V(\epsilon x+x_i)\psi^i_{\epsilon,\bar{t}}=&t_if(u_i)-\bar{f}(\epsilon x+x_i,t_iw_i\phi^\epsilon_{x_i}+\bar{h}^i_{\epsilon,\bar{t}})\\&+\bar{g}(\epsilon x+x_i,t_iw_i\phi^\epsilon_{x_i}-\bar{h}^i_{\epsilon,\bar{t}})-t_ig(v_i)
 +o_\epsilon(1).
 \endaligned\end{equation}
Below as the proof of \cite[Page 161, Lines 8 to 18]{HUGO-PHD}, we get (\ref{5.0.9}) yields and then (\ref{5.43}) follows.

{\bf Step 2}: In this step, we shall study the asymptotic behavior of $\psi^i_{\epsilon,\bar{t}}$ in $\Lambda^\epsilon_i$. For $w_i=\frac{u_i+v_i}{2}$, denote $h_{t_i}=h_{V(x_i)}(t_iw_i)$ for simplicity and $\psi_{t_i}=h_{t_i}-t_i\frac{u_i-v_i}{2}$.
As (\ref{5.41}) we have
$$-2\Delta h_{t_i}+2V(x_i)h_{t_i}=-f(t_iw_i+h_{t_i})+g(t_iw_i-h_{t_i}).$$
Then by (\ref{5.0.12}) we get
\begin{equation*}\label{5.46}
-2\Delta \psi_{t_i}+2V(x_i)\psi_{t_i}=t_if(u_i)-f(t_iu_i+\psi_{t_i})
+g(t_iv_i-\psi_{t_i})-
t_ig(v_i) \text{\ in\ }\mathbb{R}^2.
\end{equation*}
Together with (\ref{5.45}) and (\ref{5.43}), and taking similar arguments of (5.48) in \cite{HUGO-PHD}, we obtain
\begin{equation}\label{5.48}
\int_{\tilde{\Lambda}^\epsilon_{x_i}}[|\nabla(\psi^i_{\epsilon
,\bar{t}}-\psi_{t_i})|^2+V(\epsilon x+x_i)(\psi^i_{\epsilon
,\bar{t}}-\psi_{t_i})^2]\rightarrow0, \ \text{as}\ \epsilon\rightarrow0, \ \text{uniformly in }\bar{t}.
\end{equation}

{\bf Step 3}. In this step, we show (\ref{5.5.0}) and (\ref{5.5.1}). Now we consider the function $\theta_{i,\epsilon}(\bar{t}):=\langle I'_\epsilon(w_{\epsilon,\bar{t}}),w_{\epsilon,\bar{t}}\phi^\epsilon_i\rangle$.
By (\ref{2.2.2}) with $\Phi_\epsilon$ and $h_\epsilon$ replacing by $\Psi_\epsilon$ and $\bar{h}_\epsilon$ respectively we deduce
$$\aligned
\theta_{i,\epsilon}(\bar{t})
=&2t^2_i\bigl(w_{i}\phi^\epsilon_{x_i},w_{i} (\phi^\epsilon_{x_i})^2\bigr)_{V(\epsilon x+x_i)}-\int_{\mathbb{R}^2}\bar{f}_0(\epsilon x+x_i,t_iu_{i}\phi^\epsilon_{x_i}+\psi^i_{\epsilon,\bar{t}})t_iw_{i}
(\phi^\epsilon_{x_i})^2\\&-
\int_{\mathbb{R}^2}\bar{g}_0(\epsilon x+x_i,t_iv_{i}\phi^\epsilon_{x_i}-\psi^i_{\epsilon,\bar{t}})
t_iw_{i}(\phi^\epsilon_{x_i})^2,
\endaligned$$
where $(u,v)_{V(\epsilon x+x_i)}=\int_{\mathbb{R}^2}[\nabla u\nabla v+V(\epsilon x+x_i)uv]$.
Define
$$\aligned\theta_i(t_i)=\langle J'_{V(x_i)}(t_iw_{i}),w_{i}\rangle=2t_i\|w_i\|^2_{V(x_i)}
-\int_{\mathbb{R}^2}f(t_iu_i+\psi_{t_i})w_i-
\int_{\mathbb{R}^2}g(t_iv_i-\psi_{t_i})w_i.\endaligned$$
Then by (\ref{5.48}) we get $\theta_{i,\epsilon}(\bar{t})=t_i\theta_i(t_i)+o_\epsilon(1)$.
Recall that $\theta_i\in C^1$, $\theta_i(1)=0$ and
$\theta'_i(1)=\langle J''_{V(x_i)}(w_i)w_i,w_i\rangle<0$. Then
$$\theta_i(t_i)=\theta_i(1)+\theta'_i(1)(t_i-1)+o(t_i-1)
=\theta'_i(1)(t_i-1)+o(t_i-1).$$
Choose $\mu$ such that $o(t_i-1)\leq |\theta'_i(1)||t_i-1|/2$ for $t_i\in[1-\mu,1+\mu]$. Applying Miranda's theorem, there exists $\epsilon_0>0$ such that for any $\epsilon\in(0,\epsilon_0)$, there is $\bar{t}_\epsilon=(t_{1,\epsilon}, \cdots, t_{k,\epsilon})\in[1-\mu.1+\mu]^k$ such that $\theta_{i,\epsilon}({t}_{1,\epsilon},\cdots,{t}_{k,\epsilon})=0$ for every $i=1,\cdots,k$. Moreover, by changing $\mu$, we can take ${t}_{i,\epsilon}\rightarrow1$. Then
 $h_{{t}_{i,\epsilon}}=h_{V(x_i)}({t}_{i,\epsilon}w_i)\rightarrow h_{V(x_i)}(w_i)=(u_i-v_i)/2$ and $\psi_{{t}_{i,\epsilon}}=h_{{t}_{i,\epsilon}}-t_i\frac{u_i-v_i}{2}
 \rightarrow 0$.
Therefore, by  (\ref{5.48}) we get
$$\aligned
&2\int_{\Lambda^\epsilon_i}({w}^2_{\epsilon,\bar{t}_\epsilon}+
\bar{h}^2_\epsilon(w_{\epsilon,\bar{t}_\epsilon}))
=\int_{\Lambda^\epsilon_i}[({u}_{\epsilon,\bar{t}_\epsilon}+
\psi_{\epsilon,\bar{t}_\epsilon})^2+({v}_{\epsilon,\bar{t}_\epsilon}-
\psi_{\epsilon,\bar{t}_\epsilon})^2]\\
&
\geq\int_{\Lambda^\epsilon_{x_i}}\bigl[
(t_{i,\epsilon}u_i
+\psi_{{t}_{i,\epsilon}})^2
+(t_{i,\epsilon}v_i-\psi_{{t}_{i,\epsilon}})^2\bigr]+o_\epsilon(1)\geq\eta,
\endaligned$$
for some $\eta>0$.

For the above $t_{i,\epsilon}$ and $\bar{t}_\epsilon$,
by (\ref{5.43}) and (\ref{5.48}) we obtain
\begin{equation*}\aligned
&I_\epsilon(w_{\epsilon,\bar{t}_\epsilon})
=\Psi_\epsilon(\Sigma^k_{i=1}t_{i,\epsilon}u_{i,\epsilon}+
\psi_{\epsilon,\bar{t}_\epsilon},\Sigma^k_{i=1}t_{i,\epsilon}
v_{i,\epsilon}-\psi_{\epsilon,\bar{t}_\epsilon})\\
=&\Sigma^k_{i=1}\Phi_{V_i}(t_{i,\epsilon} u_i+\psi_{t_{i,\epsilon}}, t_{i,\epsilon} v_i-\psi_{t_{i,\epsilon}})+o_\epsilon(1)=\Sigma^k_{i=1}\Phi_{V_i}( u_i, v_i)+o_\epsilon(1)=\Sigma^k_{i=1}c_i+o_\epsilon(1).\endaligned\end{equation*}
Hence,
letting $\bar{w}_\epsilon:=w_{\epsilon,\bar{t}_\epsilon}$, the conclusions follow directly.\ \ \ \ $\Box$

\subsection{Proof of Theorem \ref{t1.4}}
In this subsection, we are in a position to prove Theorem \ref{t1.4}.

\begin{lemma}\label{l5.22}For small $\epsilon>0$, there exists $w_\epsilon\in M_\epsilon$ such that $I_\epsilon(w_\epsilon)=d_\epsilon$. Moreover, $(u_\epsilon,v_\epsilon)
:=(w_\epsilon+\bar{h}_\epsilon(w_\epsilon),w_\epsilon-\bar{h}_\epsilon(w_\epsilon))$ is a solution of (\ref{7.0}) and $\Psi_\epsilon(u_\epsilon,v_\epsilon)=d_\epsilon$.
\end{lemma}
{\bf Proof}: By Lemma \ref{p5.15}, $d_\epsilon\geq\rho_0$. The Ekeland variational principle implies that there exists a constrained (PS)$_{d_\epsilon}$ sequence $w_n\in M_\epsilon$ of $I_\epsilon|_{M_\epsilon}$. Letting $u_n:=w_n+\bar{h}_\epsilon(w_n)$ and $v_n:=w_n-\bar{h}_\epsilon(w_n)$, as Lemma \ref{l2.5}-(2) we know,
$(u_n,v_n)\in \tilde{M}_\epsilon$ is a (PS)$_{d_\epsilon}$ sequence of $\Psi_\epsilon$, where $\tilde{M}_\epsilon$ is given in (\ref{5.0.11}).
Then there exist $\lambda^n_i\in\mathbb{R}$, $i=1,\cdots,k$, and $\psi_n\in H^1(\mathbb{R}^2)$ such that, for any $\zeta,\xi\in H^1(\mathbb{R}^2)$,
$\Psi_\epsilon(u_n,v_n)\rightarrow c_\epsilon,$ and
\begin{equation*}\label{5.52.0}\aligned&
\langle{\Psi}'_\epsilon(u_n,v_n),
(\zeta,\xi)\rangle=\langle{\Psi}'_\epsilon(u_n,v_n),
(\Sigma^k_{i=1}\lambda^n_i\phi^\epsilon_i\zeta,
\Sigma^k_{i=1}\lambda^n_i\phi^\epsilon_i\xi)\rangle
\\&+\langle{\Psi}''_\epsilon(u_n,v_n)
(\Sigma^k_{i=1}\lambda^n_i\phi^\epsilon_i u_n+\psi_n,\Sigma^k_{i=1}\lambda^n_i\phi^\epsilon_i v_n-\psi_n),(\zeta,\xi)\rangle+o_n(1)\|(\zeta,\xi)\|_{E}.\endaligned
\end{equation*}
From the above equalities with $\zeta:=\zeta_n=\Sigma^k_{i=1}\lambda^n_i\phi^\epsilon_iu_n+\psi_n$, $\xi:=\xi_n=\Sigma^k_{i=1}\lambda^n_i\phi^\epsilon_iv_n-\psi_n$, and Lemma \ref{l5.3} it follows that
\begin{equation}\label{5.52.1}\aligned0=
\langle{\Psi}'_\epsilon(u_n,v_n),(\zeta_n,\xi_n) \rangle
\leq\Sigma^k_{i=1}(\lambda^n_i)^2(-\eta)-\frac k2\|\psi_n\|^2_\epsilon+o_n(1)\|(\zeta_n,\xi_n)\|_{1,\epsilon}.
\endaligned\end{equation}
From Lemma \ref{p5.14} we know $d_\epsilon\leq 2\Sigma^k_{i=1} c_i$. Then  Lemma \ref{p5.15} implies that $\|(u_n,v_n)\|\leq C$, where $C$ is independent on $n$. Then
$$\aligned o_n(1)\|(\zeta_n,\xi_n)\|_{1,\epsilon}&=o_n(1)
\Sigma^k_{i=1}(\lambda^n_i)^2+o_n(1)\|\psi_n\|^2_\epsilon,
\endaligned$$
together with (\ref{5.52.1}) we deduce
 $\lambda^n_i\rightarrow0$ for any $i=1,\cdots,k$ and $\psi_n\rightarrow0$ in $H^1(\mathbb{R}^2)$.
So $(u_n,v_n)$ is a PS sequence of $\Psi_\epsilon$. Similar to the argument of Lemma \ref{l2.2}, we get $(u_n,v_n)\rightarrow(u_\epsilon,v_\epsilon)$ in $E$. Then $\Psi'_\epsilon(u_\epsilon,v_\epsilon)=0$. So $I'_\epsilon(w_\epsilon)=0$ with $w_\epsilon=\frac{u_\epsilon+v_\epsilon}{2}$, $w_\epsilon\in M_\epsilon$, and $I_\epsilon(w_\epsilon)=\Psi_\epsilon(u_\epsilon,v_\epsilon)=d_\epsilon$.\ \ \ \ \ $\Box$

From now on we consider the solution $(u_\epsilon,v_\epsilon)$ of system (\ref{7.0}) obtained by Lemma \ref{l5.22}.
The rest of the proof of Theorem \ref{t1.4} will
be carried out in a series of lemmas.

\begin{lemma}\label{l5.23}
$d_\epsilon=\Sigma^k_{i=1}c_i+o_\epsilon(1)$ and $\Psi^i_\epsilon(u_\epsilon,v_\epsilon)=c_i+o_\epsilon(1)$, for any $i=1,\cdots,k$,
where
\begin{equation}\label{5.53.0}
\Psi^i_\epsilon(u,v):=\int_{\tilde{\Lambda}^\epsilon_i}\bigl[\nabla u\nabla v+V(\epsilon x)uv-\bar{F}_0(\epsilon x,u)-\bar{G}_0(\epsilon x,v)\bigr].
\end{equation}
Moreover, there exist $x^\epsilon_i\in \Lambda^\epsilon_i$ and $(\bar{u}_i,\bar{v}_i)$ ground states of system (\ref{bu2}) such that $\|u_\epsilon-\Sigma^k_{i=1}\bar{u}_i(\cdot-x^\epsilon_i)\|\rightarrow0$, and
$\|v_\epsilon-\Sigma^k_{i=1}\bar{v}_i(\cdot-x^\epsilon_i)\|\rightarrow0$.\end{lemma}
{\bf Proof}: We divide the proof into several steps.

{\bf Step 1}. Let $\xi$ be a cut-off function such that $\xi=0$ in $\Lambda$ and $\xi=1$ in $\mathbb{R}^N\backslash{\Lambda'}$ and let $\xi_\epsilon(x)=\xi(\epsilon x)$. Note that
 Testing $\Psi'_\epsilon(u_\epsilon,v_\epsilon)=0$ with $(u_\epsilon\xi_\epsilon,v_\epsilon\xi_\epsilon)$, we get
 $$\aligned
&\int_{\mathbb{R}^2}
(|\nabla u_\epsilon|^2+V(\epsilon x)u^2_\epsilon)\xi_\epsilon+\int_{\mathbb{R}^2}
(|\nabla v_\epsilon|^2+V(\epsilon x)v^2_\epsilon)\xi_\epsilon\\
=&\int_{\mathbb{R}^2\backslash\Lambda^\epsilon}[\bar{f}_0(\epsilon x,u_\epsilon)v_\epsilon+\bar{g}_0(\epsilon x,v_\epsilon)u_\epsilon]\xi_\epsilon-
\epsilon\int_{\mathbb{R}^2}\nabla u_\epsilon u_\epsilon\nabla\xi_\epsilon-
\epsilon\int_{\mathbb{R}^2}\nabla v_\epsilon v_\epsilon\nabla\xi_\epsilon\\
\leq&\delta\int_{\mathbb{R}^2}(u^2_\epsilon+v^2_\epsilon)
\xi_\epsilon+\epsilon\|(u_\epsilon,v_\epsilon)\|_{1,\epsilon}. \endaligned$$
Then
\begin{equation}\label{5.53.2}\int_{\mathbb{R}^2
\backslash{{(\Lambda')}^\epsilon}}[|\nabla u_\epsilon|^2+V(\epsilon x)u^2_\epsilon+|\nabla v_\epsilon|^2+V(\epsilon x)v^2_\epsilon]=o_\epsilon(1).\end{equation}
Hence
\begin{equation}\label{claim1}\Psi^i_\epsilon(u_\epsilon,v_\epsilon)=
\Psi_\epsilon(u_\epsilon\phi^\epsilon_i,v_\epsilon\phi^\epsilon_i)
+o_\epsilon(1),\ \Psi_\epsilon(u_\epsilon,v_\epsilon)=
\Sigma^k_{i=1}\Psi^i_\epsilon(u_\epsilon,v_\epsilon)+o_\epsilon(1).
\end{equation}
Using Lemma \ref{p5.14} we get $d_\epsilon\leq \Sigma^k_{i=1}c_i+o_\epsilon(1)$. So it suffices to show \begin{equation}\label{5.55}\Psi^i_\epsilon(u_\epsilon,v_\epsilon)\geq c_i+o_\epsilon(1).\end{equation}

{\bf Step 2}. We claim that there exist
$x^\epsilon_i\in \Lambda^\epsilon_i$ and $R_0,r>0$ such that
\begin{equation}\label{5.57}
\int_{B_{R_0}(x^\epsilon_i)}u^2_\epsilon\geq r.\end{equation}
Otherwise we assume, for any $r>0$,
$\sup_{y\in \Lambda^{\epsilon}_i}\int_{B_r(y)}u^2_\epsilon(x)\rightarrow0$, as $\epsilon\rightarrow0.$
Consider a cut-off function $\phi^\epsilon_R\in[0,1]$ satisfying $|\nabla \phi^\epsilon_R|\leq\frac{C}{R}$, $\phi^\epsilon_R(x)=1$ in $N_R(\Lambda^\epsilon_i):=\{x\in\mathbb{R}^2: dist(x,\Lambda^\epsilon_i)\leq R\}$ and $\phi^\epsilon_R(x)=0$ in $\mathbb{R}^2\backslash N_{2R}(\Lambda^\epsilon_i)$.
Taking similar arguments of (5.56) in \cite{HUGO-PHD}, for each fixed $R>0$, we get $\int_{N_R(\Lambda^\epsilon_i)}|u_\epsilon|^q\rightarrow0$ with any $q>2$.
Testing $\Psi'_\epsilon(u_\epsilon,v_\epsilon)=0$ with $(u_\epsilon\phi^\epsilon_R, v_\epsilon\phi^\epsilon_R)$ and using (\ref{2.1.0}) we infer
$$\aligned
&\int_{\mathbb{R}^2}(|\nabla u_\epsilon|^2+V(\epsilon x)u^2_\epsilon)\phi^\epsilon_R+(|\nabla v_\epsilon|^2+V(\epsilon x)v^2_\epsilon)\phi^\epsilon_R\\
=&\int_{\mathbb{R}^2}\bar{f}_0(\epsilon x,u_\epsilon)v_\epsilon\phi^\epsilon_R+
\int_{\mathbb{R}^2}\bar{g}_0(\epsilon x,v_\epsilon)u_\epsilon\phi^\epsilon_R-
\int_{\mathbb{R}^2}(u_\epsilon\nabla u_\epsilon\nabla{\phi^\epsilon_R}+
v_\epsilon\nabla v_\epsilon\nabla{\phi^\epsilon_R})\\
\leq&\delta\int_{\mathbb{R}^2}(u^2_\epsilon+v^2_\epsilon)\phi^\epsilon_R
+C_\delta\int_{N_{2R}(\Lambda^\epsilon_i)}\bigl[|u_\epsilon|(e^{\alpha u^2_\epsilon}-1)|v_\epsilon|+|v_\epsilon|(e^{\alpha v^2_\epsilon}-1)|u_\epsilon|\bigr]+\frac{C}{R}.
\endaligned$$
Note that $d_\epsilon\leq \Sigma^k_{i=1} c_i+o_\epsilon(1)$. Similar to the argument of Lemma \ref{l2.2} we get $(u_\epsilon,v_\epsilon)$ is bounded in $E$. Similar to (\ref{4.22.0}), taking suitable $\alpha$ in (\ref{2.1.0}) we get $\int_{\mathbb{R}^2}(e^{3\alpha u^2_\epsilon}-1)\leq C$ and $\int_{\mathbb{R}^2}(e^{3\alpha v^2_\epsilon}-1)\leq C$. Then choosing $\delta\leq\frac{\inf_{\mathbb{R}^2}V}{2}$, we get
$$\aligned
&\int_{\mathbb{R}^2}(|\nabla u_\epsilon|^2+V(\epsilon x)u^2_\epsilon)\phi^\epsilon_R+(|\nabla v_\epsilon|^2+V(\epsilon x)v^2_\epsilon)\phi^\epsilon_R\\
\leq&C \bigl(\int_{N_{2R}(\Lambda^\epsilon_i)}|u_\epsilon|^3\bigr)^{\frac13}
\Bigl[\bigl(\int_{\mathbb{R}^2}(e^{3\alpha u^2_\epsilon}-1)\bigr)^{\frac13}
+\bigl(\int_{\mathbb{R}^2}(e^{3\alpha v^2_\epsilon}-1)\bigr)^{\frac13}\Bigr]|v_\epsilon|_3+\frac{C}{R}=o_\epsilon(1)+\frac{C}{R}.
\endaligned$$
Taking $R$ large and afterwards $\epsilon$ small,
we obtain $\int_{\Lambda^\epsilon_i}(u^2_\epsilon +v^2_\epsilon)\rightarrow0$. However, by Lemma \ref{p5.15} we know
$$\int_{\Lambda^\epsilon_i}(w^2_\epsilon+\bar{h}^2_\epsilon
(w_\epsilon))\geq\eta_0,$$
where $w_\epsilon=\frac{u_\epsilon+v_\epsilon}{2}$ and $\bar{h}_\epsilon(w_\epsilon)=\frac{u_\epsilon-v_\epsilon}{2}$. This is a contradiction. Hence (\ref{5.57}) holds.
Up to a subsequence, we assume $\epsilon x^\epsilon_i\rightarrow y_i\in \bar{\Lambda}_i$.

{\bf Step 3}.
Let $\bar{u}_\epsilon(x)=u_\epsilon(x+x^\epsilon_i)$, $\bar{v}_\epsilon(x)=v_\epsilon(x+x^\epsilon_i)$. Then $(\bar{u}_\epsilon,\bar{v}_\epsilon)$ satisfies
 \begin{equation}\label{5.1}\aligned
\left\{ \begin{array}{lll}
-\Delta\bar{u}_\epsilon+V(\epsilon x+\epsilon x^\epsilon_i)\bar{u}_\epsilon=\bar{g}_0(\epsilon x+\epsilon x^\epsilon_i,\bar{v}_\epsilon)\ & \text{in}\quad \mathbb{R}^2,\\
-\Delta\bar{v}_\epsilon+V(\epsilon x+\epsilon x^\epsilon_i)\bar{v}_\epsilon=\bar{f}_0(\epsilon x+\epsilon x^\epsilon_i,\bar{u}_\epsilon)\ & \text{in}\quad \mathbb{R}^2.
\end{array}\right.\endaligned
\end{equation}
Denote $\bar{\phi^\epsilon_i}(x):=\phi^\epsilon_i(x+x^\epsilon_i)
=\phi_i(\epsilon x+\epsilon x^\epsilon_i)$. Recall that $({u}_\epsilon,{v}_\epsilon)$ is bounded in $E$, together with (\ref{5.57}) we may assume $(\bar{u}_\epsilon,\bar{v}_\epsilon)
\rightharpoonup(\bar{u}_i,\bar{v}_i)\neq(0,0)$ in $E$. Then $(\bar{u}_\epsilon\bar{\phi^\epsilon_i},
\bar{v}_\epsilon\bar{\phi^\epsilon_i})\rightharpoonup(\bar{u}_i,\bar{v}_i)$ in $E$.
Since $\chi_{\Lambda}(\epsilon x+\epsilon x^\epsilon_i)$ is bounded, we assume $\chi_{\Lambda}(\epsilon x+\epsilon x^\epsilon_i)\rightharpoonup \chi$ with $\chi\in [0,1]$ in $L^p_{loc}(\mathbb{R}^2)$ for some $p>2$.
Then $(\bar{u}_i,\bar{v}_i)\neq(0,0)$ satisfies
\begin{equation}\label{5.58.1}
-\Delta \bar{u}_i+V(y_i)\bar{u}_i={g}_\infty(x,\bar{v}_i),\
-\Delta \bar{v}_i+V(y_i)\bar{v}_i={f}_\infty(x,\bar{u}_i)\quad  \text{in}\quad \mathbb{R}^2,
\end{equation}
where ${g}_\infty(x,s)=\chi(x)g(s)+(1-\chi(x))\tilde{g}(s)$, ${f}_\infty(x,s)=\chi(x)f(s)+(1-\chi(x))\tilde{f}(s)$.
Denote  ${G}_\infty(x,s):=\int^s_0 {g}_\infty(x,t)dt$, ${F}_\infty(x,s):=\int^s_0{f}_\infty(x,t)dt$ and the functional
$$\bar{\Psi}_{V(y_i)}({u},{v}):=\int_{\mathbb{R}^2}(\nabla u\nabla v+V(y_i)uv)-\int_{\mathbb{R}^2}({F}_\infty(x,u)+{G}_\infty(x,v)).$$
Denote the functional of (\ref{5.1}) by $$\bar{\Psi}_{\epsilon}(u,v)=\int_{\mathbb{R}^2}[\nabla u\nabla v+V(\epsilon x+\epsilon x^\epsilon_i)uv]-\int_{\mathbb{R}^2}[\bar{F}_0(\epsilon x+\epsilon x^\epsilon_i,u)+\bar{G}_0(\epsilon x+\epsilon x^\epsilon_i,v)].$$
Since $(\bar{u}_\epsilon\bar{\phi^\epsilon_i},
\bar{v}_\epsilon\bar{\phi^\epsilon_i})\rightharpoonup(\bar{u}_i,\bar{v}_i)$ in $E$, from Fatou Lemma it follows that
$$\aligned&\bar{\Psi}_{V(y_i)}(\bar{u}_i,\bar{v}_i)-\frac12\langle
\bar{\Psi}'_{V(y_i)}(\bar{u}_i,\bar{v}_i),(\bar{v}_i,\bar{u}_i)\rangle\\
\leq&\bar{\Psi}_{\epsilon}(\bar{u}_\epsilon\bar{\phi^\epsilon_i},\bar{v}_\epsilon\bar{\phi^\epsilon_i})-\frac12
\langle\bar{\Psi}'_{\epsilon}(\bar{u}_\epsilon\bar{\phi^\epsilon_i},
\bar{v}_\epsilon\bar{\phi^\epsilon_i}), (\bar{v}_\epsilon\bar{\phi^\epsilon_i},
\bar{u}_\epsilon\bar{\phi^\epsilon_i})\rangle.
\endaligned$$
Moreover, by (\ref{5.53.2}) we get $$\aligned
&\langle\bar{\Psi}'_{\epsilon}(\bar{u}_\epsilon\bar{\phi^\epsilon_i},
\bar{v}_\epsilon
\bar{\phi^\epsilon_i}), (\bar{v}_\epsilon
\bar{\phi^\epsilon_i},\bar{u}_\epsilon\bar{\phi^\epsilon_i})\rangle
=\langle\bar{\Psi}'_{\epsilon}(\bar{u}_\epsilon,\bar{v}_\epsilon), (\bar{v}_\epsilon\bar{\phi^\epsilon_i},\bar{u}_\epsilon\bar{\phi^\epsilon_i})
\rangle+o_\epsilon(1)=o_\epsilon(1).
\endaligned$$
Hence
\begin{equation}\label{5.51.0}\aligned
&\bar{\Psi}_{V(y_i)}(\bar{u}_i,\bar{v}_i)=
\bar{\Psi}_{V(y_i)}(\bar{u}_i,\bar{v}_i)-\frac12\langle
\bar{\Psi}'_{V(y_i)}(\bar{u}_i,\bar{v}_i),(\bar{v}_i,\bar{u}_i)\rangle
\\ \leq&\bar{\Psi}_{\epsilon}(\bar{u}_\epsilon\bar{\phi^\epsilon_i},\bar{v}_\epsilon\bar{\phi^\epsilon_i})
+o_\epsilon(1)={\Psi}_{\epsilon}({u}_\epsilon{\phi^\epsilon_i},{v}_\epsilon{\phi^\epsilon_i})
+o_\epsilon(1)={\Psi}^i_{\epsilon}
({u}_\epsilon,{v}_\epsilon)+o_\epsilon(1),
\endaligned\end{equation}
where we have used (\ref{claim1}).

{\bf Step 4}. In this step, we show $y_i\in \Lambda_i$ and  $(\bar{u}_i,\bar{v}_i)$ solves system (\ref{bu2}). Taking similar arguments as (\ref{4.17}) and combining with (\ref{5.51.0}) we get \begin{equation*}\label{5.55.0}c_i\leq c_{V(y_i)}\leq\bar{\Psi}_{V(y_i)}(\bar{u}_i,\bar{v}_i)\leq
\liminf_{\epsilon\rightarrow0}\Psi^i_\epsilon(u_\epsilon,v_\epsilon).\end{equation*}
Then (\ref{5.55}) yields and $\Psi^i_\epsilon(u_\epsilon,v_\epsilon)=c_i+o_\epsilon(1)$. Hence, $V(y_i)=\inf_{\Lambda_i}V=V_i$, $y_i\in \Lambda_i$ and $(\bar{u}_i,\bar{v}_i)$ solves system (\ref{bu2}).

{\bf Step 5}. In this step, we show $\|u_\epsilon(x)-\Sigma^k_{i=1}\bar{u}_i(x-x^\epsilon_i)\|\rightarrow0$ and $\|v_\epsilon(x)-\Sigma^k_{i=1}\bar{v}_i(x-x^\epsilon_i)\|\rightarrow0$.
In fact,
let $$w^1_\epsilon(\cdot):=(u_\epsilon\phi^\epsilon_i)(\cdot)-\bar{u}_i
(\cdot-x^\epsilon_i)\phi^\epsilon_i(\cdot),\quad
w^2_\epsilon(\cdot):=(v_\epsilon\phi^\epsilon_i)(\cdot)
-\bar{v}_i(\cdot-x^\epsilon_i)\phi^\epsilon_i(\cdot).$$
It is easy to see that $w^1_\epsilon, w^2_\epsilon\rightharpoonup0$ in $H^1(\mathbb{R}^2)$ and then
\begin{equation}\label{5.0.7}
\Psi_\epsilon(u_\epsilon\phi^\epsilon_i,v_\epsilon\phi^\epsilon_i)
-\Psi_\epsilon(\bar{u}_i(\cdot-x^\epsilon_i)\phi^\epsilon_i,
\bar{v}_i(\cdot-x^\epsilon_i)\phi^\epsilon_i)-\Psi_\epsilon(w^1_\epsilon,w^2_\epsilon)\rightarrow0.
\end{equation}
\begin{equation}\label{5.0.8}
\Psi'_\epsilon(u_\epsilon\phi^\epsilon_i,v_\epsilon\phi^\epsilon_i)
-\Psi'_\epsilon(\bar{u}_i(\cdot-x^\epsilon_i)\phi^\epsilon_i,
\bar{v}_i(\cdot-x^\epsilon_i)\phi^\epsilon_i)-\Psi'_\epsilon(w^1_\epsilon,w^2_\epsilon)\rightarrow0\ \text{in}\ H^{-1}(\mathbb{R}^2).
\end{equation}
Moreover, $\bar{u}_i(\cdot-x^\epsilon_i)(\phi^\epsilon_i-1)\rightarrow0$ and
$\bar{v}_i(\cdot-x^\epsilon_i)(\phi^\epsilon_i-1)\rightarrow0$ in $H^1(\mathbb{R}^2)$, and
$$\aligned&\int_{\mathbb{R}^2}V(\epsilon x+\epsilon x^\epsilon_i)\bar{u}_i\bar{v}_i\rightarrow\int_{\mathbb{R}^2}V_i\bar{u}_i\bar{v}_i,\ \int_{\mathbb{R}^2}\bar{F}_0(\epsilon x+\epsilon x^\epsilon_i,\bar{u}_i)\rightarrow\int_{\mathbb{R}^2}{F}(\bar{u}_i),\\& \int_{\mathbb{R}^2}\bar{G}_0(\epsilon x+\epsilon x^\epsilon_i,\bar{v}_i)\rightarrow\int_{\mathbb{R}^2}{G}(\bar{v}_i).\endaligned$$ Combining with (\ref{5.0.7}) and (\ref{5.51.0}) (which is an equality and $\bar{\Psi}_{V(y_i)}=\Phi_{V_i}$ by Step 4) we get
\begin{equation}\label{5.54.0}\aligned
&\Psi_\epsilon(w^1_\epsilon,w^2_\epsilon)=
\Psi_\epsilon(u_\epsilon\phi^\epsilon_i,v_\epsilon\phi^\epsilon_i)
-\Psi_\epsilon(\bar{u}_i(\cdot-x^\epsilon_i)\phi^\epsilon_i,
\bar{v}_i(\cdot-x^\epsilon_i)\phi^\epsilon_i)+o_\epsilon(1)\\
=&{\Psi}^i_\epsilon(u_\epsilon, v_\epsilon)-\Phi_{V_i}(\bar{u}_i,\bar{v}_i)+o_\epsilon(1)
=o_\epsilon(1).\endaligned\end{equation}
By (\ref{5.53.2}) we obtain
$$\langle\Psi'_\epsilon
(u_\epsilon\phi^\epsilon_i,v_\epsilon\phi^\epsilon_i),
(w^2_\epsilon,w^1_\epsilon)\rangle
=\langle\Psi'_\epsilon
(u_\epsilon,v_\epsilon),
(w^2_\epsilon,w^1_\epsilon)\rangle+o_\epsilon(1)
=o_\epsilon(1).
$$
Then in view of (\ref{5.0.8}) we have
\begin{equation}\label{5.10}\aligned \langle\Psi'_\epsilon(w^1_\epsilon,w^2_\epsilon),
(w^2_\epsilon,w^1_\epsilon)\rangle
&=o_\epsilon(1)-\langle\Psi'_\epsilon(\bar{u}_i(\cdot-x^\epsilon_i)\phi^\epsilon_i,
\bar{v}_i(\cdot-x^\epsilon_i)\phi^\epsilon_i), (w^2_\epsilon,w^1_\epsilon)\rangle\\
&=\langle\Phi'_{V_i}(\bar{u}_i,\bar{v}_i),(w^2_\epsilon(\cdot+x^\epsilon_i),
w^1_\epsilon(\cdot+x^\epsilon_i))\rangle
+o_\epsilon(1)=o_\epsilon(1).
\endaligned\end{equation}
From $\langle \Psi'_\epsilon(w^1_\epsilon,w^2_\epsilon),(w^2_\epsilon,w^1_\epsilon)
\rangle\rightarrow0$ and (\ref{5.19}) we obtain
\begin{equation}\label{5.0.5}\aligned
\|w^1_\epsilon\|^2_\epsilon+\|w^2_\epsilon\|^2_\epsilon=&
\int_{\mathbb{R}^2}\bigl[\bar{f}_0(\epsilon x,w^1_\epsilon)w^2_\epsilon+\bar{g}_0(\epsilon x,w^2_\epsilon)w^1_\epsilon
\bigr]
\\ \leq&\delta\int_{\mathbb{R}^2} \bigl[(w^1_\epsilon)^2+(w^2_\epsilon)^2\bigr]+C_\delta\int_{\Lambda^\epsilon_i} \bigl[f(w^1_\epsilon)w^1_\epsilon+g(w^2_\epsilon)w^2_\epsilon\bigr].
\endaligned\end{equation}
If $\|(w^1_\epsilon,w^2_\epsilon)\|_{1,\epsilon}\geq\eta$ for some $\eta>0$ independent on $\epsilon$, letting $\delta<\frac{\inf_{\mathbb{R}^2}V}{2}$ in (\ref{5.0.5}) we infer
$$\int_{\Lambda^\epsilon_i}\bigl(f(w^1_\epsilon)w^1_\epsilon+g(w^2_\epsilon)w^2_\epsilon\bigr)
\geq\eta.$$
On the other hand, from Fatou lemma and Lemma \ref{l2.0} it follows that, for small enough $\epsilon>0$
$$\aligned
2\Psi_\epsilon(w^1_\epsilon,w^2_\epsilon)-\langle
\Psi'_\epsilon(w^1_\epsilon,w^2_\epsilon),(w^2_\epsilon,w^1_\epsilon)\rangle
\geq\delta'\int_{\Lambda^\epsilon_i}[f(w^1_\epsilon)w^1_\epsilon
+g(w^2_\epsilon)w^2_\epsilon]+o_\epsilon(1)\geq
\delta'\eta/2,
\endaligned$$ contradicts (\ref{5.10}) and (\ref{5.54.0}).
 Therefore, $\|(w^1_\epsilon,w^2_\epsilon)\|_{1,\epsilon}\rightarrow0$, and then
 $u_\epsilon\phi^\epsilon_i-\bar{u}_i(\cdot-x^\epsilon_i)\rightarrow0$ and $v_\epsilon\phi^\epsilon_i-\bar{v}_i(\cdot-x^\epsilon_i)\rightarrow0$ in $H^1(\mathbb{R}^2)$.
Together with $\|u_\epsilon-\Sigma^k_{i=1}u_\epsilon\phi^\epsilon_i\|\rightarrow0$ using (\ref{5.53.2}), we have
$u_\epsilon-\Sigma^k_{i=1}\bar{u}_i(\cdot-x^\epsilon_i)\rightarrow0$ in  $H^1(\mathbb{R}^2).$ Similarly $v_\epsilon-\Sigma^k_{i=1}\bar{v}_i(\cdot-x^\epsilon_i)\rightarrow0$ in $H^1(\mathbb{R}^2)$.
\ \ \ \ \ $\Box$
\begin{lemma}\label{l5.27}
\noindent(1) If $P_\epsilon\in \bar{\Lambda}^\epsilon_i$ satisfies
$\liminf_{\epsilon\rightarrow0}\max\{|u_\epsilon(P_\epsilon)|, |v_\epsilon(P_\epsilon)|\}>0$, then $\lim_{\epsilon\rightarrow0}V(\epsilon P_\epsilon)=\inf_{\Lambda_i}V$.

\noindent(2) $({u}_\epsilon,{v}_\epsilon)$ satisfies $\lim_{\epsilon\rightarrow0}\sup_{\mathbb{R}^2
\backslash{\Lambda^\epsilon}}\{|u_\epsilon|, |v_\epsilon|\}=0$. In particular, $({u}_\epsilon,{v}_\epsilon)$ solves system (\ref{2.1}).

\noindent(3) $\liminf_{\epsilon\rightarrow0}\min\{\sup_{\Lambda^\epsilon_i}|u_\epsilon|,
\sup_{\Lambda^\epsilon_i}|v_\epsilon|\}>0$ for all $i=1,\cdots,k$.\end{lemma}
{\bf Proof}: (1) Using similar arguments of Lemma \ref{l5.23}, we get desired results.

(2) Firstly, we claim that
\begin{equation}\label{5.0.3}\sup_{\partial {\Lambda^\epsilon_i}}|u_\epsilon|\rightarrow0,\ \sup_{\partial {\Lambda^\epsilon_i}}|v_\epsilon|\rightarrow0, \quad\forall i=1,\cdots,k.\end{equation}
Otherwise, $\liminf_{\epsilon\rightarrow0}\sup_{\partial {\Lambda^\epsilon_i}}\{|u_\epsilon|, |v_\epsilon|\}>0$. Then there exists $z_\epsilon\in \partial {\Lambda^\epsilon_i}$ such that $\epsilon z_\epsilon\rightarrow\bar{z}\in \partial \Lambda_i$ and $\liminf_{\epsilon\rightarrow0}\max\{|u_\epsilon(z_\epsilon)|, |v_\epsilon(z_\epsilon)|\}>0$. By (1), we know $\lim_{\epsilon\rightarrow0}V(\epsilon z_\epsilon)=V(\bar{z})=\inf_{\Lambda_i}V$. This is impossible since $\bar{z}\in \partial {\Lambda_i}$.

In the same way as (\ref{5.0.2}), there holds
$-\Delta (|u_\epsilon|+|v_\epsilon|)\leq0$ in $\mathbb{R}^2\backslash{\Lambda^\epsilon}$. Applying the maximum principle and using (\ref{5.0.3}) we get $\lim_{\epsilon\rightarrow0}
\sup_{\mathbb{R}^2\backslash{\Lambda^\epsilon}}
\{|u_\epsilon|,|v_\epsilon|\}=0.$
Then $\bar{f}(\epsilon x,u_\epsilon)=f(u_\epsilon)$ and $\bar{g}(\epsilon x,v_\epsilon)=g(v_\epsilon)$. Consequently, $(u_\epsilon,v_\epsilon)$ solves (\ref{2.1}).

(3) Assume that $\sup_{\Lambda^\epsilon_i}|u_\epsilon|\rightarrow0$. Combining the conclusion (2) we get $\int_{B_{r}(x^\epsilon_i)}u^2_\epsilon\rightarrow0$ for any fixed $r>0$, contradicts with (\ref{5.57}). Therefore,
$\sup_{\Lambda^\epsilon_i}|u_\epsilon|\geq\rho>0$ for small $\epsilon$. Reasoning exactly as (\ref{5.57}), there exists $y^\epsilon_i\in \Lambda^\epsilon_i$ and $R'_0, \rho'_0>0$ such that
$\int_{B_{R'_0}(y^\epsilon_i)}v^2_\epsilon\geq\rho'_0>0$. Then as before $\sup_{\Lambda^\epsilon_i}|v_\epsilon|\geq\rho>0$ for small $\epsilon$.\ \ \ \ $\Box$

\begin{lemma}\label{l5.29}\noindent(1) For small enough $\epsilon>0$, $|u_\epsilon|$, $|v_\epsilon|$ and $|u_\epsilon|+|v_\epsilon|$ have one local maximum point $x^\epsilon_{i}$, $y^\epsilon_{i}$ and $z^\epsilon_{i}$ in $\Lambda^\epsilon_i$ respectively for each $i=1,\cdots,k$.\\
\noindent(2) $|u_\epsilon|+|v_\epsilon|$ does not admit local maximum points in $\mathbb{R}^2\backslash{\Lambda^\epsilon}$ for sufficiently small $\epsilon$.\end{lemma}
{\bf Proof}: (1) Using Lemma \ref{l5.27}, the conclusion follows directly.

(2) Treating as (\ref{5.0.1}) we get
\begin{equation}\label{bu7}-\Delta(|u_\epsilon|+|v_\epsilon|)+V(\epsilon x)(|u_\epsilon|+|v_\epsilon|)\leq |f(u_\epsilon)|+|g(v_\epsilon)|\ \ \text{in}\ \mathbb{R}^2.\end{equation}
If $z_\epsilon\in \mathbb{R}^2\backslash{\Lambda^\epsilon}$ is a local maximum point of $|u_\epsilon|+|v_\epsilon|$, then by (\ref{bu7}) and (H$_1$) we infer that $|u_\epsilon(z_\epsilon)|+|v_\epsilon(z_\epsilon)|\geq\rho$ for some $\rho>0$. This contradicts Lemma \ref{l5.27}-(2).\ \ \ $\Box$

\begin{lemma}\label{l5.14} There exist $C_0, c_0>0$ such that
$$|u_\epsilon(x)|+|v_\epsilon(x)|\leq C_0e^{-c_0|x-z^\epsilon_i|},\quad\text{for every}\ x\in\mathbb{R}^2\backslash{\cup_{j\neq i}\Lambda^\epsilon_j}.$$
\end{lemma}
{\bf Proof}:  Take $\Lambda''_i\Subset \Lambda_i$ such that
$\inf_{\Lambda_i}V=\inf_{\Lambda''_i}V<\inf_{\partial\Lambda''_i}V$.  Reasoning as
Lemma \ref{l5.27} we have $\sup_{\mathbb{R}^2\backslash{(\Lambda'')^\epsilon}}
(|u_\epsilon|+|v_\epsilon|)\rightarrow0$ with $\Lambda'':=\cup^k_{i=1}\Lambda''_i$,
 and $|u_\epsilon|+|v_\epsilon|$ does not admit any local maximum points
over $\mathbb{R}^2\backslash{(\Lambda'')^\epsilon}$  for sufficiently small $\epsilon$.

For fixed $i\in\{1,\cdots,k\}$, as in Lemma \ref{l5.23}, define $\bar{u}_\epsilon(x)={u}_\epsilon(x+z^\epsilon_{i})$ and $\bar{v}_\epsilon(x)={v}_\epsilon(x+z^\epsilon_{i})$ we have
$(\bar{u}_\epsilon,\bar{v}_\epsilon)\rightharpoonup (\bar{u}_i,\bar{v}_i)\neq(0,0)$
in $E$, and $(\bar{u}_\epsilon,\bar{v}_\epsilon)\rightarrow (\bar{u}_i,\bar{v}_i)$ in $C^2_{loc}(\mathbb{R}^2)\times C^2_{loc}(\mathbb{R}^2)$. Moreover, $(\bar{u}_i,\bar{v}_i)$ solves
system (\ref{bu2}).
Choose constant $a>0$ such that $$b:=\inf_{\mathbb{R}^2}V-f(a)/a+f(-a)/a-g(a)/a+g(-a)/a>0.$$
By Lemma \ref{lbu}, we know $\bar{u}_i,\bar{v}_i$ decay exponentially to $0$ as $|x|\rightarrow\infty$. Then there exist $R_0, \epsilon_0>0$ such that for $0<\epsilon<\epsilon_0$, $|\bar{u}_\epsilon(x)|+|\bar{v}_\epsilon(x)|\leq a$, for all $|x|=R_0.$
Together with Lemma \ref{l5.27} we infer
$|\bar{u}_{\epsilon}(x)|+|\bar{v}_{\epsilon}(x)|\leq a$,  for all $x\in \cup_{j\neq i}[\partial(\Lambda''_j)^\epsilon-z^\epsilon_i]\cup \partial B_{R_0}(0).$
Recall that $|u_\epsilon|+|v_\epsilon|$ does not admit any local maximum
over $\mathbb{R}^2\backslash{(\Lambda'')^\epsilon}$  for sufficiently small $\epsilon$.
Then
\begin{equation*}\label{bu4}|\bar{u}_{\epsilon}(x)|+|\bar{v}_{\epsilon}(x)|\leq a, \quad\forall x\in \Omega:=\cap_{j\neq i}\bigl((\Lambda''_j)^\epsilon-z^\epsilon_i\bigr)^c\cap  B^c_{R_0}(0).\end{equation*}
By (\ref{bu7}) we obtain
$$-\Delta(|\bar{u}_\epsilon|+|\bar{v}_\epsilon|)+b(|\bar{u}_\epsilon|+|\bar{v}_\epsilon|)\leq0\quad\text{in} \ \Omega,\ \text{and}\ 0\leq |\bar{u}_\epsilon|+|\bar{v}_\epsilon|\leq a, \quad\text{in}\ \partial \Omega.$$
Below taking similar arguments of \cite[Lemma 5.32]{HUGO-PHD}, we get the desired results.\ \ \ \ $\Box$

{\bf Proof of Theorem \ref{t1.4}}: From Lemmas \ref{l5.22}, \ref{l5.23}, \ref{l5.27}, \ref{l5.29}, \ref{l5.14}, the conclusions follow directly.\ \ \ \ $\Box$

\subsection{Proof of Theorem \ref{t1.5}}

In this subsection, we consider positivity of solutions and uniqueness of local maximum points of solutions in $\Lambda_i$ under additional condition (H$_5$).

\begin{lemma}\label{l5.31} Under the assumptions of Theorem \ref{t1.5}, there holds $u_\epsilon>0$ and $v_\epsilon>0$. Moreover, there exists $\epsilon_0>0$ such that for every $0<\epsilon<\epsilon_0$,  the local maximum points of $u_\epsilon$, $v_\epsilon$ and  $u_\epsilon+v_\epsilon$ in $\Lambda^\epsilon_i$  are unique for each $i=1,\cdots,k$. Moreover, given $\tilde{x}_{i,\epsilon}$ a local maximum of $u_\epsilon$ in $\Lambda^\epsilon_i$ and $\tilde{y}_{i,\epsilon}$ a local maximum of $v_\epsilon$ in $\Lambda^\epsilon_i$, then $\tilde{x}_{i,\epsilon}-\tilde{y}_{i,\epsilon}\rightarrow0$ as $\epsilon\rightarrow0$.
\end{lemma}
{\bf Proof}: By (H$_5$), it is easy to see that $u_\epsilon>0$ and $v_\epsilon>0$. Let $x_{i,\epsilon}, y_{i,\epsilon}\in \Lambda^\epsilon_i$ be local maximum points of $u_\epsilon+v_\epsilon$. By Lemma \ref{l5.27}-(1) and Lemma \ref{l5.27}-(3), for some $r>0$ we have
$u_\epsilon(x_{i,\epsilon})+v_\epsilon(x_{i,\epsilon})\geq r$, $u_\epsilon(y_{i,\epsilon})+v_\epsilon(y_{i,\epsilon})\geq r$
and $V(\epsilon x_{i,\epsilon}), V(\epsilon y_{i,\epsilon})\rightarrow \inf_{\Lambda_i}V=V_i$. Moreover, we assume $\epsilon x_{i,\epsilon}\rightarrow \bar{x}\in \Lambda_i, \epsilon y_{i,\epsilon}\rightarrow \bar{y}\in \Lambda_i$. Define
$(\bar{u}_\epsilon,\bar{v}_\epsilon)(\cdot)=({u}_\epsilon,{v}_\epsilon)(\cdot+x_{i,\epsilon})$, and $(\tilde{{u}}_\epsilon,\tilde{{v}}_\epsilon)(\cdot)=({u}_\epsilon,{v}_\epsilon)
(\cdot+y_{i,\epsilon}).$
Then as in the proof of Lemma \ref{l5.23} we deduce $(\bar{u}_\epsilon,\bar{v}_\epsilon)\rightharpoonup(\bar{u},\bar{v})\neq(0,0)$ in $E$, $(\tilde{u}_\epsilon,\tilde{v}_\epsilon)\rightharpoonup(\tilde{u},\tilde{v})
\neq(0,0)$ in $E$, $(\bar{u}_\epsilon,\bar{v}_\epsilon)\rightarrow(\bar{u},\bar{v})$, and $(\tilde{u}_\epsilon,\tilde{v}_\epsilon)\rightarrow(\tilde{u},\tilde{v})$ in $C^1_{loc}(\mathbb{R}^2)\times C^1_{loc}(\mathbb{R}^2)$,
and $(\bar{u},\bar{v})$, $(\tilde{u},\tilde{v})$ solves (\ref{bu2}).

Take $z_{i,\epsilon}=y_{i,\epsilon}-x_{{i,\epsilon}}$. We claim that $z_{i,\epsilon}$ is bounded. Argue by contradiction, we assume $|z_{i,\epsilon}|\rightarrow+\infty$. Then for any $R>0$, $B_R(x_{{i,\epsilon}})$ and $B_R(y_{{i,\epsilon}})$ are disjoint for small enough $\epsilon$. Hence
\begin{equation}\label{5.63.0}\aligned
\Psi^i_{\epsilon}(u_\epsilon,v_\epsilon)&=\int_{B_R(0)}\kappa_\epsilon
+\int_{B_R(0)}\tilde{\kappa}_\epsilon+
\int_{\tilde{\Lambda}^\epsilon_i\backslash(B_R(0)\cup B_R(z_{i,\epsilon}))}\kappa_\epsilon,\endaligned\end{equation}
where $\Psi^i_{\epsilon}$ is given in (\ref{5.53.0}) and
$$\aligned
\kappa_\epsilon&=\nabla {\bar{u}_\epsilon}\nabla {\bar{v}_\epsilon}+V(\epsilon x+\epsilon x_{i,\epsilon}){\bar{u}_\epsilon}{\bar{v}_\epsilon}-F(\bar{u}_\epsilon)-
G(\bar{v}_\epsilon),\\
\tilde{\kappa}_\epsilon&=\nabla {\tilde{u}_\epsilon}\nabla {\tilde{v}_\epsilon}+V(\epsilon x+\epsilon y_{i,\epsilon}){\tilde{u}_\epsilon}{\tilde{v}_\epsilon}-F(\tilde{u}_\epsilon)-
G(\tilde{v}_\epsilon).
\endaligned$$
For each $R>0$ fixed, note that $(\bar{u}_\epsilon,\bar{v}_\epsilon)\rightarrow(\bar{u},\bar{v})$ in $C^1_{loc}(\mathbb{R}^2)\times C^1_{loc}(\mathbb{R}^2)$, we have
\begin{equation*}\label{5.60}\aligned
\lim_{\epsilon\rightarrow0}
\int_{B_R(0)}\kappa_\epsilon&={\Phi}_{V(\bar{x})}(\bar{u},\bar{v})
-\int_{B^c_R(0)}[\nabla \bar{u}\nabla \bar{v}+V(\bar{x})\bar{u}\bar{v}-{F}(\bar{u})-{G}(\bar{v})].
\endaligned\end{equation*}
Then for any given $\varrho>0$, there exists  large $R>0$ such that
$\lim_{\epsilon\rightarrow0}
\int_{B_R(0)}\kappa_\epsilon\geq {\Phi}_{V(\bar{x})}(\bar{u},\bar{v})-\varrho\geq c_i-\varrho.$
Similarly
$\lim_{\epsilon\rightarrow0}
\int_{B_R(0)}\tilde{\kappa}_\epsilon\geq {\Phi}_{V(\bar{y})}(\tilde{u},\tilde{v})\geq c_i-\varrho.$

On the other hand, take a cut-off function $0\leq\phi_R\leq1$  such that $|\nabla \phi_R|\leq C/R$, $\phi_R=1$ in $\mathbb{R}^2\backslash ({B_{R}(0)\cup B_R(z_{i,\epsilon})})$, $\phi_R=0$ in $B_{R/2}(0)\cup B_{R/2}(z_{i,\epsilon})$. Let $K_i$ be a set such that $\tilde{\Lambda}_i\Subset K_i$ and $K_i\cap\tilde{\Lambda}_j=\emptyset$ for $j\neq i$. Take $\xi_i\in[0,1]$ such that $\xi_i=1$ in $\tilde{\Lambda}_i$ and $\xi_i=0$ in $\mathbb{R}^2\backslash{K_i}$. Denote $\xi^\epsilon_i(x)=\xi_i(\epsilon x+\epsilon x_{i,\epsilon})$.
Recall that
\begin{equation*}
-\Delta \bar{u}_\epsilon+V(\epsilon x+\epsilon x_{i,\epsilon})\bar{u}_\epsilon={g}(\bar{v}_\epsilon),\
-\Delta \bar{v}_\epsilon+V(\epsilon x+\epsilon x_{i,\epsilon})\bar{v}_\epsilon={f}(\bar{u}_\epsilon)\quad  \text{in}\ \mathbb{R}^2.
\end{equation*}
Testing the above system with $(\bar{v}_\epsilon\xi^\epsilon_i\phi_R,\bar{u}_\epsilon\xi^\epsilon_i\phi_R)$
and letting $\tilde{\Lambda}^\epsilon_{{i,\epsilon}}:=
\tilde{\Lambda}^\epsilon_{x_{i,\epsilon}}$, ${K}^\epsilon_{{i,\epsilon}}:=
{K}^\epsilon_{x_{i,\epsilon}}$ and
$$\Omega_1:={\tilde{\Lambda}^\epsilon_{i,\epsilon}}\backslash (B_R(0)\cup B_R(z_{i,\epsilon})),\ \Omega_2:=K^\epsilon_{i,\epsilon}\backslash{{
\tilde{\Lambda}^\epsilon_{i,\epsilon}}}\cup {(B_R(0)\backslash B_{\frac R2}(0))}\cup {(B_R(z_{i,\epsilon})\backslash B_{\frac R2}(z_{i,\epsilon}))},$$
 for simplicity,
we deduce
$$\aligned
2\int_{\Omega_1}\kappa_\epsilon
\geq&\int_{\Omega_1}\bigl[2\nabla \bar{u}_\epsilon\nabla \bar{v}_\epsilon\xi^\epsilon_i\phi_R+2V(\epsilon x+x^\epsilon_i)\bar{u}_\epsilon\bar{v}_\epsilon\xi^\epsilon_i\phi_R-{f}(\bar{v}_\epsilon)\bar{v}_\epsilon
\xi^\epsilon_i\phi_R-g(\bar{u}_\epsilon)\bar{u}_\epsilon\xi^\epsilon_i\phi_R\bigr]\\
=&\int_{\Omega_1}\bigl[\nabla \bar{u}_\epsilon\nabla (\bar{v}_\epsilon\xi^\epsilon_i\phi_R)+\nabla \bar{v}_\epsilon\nabla (\bar{u}_\epsilon\xi^\epsilon_i\phi_R)+2V(\epsilon x+\epsilon x^\epsilon_i)\bar{u}_\epsilon\bar{v}_\epsilon\xi^\epsilon_i\phi_R\\
&-{f}(\bar{v}_\epsilon)\bar{v}_\epsilon
\xi^\epsilon_i\phi_R-{g}(\bar{u}_\epsilon)\bar{u}_\epsilon\xi^\epsilon_i
\phi_R\bigr]+o_\epsilon(1)+o_R(1)\\
=&-\int_{\Omega_2}\bigl[\nabla \bar{u}_\epsilon\nabla (\bar{v}_\epsilon\xi^\epsilon_i\phi_R)+\nabla \bar{v}_\epsilon\nabla (\bar{u}_\epsilon\xi^\epsilon_i\phi_R)+2V(\epsilon x+\epsilon x^\epsilon_i)\bar{u}_\epsilon\bar{v}_\epsilon\xi^\epsilon_i\phi_R
\\
&-{f}(\bar{v}_\epsilon)\bar{v}_\epsilon\xi^\epsilon_i\phi_R
-{g}(\bar{u}_\epsilon)
\bar{u}_\epsilon\xi^\epsilon_i\phi_R\bigr]+o_\epsilon(1)+o_R(1).
\endaligned$$
Note that $(\bar{u}_\epsilon,\bar{v}_\epsilon)\rightarrow(\bar{u},\bar{v})$ and $(\tilde{{u}}_\epsilon,\tilde{{v}}_\epsilon)\rightarrow(\tilde{{u}},\tilde{{v}})$ in $C^1_{loc}(\mathbb{R}^2)\times C^1_{loc}(\mathbb{R}^2)$, together with (\ref{5.53.2}), for the above $\varrho>0$ we get for large $R>0$
$\liminf_{\epsilon\rightarrow0}
\int_{\Omega_1}\kappa_\epsilon\geq-\varrho$.
Therefore, from (\ref{5.63.0}) and Lemma \ref{l5.23} we know
$$-3\varrho+2c_i\leq \Psi^i_{\epsilon}(u_\epsilon,v_\epsilon)=c_i+o_\epsilon(1),$$
for every $\varrho>0$. This is impossible.

Let $z_0$ be such that $z_{i,\epsilon}\rightarrow z_0$ as $\epsilon\rightarrow0$. Then
$\bar{u}_\epsilon(x)+\bar{v}_\epsilon(x)
=\tilde{u}_\epsilon(x-z_{i,\epsilon})+\tilde{v}_\epsilon(x-z_{i,\epsilon})$, by the strongly convergence of $\bar{u}_\epsilon\rightarrow\bar{u}$ and  $\bar{v}_\epsilon\rightarrow\bar{v}$ in  $C^1_{loc}(\mathbb{R}^2)$ we have
$\bar{u}(x)+\bar{v}(x)=
\tilde{u}(x-z_0)+\tilde{v}(x-z_0)$. Since the maximum points of $\bar{u}+\bar{v}$ and $\tilde{u}+\tilde{v}$ are both $0$, we have $z_0=0$.

Observe that $$\nabla({\bar{u}_\epsilon}+{\bar{v}_\epsilon})(0)=
\nabla({\bar{u}_\epsilon}+{\bar{v}_\epsilon})
(y_{i,\epsilon}-x_{i,\epsilon})=0.$$
By Lemma \ref{lbu}, we know $\Delta(\bar{u}+\bar{v})(0)<0$, $\bar{u}$ is radial and $\bar{u}'(0)=0$, $\bar{u}''(r)<0$ for $r=|x|$ small. On the other hand, since $\bar{u}_\epsilon\in C^2$ and $\bar{u}_\epsilon\rightarrow \bar{u}$ in $C^2_{loc}(\mathbb{R}^2)$ as $\epsilon\rightarrow0$. It follows from \cite[Lemma 4.2]{NIWM} that $y_{i,\epsilon}=x_{i,\epsilon}$ for $\epsilon>0$ sufficiently small. Taking similar arguments as above we can conclude that $\tilde{x}_{i,\epsilon}-\tilde{y}_{i,\epsilon}\rightarrow0$ as $\epsilon\rightarrow0$.
\ \ \ \ $\Box$

{\bf Proof of Theorem \ref{t1.5}}: From Lemma \ref{l5.31}, the conclusions follow directly.\ \ \ \ $\Box$
\section{Proof of Theorem \ref{t1.6}}
\renewcommand{\theequation}{6.\arabic{equation}}
In this section, we shall prove Theorem \ref{t1.6} and we only show the first part of Theorem \ref{t1.6} since minor modifications are enough to give the second part of  Theorem \ref{t1.6}.

Given $n\in\mathbb{N}$, define the truncated function
\begin{equation*}\aligned f_n(s)=
\left\{ \begin{array}{lll}
f(s)\ & \text{for}\ -n\leq s\leq n\\
\frac{f'(n)}{e^n}e^s+f(n)-f'(n)\ & \text{for}\ s>n\\
\frac{f'(-n)}{e^{-n}}e^s+f(-n)-f'(-n)\ & \text{for}\ s<-n
\end{array}\right.
\endaligned
\end{equation*}
and define $g_n(s)$ in the same way. Moreover, denote $F_n(s):=\int^s_0 f_n(t)dt$ and $G_n(s):=\int^s_0 g_n(t)dt$.
It is easy to see that $f_n, g_n$ satisfy (H$_1$), (H$_3$)-(i), (H$_6$) and (H$_7$). Then for a fixed n, thanks to Theorem \ref{t1.4} there exists $\epsilon_{0,n}>0$ such that for $0<\epsilon<\epsilon_{0,n}$,
there are $(\varphi_\epsilon,\psi_\epsilon)$ satisfying the conclusions of Theorem \ref{t1.4} with the quantities appearing in
Theorem \ref{t1.4} depend of $n$. Then
$u_\epsilon:=\varphi_\epsilon(\epsilon x), v_\epsilon:=\psi_\epsilon(\epsilon x)$ solves
\begin{equation*}\label{8.3}-\Delta u+V(\epsilon x)u=g_n(v),\ -\Delta v+V(\epsilon x)v=f_n(u),\  \text{in}\ \mathbb{R}^2.\end{equation*}
 Moreover, by Lemma \ref{l5.23} we have
\begin{equation}\label{8.3.0}\mathcal{I}^n_\epsilon(u,v)
=\Sigma^k_{i=1}c_{i,n}+o_n(1),\quad\text{as}\ \epsilon\rightarrow0,\end{equation}
where
$$\mathcal{I}^n_\epsilon(u,v)=\int_{\mathbb{R}^2}\bigl[\nabla u\nabla v+V(\epsilon x)uv-F_n(u)-G_n(v)\bigr],\ \forall(u,v)\in E,$$
and $c_{i,n}$ is the energy of ground states of the problem
\begin{equation*}-\Delta u+V(x_i)u=g_n(v),\ -\Delta v+V(x_i)v=f_n(u), \  \text{in}\ \mathbb{R}^2.\end{equation*}

We also introduce the auxiliary function
\begin{equation*}\label{8.2.0}\aligned h_f(s)=
\left\{ \begin{array}{lll}
f(s)\ & \text{for}\ -1\leq s\leq 1\\
\frac{f'(1)}{1+\delta'}s^{1+\delta'}+f(1)-\frac{f'(1)}{1+\delta'}\ & \text{for}\ s>1\\
\frac{f'(-1)}{1+\delta'}|s|^{\delta'}s+f(-1)+\frac{f'(-1)}{1+\delta'}\ & \text{for}\ s<-1
\end{array}\right.
\endaligned
\end{equation*}
In view of (H$_6$), we have $f(s)/|s|^{1+\delta'}$ is increasing. Then  $f'(s)\geq(1+\delta')f(s)/s\geq(1+\delta')f(1)s^{\delta'}$ for $s\geq1$ and $f'(s)\geq(1+\delta')f(s)/s\geq(1+\delta')f(-1)|s|^{\delta'-1}s$. It is straight to check that
$$\lambda_1 h'_f\leq f'_n\ \text{for every}\ 0<\lambda_1<\min\{1,(1+\delta')f(1)/f'(1),-(1+\delta')f(-1)/f'(-1)\},$$
and so $\lambda_1 h_f\leq f_n$. In the same way as $h_f$, we can define a function $h_g$ and there exists $\lambda_2>0$ such that $\lambda_2 h'_g\leq g'_n$ and so $\lambda_2 h_g\leq g_n$.
Let $c_{\lambda h_f,\lambda h_g}$ be the energy of ground states of the problem
\begin{equation*}-\Delta u+V(x_i)u=\lambda h_g(v),\ -\Delta v+V(x_i)v=\lambda h_f(u),\ \text{in}\ \mathbb{R}^2.\end{equation*}Then taking similar arguments of Lemma \ref{l5.0}, we have $c_{i,n}\leq c_{\lambda h_f,\lambda h_g}$. Together with (\ref{8.3.0}), we have the following lemma.
\begin{lemma}\label{l8.1} For any $n\in\mathbb{N}$, there exists $\epsilon_{0,n}>0$ such that, for every $\epsilon\in(0,\epsilon_{0,n})$, we have $\mathcal{I}^n_\epsilon(u_\epsilon,v_\epsilon)\leq C_0$ for some $C_0>0$ independent of $n$ and $\epsilon$.
\end{lemma}

\begin{lemma}\label{l8.2}  Given $\rho>0$, $i\in\{1,\cdots,k\}$ and $n\in\mathbb{N}$, there exists $\epsilon_{0,n}>0$ such that, for $\epsilon\in(0,\epsilon_{0,n})$ we have $\varphi_\epsilon(x)$, $\psi_\epsilon(x)\leq1$ for all $x\in \mathbb{R}^2\setminus{\cup_{j\neq i}\Lambda_j}$ such that $|x-z_{i,\epsilon}|\geq\rho$, where $z_{i,\epsilon}$ and $(\varphi_\epsilon,\psi_\epsilon)$ are given in Theorem \ref{t1.4}.
\end{lemma}
{\bf Proof}: By Theorem \ref{t1.4} (2)-(iii) we have $|\varphi_\epsilon(x)|$, $|\psi_\epsilon(x)|\leq C_ne^{-\frac{c_n}{\epsilon}|x-z_{i,\epsilon}|}$ for all $x\in\mathbb{R}^2\backslash{\cup_{j\neq i}\Lambda_j}$. Choose $\epsilon_{0,n}\leq \rho c_n/\log C_n$, the conclusion follows directly.\ \ \ \ $\Box$

Consider $\rho>0$ be
such that $B_{2\rho}(z_{i,\epsilon})\subseteq \Lambda_i$ and a cut-off function $\xi_i$  such that
\begin{equation}\label{bu5}\xi_i=1\ \text{in} \ B_{\rho}(z_{i,\epsilon}),\
\xi_i=0\ \text{in} \ \mathbb{R}^2\backslash{B_{2\rho}(z_{i,\epsilon})},\end{equation} and denote $\xi_{i,\epsilon}(x)=\xi_i(\epsilon x+z_{i,\epsilon})$. We also consider the functions
$\bar{v}_\epsilon(x)=\varphi_\epsilon(\epsilon x+z_{i,\epsilon})$, $\bar{u}_\epsilon(x)=\psi_\epsilon(\epsilon x+z_{i,\epsilon}),$
which solve
\begin{equation}\label{8.6}
\aligned
\left\{ \begin{array}{lll}
-\Delta(\bar{u}_\epsilon \xi_{i,\epsilon})+V(\epsilon x+z_{i,\epsilon})\bar{u}_\epsilon \xi_{i,\epsilon}&=g_n(\bar{v}_\epsilon)\xi_{i,\epsilon}-
\bar{u}_\epsilon\Delta\xi_{i,\epsilon}-2\nabla\bar{u}_\epsilon\nabla \xi_{i,\epsilon},\\
-\Delta(\bar{v}_\epsilon \xi_{i,\epsilon})+V(\epsilon x+z_{i,\epsilon})\bar{v}_\epsilon \xi_{i,\epsilon}&=f_n(\bar{u}_\epsilon)\xi_{i,\epsilon}-
\bar{v}_\epsilon\Delta\xi_{i,\epsilon}-2\nabla\bar{v}_\epsilon\nabla \xi_{i,\epsilon}.
\end{array}\right.
\endaligned
\end{equation}

\begin{lemma}\label{l8.3}Given $n\in\mathbb{N}$,  there exists $\epsilon_{0,n}$ such that for $0<\epsilon<\epsilon_{0,n}$ we have
$$\|\bar{u}_\epsilon \xi_{i,\epsilon}\|+\|\bar{v}_\epsilon \xi_{i,\epsilon}\|\leq C,$$
where $C>0$ is independent of $n$ and $\epsilon$.
\end{lemma}
{\bf Proof}: Observe that$$\aligned \mathcal{I}^n_\epsilon(u_\epsilon,v_\epsilon)-\frac12\langle (\mathcal{I}^n_\epsilon)'(u_\epsilon,v_\epsilon),(v_\epsilon,u_\epsilon)\rangle
\geq
\bigl(\frac12-\frac{1}{1+\delta'}\bigr)\int_{\mathbb{R}^2}[f_n(u_\epsilon)
u_\epsilon+g_n(v_\epsilon)v_\epsilon].\endaligned$$
Then by Lemma \ref{l8.1} we obtain\begin{equation}\label{8.5} \int_{\mathbb{R}^2}
[f_n(u_\epsilon)u_\epsilon+g(v_\epsilon)v_\epsilon]\leq C_0,\end{equation}
where $C_0$ is independent of $n$ and $\epsilon$.

 Testing the second equation of (\ref{8.6}) with $\bar{v}_\epsilon \xi_{i,\epsilon}$ we have
\begin{equation}\label{8.8}\aligned\|\bar{v}_\epsilon \xi_{i,\epsilon}\|^2_{V_i}\leq&
\int_{\mathbb{R}^2}\bigl[f_n(\bar{u}_\epsilon)\xi_{i,\epsilon}
-\bar{v}_\epsilon\Delta\xi_{i,\epsilon}-2\nabla\bar{v}_\epsilon\nabla
\xi_{i,\epsilon}\bigr]\bar{v}_\epsilon \xi_{i,\epsilon}\\
&+\int_{\mathbb{R}^2}\bigl[V(z_{i,\epsilon})-V(\epsilon x+z_{i,\epsilon})\bigr](\bar{v}_\epsilon \xi_{i,\epsilon})^2,
\endaligned\end{equation}
where $V_i=\inf_{\Lambda_i}V$.
By (H$_1$), there exists small $\delta>0$ independent on $\epsilon$ and $n$, such that $|f_n(\bar{u}_\epsilon)|=|f(\bar{u}_\epsilon)|\leq \frac14|\bar{u}_\epsilon|$ when $|\bar{u}_\epsilon|\leq\delta$. Then
\begin{equation}\label{8.7}\aligned\Bigl|
\int_{\mathbb{R}^2}f_n(\bar{u}_\epsilon)\xi^2_{i,\epsilon}
\bar{v}_\epsilon\Bigr|\leq\frac14\int_{\{|\bar{u}_\epsilon|\leq\delta\}}
|\bar{v}_\epsilon||\bar{u}_\epsilon|\xi^2_{i,\epsilon}+\Bigl|
\int_{\{|\bar{u}_\epsilon|\geq\delta\}}f_n(\bar{u}_\epsilon)
\xi^2_{i,\epsilon}
\bar{v}_\epsilon\Bigr|.\endaligned\end{equation}
Using (H$_8$), we have $\lim_{n\rightarrow\infty}\frac{f'(n)}{e^n}=0$. Then
there exist $C', C''>0$ independent of $n$ and $\epsilon$ such that, for any $s\geq n$, $|f_n(s)|\leq C'e^s+f(n)\leq C'e^s+f(s)\leq C''e^s$. Similarly, for any $s\leq-n$, $|f_n(s)|\leq C'e^s+|f(-n)|\leq C''e^s$. Moreover, for $\delta<|s|\leq n$, we have $|f_n(s)|=|f(s)|\leq C''e^s$.
Then \begin{equation}\label{5.12.0}|f_n(s)|\leq C_0e^{|s|}\ \text{for some}\ C_0>0\ \text{independent of}\ n\ \text{and}\ s.\end{equation}
Let
$$\Lambda^1_{n,\epsilon}:=\bigl\{
x\in \Lambda^\epsilon_i: |\bar{u}_\epsilon|\geq\delta,\
|f_n(\bar{u}_\epsilon)\xi_{i,\epsilon}|
\geq C_0e^{\frac14}
\bigr\}.$$
$$\Lambda^2_{n,\epsilon}:=\bigl\{
x\in \Lambda^\epsilon_i: |\bar{u}_\epsilon|\geq\delta, |f_n(\bar{u}_\epsilon)\xi_{i,\epsilon}|
\leq C_0e^{\frac14}
\bigr\}.$$
Then  there exists $C_1>0$ independent of $n$ and $\epsilon$ such that $|f_n(\bar{u}_\epsilon)\xi_{i,\epsilon}|\leq C_1|\bar{u}_\epsilon|$ for any $x\in \Lambda^2_{n,\epsilon}$. Applying (\ref{bu1}) with $t=\frac{\bar{v}_\epsilon\xi_{i,\epsilon}}{
\|\bar{v}_\epsilon\xi_{i,\epsilon}\|_{V_i}}$ and $s=\frac{{f}_n(\bar{u}_\epsilon)\xi_{i,\epsilon}}{C_0}$, from Trudinger-Moser inequality and (\ref{5.12.0}) we conclude
$$\aligned\Bigl|\int_{\{|\bar{u}_\epsilon|\geq \delta\}}
\frac{f_n(\bar{u}_\epsilon)\bar{v}_\epsilon\xi^2_{i,\epsilon}}
{\|\bar{v}_\epsilon\xi_{i,\epsilon}\|_{V_i}}\Bigr|
\leq& C_0\int_{\Lambda^1_{n,\epsilon}}\frac{|{f}_n(\bar{u}_\epsilon)
\xi_{i,\epsilon}|}
{C_0}\Bigl[\log\bigl(\frac{|{f}_n(\bar{u}_\epsilon)\xi_{i,\epsilon}|}{C_0}\bigr)
\Bigr]^{\frac12}\\
&+\frac12
\int_{\Lambda^2_{n,\epsilon}}
\frac{1}{C_0}{f}^2_n(\bar{u}_\epsilon)\xi^2_{i,\epsilon}+C_0
\int_{\Lambda^\epsilon_i}\bigl[e^{\bigl(\frac{\bar{v}_\epsilon\xi_{i,\epsilon}}
{\|\bar{v}_\epsilon\xi_{i,\epsilon}\|}\bigr)^2}-1
\bigr]
\\ \leq& \bigl(\delta^{-\frac12}+\frac{C_1}{2C_0}\bigr)
\int_{\Lambda^\epsilon_i}f_n(\bar{u}_\epsilon)\bar{u}_\epsilon+C_2\leq C_3,\endaligned$$
where we have used (\ref{8.5}) in the last inequality  and $C_3$ is independent of $\epsilon$ and $n$.
Hence from (\ref{8.7}) we know
$$\Bigl|\int_{\mathbb{R}^2}f_n(\bar{u}_\epsilon)\xi^2_{i,\epsilon}
\bar{v}_\epsilon\Bigr|\leq\frac14(\|\bar{u}_\epsilon\xi_{i,\epsilon}\|^2_{V_i}+
\|\bar{v}_\epsilon\xi_{i,\epsilon}\|^2_{V_i})+C_3\|\bar{v}_\epsilon
\xi_{i,\epsilon}\|_{V_i}.$$
Below adapting similar arguments of \cite[Lemma 5.36]{HUGO-PHD}, from (\ref{8.8}) we get the desired result.\ \ \ \ $\Box$

\begin{lemma}\label{l8.4}
Given $n\in\mathbb{N}$, there exists $\epsilon_{0,n}$ such that for $0<\epsilon<\epsilon_{0,n}$ we have
$|u_\epsilon|_\infty+|v_\epsilon|_\infty\leq C$, with $C>0$ independent of $n$ and $\epsilon$.
\end{lemma}
{\bf Proof}: Change $\rho$ into $2\rho$ in the definition of the cut-off function $\xi_i$  given in (\ref{bu5}) and without loss of generality we assume $B_{4\rho}(z_{i,\epsilon})\subseteq \Lambda_i$. Moreover, consider a cut-off function $\tilde{\xi}_i$ such that $\tilde{\xi}_i=1$ in $B_{2\rho}(z_{i,\epsilon})$ and $\tilde{\xi}_i=0$ in  $\mathbb{R}^2\backslash{B_{4\rho}(z_{i,\epsilon})},$ and denote $\tilde{\xi}_{i,\epsilon}=\tilde{\xi}_i(\epsilon x+z_{i,\epsilon})$.
Similar to Lemma \ref{l8.3}, we have $\|\bar{u}_\epsilon\tilde{\xi}_{i,\epsilon}\|^2+
\|\bar{v}_\epsilon\tilde{\xi}_{i,\epsilon}\|^2\leq C'$, where $C'$ is independent on $\epsilon$ and $n$.
Moreover, it is easy to see that for any $\alpha>0$, there exists $C_\alpha>0$ (independent on $n$ and $s$) such that
 $$|g_n(s)|^2\leq\frac14|s|^2+C_\alpha(e^{\alpha s^2}-1),\ \forall s.$$
Letting $\alpha=\frac{2\pi}{C'}$, by Trudinger-Moser inequality we get
$$\aligned
\int_{\mathbb{R}^2}
|g_n(\bar{v}_\epsilon){\xi}_{i,\epsilon}|^{2}
&\leq\frac14
\int_{B_{2\rho/\epsilon}(0)}
|\bar{v}_\epsilon|^{2}|{\xi}_{i,\epsilon}|^{2}+
C\int_{B_{2\rho/\epsilon}(0)}
\bigl[e^{\frac{2\pi}{C'}\bar{v}^2_\epsilon}
-1\bigr]{\xi}^{2}_{i,\epsilon}
\\&\leq\frac14C'+
C
\int_{B_{4\rho/\epsilon}(0)}
\bigl[e^{\frac{2\pi}{C'}\bar{v}^2_\epsilon{\tilde{\xi}}^2_{i,\epsilon}}
-1\bigr]
{\tilde{\xi}}^{\frac32}_{i,\epsilon}\\
&\leq \frac14C'+
C\int_{\mathbb{R}^2}
\Bigl[e^{\frac{2\pi}{C'}
\|\bar{v}_\epsilon{\tilde{\xi}}_{i,\epsilon}\|^2
\bigl(
\frac{\bar{v}_\epsilon{\tilde{\xi}}_{i,\epsilon}}
{\|\bar{v}_\epsilon{\tilde{\xi}}_{i,\epsilon}\|}\bigr)^2}
-1\Bigr]\leq C.\endaligned$$
The constants $C$ here and below are independent on $\epsilon$ and $n$.
Hence,
$$\int_{\mathbb{R}^2}
|\bar{u}_\epsilon\Delta{\xi}_{i,\epsilon}|^{2}\leq
\Bigl(\int_{B_{2\rho/\epsilon}(0)}|\bar{u}_\epsilon|^4\Bigr)^{\frac12}
\epsilon^4
\Bigl(\int_{\mathbb{R}^2}|\Delta{\xi}_{i}(\epsilon x)|^4\Bigr)^{\frac12}\leq
C'\epsilon^{3}|\Delta{\xi}_{i}|^{2}_4\leq C,$$
and
$$\int_{\mathbb{R}^2}|\nabla\bar{u}_\epsilon\nabla{\xi}_{i,\epsilon}|
^{2}\leq C\epsilon^{2}
\int_{B_{2\rho/\epsilon}(0)\backslash{B_{\rho/\epsilon}(0)}}
|\nabla\bar{u}_\epsilon|^2\leq C.$$
Then $g_n(\bar{v}_\epsilon){\xi}_{i,\epsilon}-\bar{u}_\epsilon\Delta{\xi}_{i,\epsilon}-2\nabla\bar{u}_\epsilon\nabla{\xi}_{i,\epsilon}\in L^{2}(\mathbb{R}^2)$. Applying \cite[Lemma 2.4]{Zhang-chen-zou} to the first equation of system (\ref{8.6}) we deduce $|\bar{u}_\epsilon{\xi}_{i,\epsilon}|_\infty\leq C$, and similarly $|\bar{v}_\epsilon{\xi}_{i,\epsilon}|_\infty\leq C$. Combining with Lemma \ref{l8.2}, we get $|{u}_\epsilon|_\infty+|{v}_\epsilon|_\infty\leq C$ with $C$ independent of $n$ and $\epsilon$.\ \ \ \ $\Box$

{\bf Proof of Theorem \ref{t1.6}}: In view of Lemma \ref{l8.4}, we take $n=\max\{|{u}_\epsilon|_\infty, |{v}_\epsilon|_\infty\}$ in the definitions of truncated functions $f_n$ and $g_n$, then the results follow.\ \ \ $\Box$

\end{document}